\documentclass{report}
\usepackage[utf8]{inputenc}
\usepackage{titlesec}
\usepackage{amssymb}
\usepackage{amsmath}
\usepackage{amsthm}
\usepackage{bbm}
\usepackage{esvect}
\usepackage{relsize}
\usepackage{ marvosym }
\usepackage{graphicx}
\usepackage{lipsum}
\graphicspath{ {e:/} }

\newtheorem{theorem}{Theorem}[section]

\newtheorem{corollary}[theorem]{Corollary}
\newtheorem{lemma}[theorem]{Lemma}

\theoremstyle{definition}
\newtheorem{definition}[theorem]{Definition}

\theoremstyle{remark}
\newtheorem{remark}[theorem]{\bf {Remark}}

\newtheorem{example}[theorem]{\bf {Example}}

\theoremstyle{definition}
\newtheorem*{dfn}{Definition}

\theoremstyle{theorem}
\newtheorem*{thm}{Theorem}
\newtheorem*{lm}{Lemma}
\newtheorem*{crl}{Corollary}

\newcommand{\mychapter}[2]{

\setcounter{chapter}{#1}
\setcounter{section}{0}
\chapter*{#2}
\addcontentsline{toc}{chapter}{#2}
}

\newcommand{\mysection}[2]{

\setcounter{section}{#1}
\setcounter{subsection}{0}
\subsection*{#2}

\addcontentsline{toc}{section}{\protect\numberline{}#2}

}

\usepackage{color}
\usepackage[top=2.7cm, bottom=2.3cm, left=2.5cm, right=2.5cm]{geometry}

\newcommand{\Nset}{\mathbb{N}}

\newcommand{\Rset}{\mathbb{R}}

\titleformat
{\chapter} 
[display] 
{\bfseries\huge\itshape} 
{Part \ \thechapter} 
{0.5ex} 
{
    \rule{\textwidth}{1pt}
    \vspace{1ex}
    \centering
} 
[
\vspace{-0.5ex}%
\rule{\textwidth}{1pt}
] 

\titleformat
{\section}
{\bfseries\LARGE}
{\thesection.}{0.5em}{}
[
\vspace{1.5em}
]

\titleformat
{\subsection}
{\bfseries\Large}
{\thesubsection.}{0.5em}{}
[
\vspace{1ex}%
]

\begin{document}
\Large

\pagenumbering{gobble}
\clearpage

\begin{center}

\rule{\textwidth}{2pt}
    \vspace{1ex}

\Huge {\itshape \bfseries Bourbaki-complete spaces }
\medskip

{\itshape \bfseries and  Samuel realcompactification}

\bigskip

\huge

{\it Espacios Bourbaki-completos}
\medskip

{\it y realcompactificaci\'on de Samuel}

\vspace{1ex}%
\rule{\textwidth}{2pt}

\vspace{4ex}

\large

\vspace{1.5ex}

\vspace{3ex}
\LARGE


\vspace{4ex}

\Large A thesis submitted in fulfillment of the requirements 

for the degree of Doctor of Philosophy by

\medskip

\huge {\bf Ana S. Mero\~no Moreno}
\medskip

\Large under the supervision of 

\medskip

\huge {\bf Prof. M. Isabel Garrido Carballo}

\vspace{2ex}

\Large This PhD thesis was defended at the Faculty of Mathematics

of the Universidad Complutense de Madrid on April 10th, 2019.
\end{center}

\Large

\newpage
\mbox{}
\newpage

{\bf Acknowledgments}
\medskip

First of all, I would like to thank my parents and my husband because of their moral and financial support along all these years of studies. Besides, I am very grateful to Mª Isabel Garrido for all the mathematical discussions that we have had and for understanding and guiding my thoughts and reasonings during all these years we have worked together. Also, I would like to thank her for the ungrateful work of cleaning up all the messy manuscripts I have written during my research. Finally, I do not forget Aarno Hohti and Heikki Junnila, whose brief teachings have turned out to be exceptionally useful almost everywhere in these thesis.
\newpage
\mbox{}
\newpage







\tableofcontents{}
\newpage
\mbox{}
\thispagestyle{empty}

\mbox{}
\thispagestyle{empty}
\newpage
\mysection{-2}{Abstract}
\pagenumbering{roman}
\setcounter{page}{1}

\hspace{15pt} This thesis belongs to the area of General Topology and, in particular, to the field of study of uniform spaces. It is divided in three parts where the related topics {\it Bourbaki-completeness} and {\it Samuel realcompactification} are studied. Many of the results presented here have already been published by the author in the following papers.

\smallskip

\begin{enumerate}

\item[$\bullet$] [GaMe14] M. I. Garrido and A. S. Mero\~{n}o, \textit{New types of completeness in metric spaces},  Ann. Acad. Sci. Fenn. Math. 39 (2014) 733-758.
\smallskip

\item[$\bullet$] [GaMe16] M. I. Garrido and A. S. Mero\~{n}o, {\it On paracompactness, completeness and boundedness in uniform spaces}, Topology Appl. 203 (2016) 98-107. 
\smallskip

\item[$\bullet$] [GaMe17] M. I. Garrido and A. S. Mero\~no, {\it The Samuel realcompactification of a metric space}, J. Math. Anal. Appl. 456 (2017) 1013-1039.

\item[$\bullet$] [GaMe18] M. I. Garrido and A. S. Mero\~no, {\it The Samuel realcompactification}, Topology Appl. 241 (2018) 150-161.

\item[$\bullet$] [HoJuMe19] A. Hohti, H. Junnila and A. S Mero\~{n}o, \textit{On strongly \v{C}ech-complete spaces}, to appear in Topology Appl. (2019).

\end{enumerate}

\noindent In addition, new results, originated during the writing process, have been included. 
\smallskip

In the first part of the thesis we present many results of [GaMe14], [GaMe16] and some of [HoJuMe19]. More precisely, we study {\it Bourbaki-completeness} and {\it cofinal Bourbaki-completeness} (equivalently {\it uniform strong-paracompact- ness}) of uniform spaces. In particular,  we solve several primary problems related to  products, subspaces, hyperspaces, metric spaces and fine spaces.

Observe that these new concepts are respectively an extension of {\it completeness} and {\it cofinal completeness} (equivalently, {\it uniform paracompactness}). Moreover, we show that completeness and cofinal completeness are, respectively, equivalent to Bourbaki-completeness and cofinal Bourbaki-completeness whenever we consider uniformities having a base of star-finite covers (observe that this fact was partially proved in  [GaMe18]). Thus, this kind of uniform spaces are crucial along the thesis. For instance, we show that every cofinally Bourbaki-complete uniform space satisfies that its uniformity has a star-finite base. 

In spite of the above equivalence, working with {\it Bourbaki-Cauchy filters}, that is, those filters inducing  the property of Bourbaki-completeness, instead of Cauchy filters of the star-finite modification of a uniformity, will be really useful, being particularly relevant in the second part of the thesis.

\smallskip

In the second part of the thesis some of the main results of  [HoJuMe19] are included. This work began while  the author was visiting the University of Helsinki during the summer of 2014, under the supervision of Heikki Junnila and Aarno Hohti. This paper is focused on studying a topological property called {\it strong-\v{C}ech-completeness}, and, as  consequences, characterizations of the Bourbaki-completely metrizable spaces and an universal space for Bourbaki-complete metric spaces are obtained.

However, instead of using the topological property of strong \v{C}ech-comple- teness, and since we need results that preserve certain uniform structure,  we have decided to start this second part with embeddings. More precisely, we embed the Bourbaki-complete metric spaces (and after, the Bourbaki-complete uniform spaces)  into an universal space, in such a way that the embedding is uniform in one way and its inverse preserves at least Bourbaki-Cauchy filters.  

Next, we characterized those metric spaces which are metrizable by a Bourbaki-complete metric. These are exactly the spaces which are {\it completely metrizable spaces} and {\it strongly metrizable} at the same time, as it is shown in [HoJuMe19]. Moreover, we prove also that this is equivalent to be metrizable by a complete metric satisfying that its metric uniformity has a star-finite base. Notice that, in general, in spite of the above result,  Bourbaki-complete spaces do not have a star-finite base for its uniformity, as we show in the first part of the thesis.

Finally, we study those spaces metrizable by a cofinally Bourbaki-complete metric. These are exactly the spaces which are {\it strongly paracompact} and metrizable by a cofinally complete. Observe that this result was already been published in [GaMe14]. It turns out that the spaces metrizable by a cofinally Bourbaki-complete metric are strictly stronger that those spaces metrizable by cofinally complete and a Bourbaki-complete metric at the same time.

\smallskip

The third part of the thesis is dedicated to the Samuel realcompactification of uniform spaces, which is related to the well-known Samuel compactification, and contains all the results from [GaMe18] which are the extension to uniform spaces of the results in [GaMe17]. 

The first problem presented here is the characterization of the {\it Samuel realcompact spaces}, that is, those uniform spaces which are complete with the weak uniformity induced by all the real-valued uniformly continuous functions on the space. It turns out that these are exactly the Bourbaki-complete uniform spaces having no uniform partition of Ulam-measurable cardinal. Observe that this result, which is a uniform extension of the well-known Kat\v{e}tov-Shirota theorem characterizing the realcompact spaces,  is the link between the different parts of the thesis. On the other hand, observe that, instead of proving this result like in the papers [GaMe17] or [GaMe18] we use the embeddings from the second part of the thesis.

Finally, as a natural step in order to generalize the above result, we study which kind of metric or uniform spaces satisfy that  the Samuel realcompactification and the Hewitt realcompactification are equivalent.  However, we only provided partial results to this problem. These are not published yet.

\newpage
\mysection{-1}{Resumen en castellano}

\hspace{15pt} Esta tesis pertenece al \'area  de la Topolog\'ia General y en particular, al campo de estudio de los espacios uniformes. Est\'a dividida en tres partes donde se estudian los temas relacionados de  {\it Bourbaki-completitud} y  {\it Realcompactificaci\'on de Samuel}. Muchos de los resultados presentados aqu\'i aparecen en las siguientes publicaciones del autor. 

\smallskip

\begin{enumerate}

\item[$\bullet$] [GaMe14] M. I. Garrido and A. S. Mero\~{n}o, \textit{New types of completeness in metric spaces},  Ann. Acad. Sci. Fenn. Math. 39 (2014) 733-758.
\smallskip

\item[$\bullet$] [GaMe16] M. I. Garrido and A. S. Mero\~{n}o, {\it On paracompactness, completeness and boundedness in uniform spaces}, Topology Appl. 203 (2016) 98-107. 
\smallskip

\item[$\bullet$] [GaMe17] M. I. Garrido and A. S. Mero\~no, {\it The Samuel realcompactification of a metric space}, J. Math. Anal. Appl. 456 (2017) 1013-1039.

\item[$\bullet$] [GaMe18] M. I. Garrido and A. S. Mero\~no, {\it The Samuel realcompactification}, Topology Appl. 241 (2018) 150-161.

\item[$\bullet$] [HoJuMe19] A. Hohti, H. Junnila and A. S Mero\~{n}o, \textit{On strongly \v{C}ech-complete spaces}, to appear in Topology Appl. (2019).
\end{enumerate}

\noindent Tambi\'en han sido incluidos resultados nuevos, surgidos durante el proceso de escritura de la tesis.

En la primera parte de esta memoria aparecen muchos resultados de los trabajos [GaMe14], [GaMe16] y algunos de [HoJuMe19]. M\'as precisamente, estudiamos  la {\it Bourbaki-completitud} y la {\it cofinal Bourbaki-completitud} (equivalentemente,  la {\it fuerte paracompacidad  uniforme}) en los espacios uniformes. En particular, resolvemos varios problemas b\'asicos relacionados con productos, subespacios, hiperespacios, espacios m\'etricos y espacios con la uniformidad fina. 

Se hace observar que estos conceptos nuevos son, respectivamente, extensiones de la {\it completitud} y la {\it cofinal completitud } (equivalentemente, la {\it  paracompacidad uniforme}). Adem\'as  probamos que la completitud y la cofinal completitud son respectivamente equivalentes a la Bourbaki-completitud y la cofinal Bourbaki-completitud  cuando consideramos uniformidades que tienen una base de recubriminetos estrella-finitos. Por lo tanto, este tipo de espacios uniformes es de suma importancia a lo largo de la tesis. Por ejemplo, probamos que todo espacio uniforme cofinalmente Bourbaki-completo tiene una base estrella-finita para su uniformidad.

A pesar de la equivalencia anterior, trabajar con {\it filtros  de Bourbaki-Cauchy}, es decir, aquellos filtros que inducen la propiedad de la Bourbaki-completitud, en vez de con los filtros de Cauchy de la modificaci\'on estrella-finita de la uniformidad, es realmente \'util, siendo particularmente relevante en la segunda parte de la tesis.
\smallskip

En la segunda parte de la tesis, se incluyen algunos de los resultados principales de [HoJuMe19]. Este trabajo empez\'o durante una visita del autor a la Universidad de Helsinki en el verano de 2014, bajo la supervisi\'on de Heikki Junnila y Aarno Hohti. Esta publicaci\'on se centra en el estudio de la propiedad topol\'ogica llamada {\it \v{C}ech-completitud fuerte} y, como consequencias, se obtienen caracterizaciones del los espacios Bourbaki-completamente metrizables y la existencia de un espacio universal para los espacios Bourbaki-completos.

Sin embargo, en vez de utilizar la propiedad topol\'ogica de la \v{C}ech-completi- tud, y ya que necesitamos resultados que preseven cierta estructura uniforme, hemos decido empezar esta segunda parte con ciertas inmersiones. M\'as precisamente, sumergimos los espacios m\'etricos Bourbaki-completos (y despu\'es los espacios uniformes Bourbaki-completos) en un espacio universal, de tal modo que la inmersi\'on es uniforme en un sentido y en el sentido inverso preserva al menos los filtros de Bourbaki-Cauchy.

A continuaci\'on, caracterizamos los espacios metrizables por una m\'etrica Bourbaki-completa. Estos son exactamente los espacios que son {\it completamente metrizables} y {\it fuertemente metrizables} a la vez, como se muestra en [HoJuMe19]. Adem\'as, probamos  que esto es equivalente a ser metrizable por una m\'etrica completa tal  que su uniformidad m\'etrica tiene una base estrella-finita. Se hace notar que, a pesar del resultado anterior, los espacios Bourbaki-completos no tienen  en general una base estrella-finita para su uniformidad, como mostramos en la primera parte de la tesis.

Finalmente, estudiamos los espacios metrizables por una m\'etrica cofinalmente Bourbaki-completa. Estos son exactamente los espacios que son {\it fuertemente paracompactos} y metrizables por una m\'etrica cofinalmente completa, a la vez. Obs\'ervese que este resultado est\'a ya publicado en [GaMe14]. Se deriva que los espacios metrizables por una m\'etrica cofinalmente Bourbaki-completa son estrictamente m\'as fuertes que los espacios metrizables por una m\'etrica cofinalmente completa y Bourbaki-completa a la vez.
\smallskip

La tercera parte de la tesis est\'a dedica a la realcompactificaci\'on de Samuel, que est\'a relacionada con la conocida compactificaci\'on de Samuel, y contiene todos los resultados de [GaMe18] que son a su vez una extensi\'on a los espacios uniformes de los resultados en [GaMe17].

El primer problema que se presenta aqu\'i es una caracterizaci\'on de los {\it espacios Samuel realcompactos}, es decir, aquellos espacios uniformes que son completos con la uniformidad d\'ebil inducida por todas las funciones reales uniformemente continuas sobre el espacio. Se deriva que estos son exactamente los espacios Bourbaki-completos que no poseen ninguna partici\'on uniforme de cardinal Ulam-medible. Obs\'ervese que este resultado, que es una extensi\'on uniforme del conocido teorema de Kat\v{e}tov-Shirota  que caracteriza los espacios realcompactos, es el nexo entre las diferentes partes de la tesis. Por otra parte, obs\'ervese que en vez de probar este resultdo como en las publicaciones [GaMe17] o [GaMe18], utilizamos las inmersiones de la segunda parte de la tesis.

Finalmente, damos un paso natural en la generalizaci\'on del resultado anterior y estudiamos qu\'e tipo de espacios m\'etricos o uniformes satisface que la realcompactificaci\'on de Samuel y la realcompactificaci\'on de Hewitt son equivalentes. Sin embargo, solamente proporcionamos resultados parciales. Estos no han sido publicados todav\'ia.

\newpage

\mbox{}
\thispagestyle{empty}
\newpage




\mychapter{0}{Introduction}


\hspace{15pt} In this thesis two main topics are studied: {\it Bourbaki-completeness} and {\it Samuel realcompactification.} The first one was born in \cite{merono.completeness} from  the primary problem {\it of characterizing those metric spaces such that the closure of every Bourbaki-bounded subset is compact}. The second one was first intended to be studied in the frame of metric spaces in order to prove its strict dependence of the bornological notion of {\it Bourbaki-boundedness} (see \cite[Section 7]{merono.Samuel.metric} and \cite{merono.bornology}).  

Many of the results given in this thesis come from previous publications by the author. These results are distributed in the three parts of this thesis along the following way. 

In the first part, we find the results from \cite{merono.completeness} and \cite{merono.nets} on Bourbaki-complete and cofinal Bourbaki-complete metric spaces and uniform spaces. However, in order to organize all these results we have chosen to start in the frame of uniform spaces and then move on to metric spaces as a particular case. Moreover, many proofs, as well as many results, are new since we have considered as a main character the {\it star-finite modification  of a uniformity},  which did not appear in our publications until \cite{merono.Samuel.uniform}.

In the second part, we have tried to condensate the main results of \cite{merono.hh}. But, instead of following the approach of said paper, which studies a kind of {\it \v{C}ech-complete property}, we have started with several results on embeddings of Bourbaki-complete spaces inspired by publication \cite{balogh}. On the other hand, we must notice that our embeddings also take under consideration the uniform side, since we will need it in the last part of the thesis. 

Finally, in the third part of the thesis we include the results from \cite{merono.Samuel.uniform} on the Samuel realcompactification of a uniform space. These are the extension of the results in \cite{merono.Samuel.metric} on the Samuel realcompactification of a metric space. The proofs given here depend on the embedding results of the second part of the thesis and hence, these are different from the proofs appearing in the papers.
\smallskip

Since we have reunified all the above works under a common approach, new unpublished results have arisen. These have been also included here.

Next, we expose the main concepts and results obtained, as well as their background.
 
\bigskip

\bigskip

\mysection{-5}{Bourbaki-complete uniform spaces}


\hspace{15pt} Bourbaki-bounded subsets were introduced in the frame of uniform spaces by Hejcman in \cite{hejcman}, in order to extend the notion of bounded subset of a metric space. We define these using the following notation.

Given a set  $X$, a cover $\mathcal{C}$ of $X$ and $A\subset X$, let us write:
\begin{itemize}

\item $St(A,\mathcal{C})= \bigcup \{C\in \mathcal{C}: C\cap A\neq \emptyset\}$

\item $St^{0}(A,\mathcal{C})=A$



\item $St^{m}(A,\mathcal{C})=St(St^{m-1}(A,\mathcal{C}),\mathcal{C})$, $m\in \Nset$


\end{itemize}

\begin{dfn}A subset $B$ of a uniform space $(X,\mu)$ is a {\it Bourbaki-bounded subset} of $(X,\mu)$ if for every uniform cover $\mathcal{U}\in \mu$ there exist  $m\in \Nset$ and finitely many $U_1,...,U_k\in \mathcal{U}$ such that $$B\subset \bigcup _{i=1}^{k}St^{m}(U_i,\mathcal{U}).$$ In particular $(X,\mu)$ is a {\it Bourbaki-bounded space} if it is a Bourbaki-bounded subset of itself.
\end{dfn}


Bourbaki-bounded subsets have their origin in locally convex topological vector spaces (LCTVS over the real or complex field). Recall that  von Neumann introduced  a notion of bounded subset for LCTVS in order to generalize the notion of  boundedness of the normed vector spaces, which is invariant by equivalence of norms  \cite[Remark 1.6]{hejcman}.  Precisely, the Bourbaki-bounded subsets of a LCTVS, endowed with the uniformity induced by its structure as a topological group, are exactly  the von Neumann's bounded subsets. Note that, whenever we are dealing with normed spaces, Bourbaki-bounded subsets coincide with the bounded subsets by the norm. 

In the same way that bounded subsets of normed spaces are invariant  under equivalence of norms, Bourbaki-bounded subsets are invariant under uniform equivalence, which is an advantage over the bounded subsets of metric spaces since they are not necessarily invariant for uniformly equivalent metrics. Indeed, think on the bounded subsets of $(\Rset , d_u)$, where $d_u$ is the usual euclidean metric on $\Rset$, and $(\Rset, d)$, where $d$ is the bounded metric defined by $d(x,y)=\inf\{1, d_u(x,y)\}$.
\smallskip

On the other hand, it is easily seen that {\it every totally bounded subset of a uniform space is always a Bourbaki-bounded subset}. Moreover, recall that it is a classical result that {\it a metric space is complete if and only if the closure of every totally bounded subset is compact}. Therefore, it is quite natural to try to characterize those metric spaces satisfying that the closure of every Bourbaki-bounded subset is compact. 

Since totally bounded subsets of metric spaces are characterized trough {\it Cauchy sequences} and Cauchy sequences lead to completeness of the metric spaces (by definition {\it a metric space is complete if every Cauchy sequence converges}),  our way to proceed in \cite{merono.completeness} was to give a sequence characterization of Bourbaki-bounded subsets of a metric space (Theorem \ref{metric.BB}), and then asking for the clustering of this family of sequences, called {\it Bourbaki-Cauchy sequences}. In this way we obtained the property of Bourbaki-completeness of metric spaces which is equivalent to the above stated question about the closure of Bourbaki-bounded sets (Theorem \ref{metric.BC}).

Next, it is also natural to extend the notion of Bourbaki-completeness to the frame of uniform spaces trough {\it Bourbaki-Cauchy nets}  as we do in \cite{merono.nets} or, equivalently, trough {\it Bourbaki-Cauchy filters} as in \cite{merono.hh} and \cite{merono.Samuel.uniform}. Recall that completeness of uniform spaces is defined by means of {\it Cauchy nets} (\cite{willard}) or {\it Cauchy filters} (\cite{bourbaki}). In this way, we obtain {\it Bourbaki-completeness} in uniform spaces.
\smallskip

\begin{dfn}  A filter $\mathcal{F}$  of a uniform space $(X,\mu)$ is {\it Bourbaki-Cauchy} in $X$ if for every uniform cover $\mathcal{U}\in \mu$ there is some $m\in \Nset$ and $U\in \mathcal{U}$ such that $$F\subset St^{m}(U,\mathcal{U}) \text{ for some }F\in \mathcal{F}.$$ 
\end{dfn}

\begin{dfn} A uniform space is {\it Bourbaki-complete} if every Bourbaki-Cauchy filter clusters.

\end{dfn}

It is clear that {\it every compact uniform space is Bourbaki-complete  and that every Bourbaki-complete uniform space is complete}. The main result relating all these properties is the following theorem (see Theorem \ref{Bourbaki-compact}).

\begin{thm} For a uniform space $(X,\mu)$ the following statements are equivalent:

\begin{enumerate}
\item $X$ is compact;

\item $(X,\mu)$ is complete and totally bounded;

\item $(X,\mu)$ is Bourbaki-complete and Bourbaki-bounded.
\end{enumerate}
\end{thm}

Observe that a Bourbaki-complete uniform space always satisfies that the closure of every Bourbaki-bounded subset is compact.  {\it A first example of complete space which is not Bourbaki-complete is any infinite-dimensional Banach space.} Indeed, the closed unit ball of such spaces is a Bourbaki-bounded subset (and subspace) that is not compact.

However, the reverse question, that is, if every uniform space satisfying that the closure of every Bourbaki-bounded subset is compact must be Bourbaki-complete, is false as it is expected (see Example \ref{BB=CB.notcomplete}), contrarily to what happens for metric spaces (Theorem \ref{metric.BC}). We must notice now that  LCTVS satisfying that the closure of every Bourbaki-bounded subset is compact has been also traditionally considered.  These spaces are called {\it semi-Montel}, and, whenever the LCVTS is {\it barrelled}, these are exactly the {\it Montel spaces} (see \cite{jarchow}). For instance, if we endow a LCVTS with the ``weak topology" (induced by a family of continuous functionals), then it is always a semi-Montel space.

As for the non-metrizable uniform spaces, the (semi-)Montel spaces are not necessarily complete as it is shown in the difficult example by K\"omura  \cite[Section 5]{komura}. However, by Theorem \ref{metric.BC} {\it every semi-Montel Frechet space (that is, metrizable LCTVS) space is Bourbaki-complete}.    In spite of this precedent in the frame of LCTVS, we must make clear now that we will not talk of LCTVS anymore. 
\smallskip

At first, in the papers \cite{merono.completeness} and \cite{merono.nets}, we treated Bourbaki-completeness as a property stronger than usual completeness ({\it every Cauchy filter is a Bourbaki-Cauchy filter}) but not necessarily ``linked" to completeness. It was not until \cite{merono.Samuel.uniform} that we realized that Bourbaki-completeness of a uniform space $(X,\mu)$ is equivalent to completeness of the {\it star-finite modification} of the uniform space $(X,s_f\mu)$ (Theorem \ref{star-finite2}). By a {\it modification} of a uniform space $(X,\mu)$ we mean a uniformity $\nu$ compatible with the topology of $X$ having as a base or subbase of the uniformity a subfamily of covers from the uniformity $\mu$.

\begin{dfn} A cover $\mathcal{C}$ of a set $X$ is {\it star-finite} if for every  $C\in \mathcal{C}$ there are at most finitely many $C'\in \mathcal{C}$ such that $C$ meets $C'$.
\end{dfn}

That the family of all the star-finite covers from $\mu$ is a base for an admissible uniformity is well-known (see \cite{isbellbook}). In addition, from an old result by Nj\aa stad in \cite{njastad}, it can be deduced that {\it a subset of a uniform space $(X,\mu)$ is Bourbaki-bounded if and only if  it is totally bounded in $(X, s_f\mu)$} (Theorem \ref{BB.sf}). Moreover, we also prove the following results.

\begin{lm} {\rm (Theorem \ref{star-finite1})} For a uniform space $(X,\mu)$ the following statements are equivalent:
\begin{enumerate}
\item every Cauchy filter of the uniform space $(X,s_f\mu)$ is a Bourbaki-Cauchy filter of $(X,\mu)$;

\item every Bourbaki-Cauchy ultrafilter of $(X,\mu)$ is a Cauchy ultrafilter of the uniform space $(X, s_f\mu)$.
\end{enumerate}
\end{lm}

\begin{thm} {\rm (Theorem \ref{star-finite2})}A uniform space $(X,\mu)$ is Bourbaki-complete if and only if $(X,s_f\mu)$ is complete.

\end{thm}

In spite of the above equivalence, Bourbaki-Cauchy filters of a uniform space $(X,\mu)$ are revealed as a useful tool, and this is because they are, in many situations, more tangible and easy to work with than Cauchy filters of  the star-finite modification $(X, s_f\mu)$. In order to support the tangibility and the strength of Bourbaki-completeness and the Bourbaki-Cauchy filters, we study the following interesting problem.

Consider the {\it point-finite modification} $(X, p_f\mu )$ of the uniform space $(X,\mu)$, that is, the admissible uniformity on $X$ having as base all the point-finite covers from $\mu$ (\cite{isbellbook}).

\begin{dfn} A cover $\mathcal{C}$ of a set $X$ is {\it point-finite} if for every $x\in X$ there exists at most finitely many $C\in \mathcal{C}$ such that $x\in C$.
\end{dfn}

It is difficult  to show that not every complete metric space $(X, d)$ satisfies that $(X, p_f \mu _d)$ is complete (where $\mu _d$ denotes the metric uniformity). It was proved by Pelant in \cite{pelant.complete} that the Banach space $(\ell _{\infty}(\omega _1), ||\cdot ||)$ of all the bounded real-valued functions $f: \omega _1 \rightarrow \Rset$ over a set having as a cardinal the first uncountable ordinal $\omega _1$, endowed with the norm of the supremum, is such a space. On the other hand, it is clear that $s_f\mu \leq p_f\mu$ so every Bourbaki-complete uniform space $(X,\mu)$ satisfies that $(X, p_f\mu)$ is complete but the reverse implication is not true. Indeed, any separable infinite-dimensional Banach space, as $\ell _2(\Nset)$, is a counterexample, as every separable uniform space $(X,\mu)$ has a point-finite base for its uniformity, that is $\mu=p_f\mu$ (see \cite{vidossich}).

By all the foregoing, we ask if we can proceed backwards as we have done with Bourbaki-completeness in order to prove that $(X, p_f \mu _d)$ is complete if and only if the closure of every totally bounded set in $(X, p_f \mu _d)$   is compact. However this does not work as we will see in Theorem \ref{tot.bounded.point.finite}: {\it a subset $B$ of a uniform space $(X,\mu)$ is totally bounded if and only if every point-finite uniform cover of $X$ contains a finite subcover of $B$.} Therefore, we cannot characterize the completeness of $(X, p_f \mu _d)$ trough its totally bounded subsets. This is different to what happens to the completeness of $(X, s_f \mu _d)$, since, as we have previously said, the totally bounded sets of  $(X, s_f \mu _d)$ are exactly the Bourbaki-bounded subsets of $(X,d)$ (see Theorem \ref{BB.sf}, Theorem \ref{metric.BC} and Theorem \ref{star-finite2}).


\medskip


In addition, we will study subspaces, products and hyperspaces of Bourbaki-complete uniform spaces (Subsection 1.2.4), Bourbaki-completeness of the fine uniform space $(X,{\tt u})$ (Subsection 1.2.3), {\it sequential Bourbaki-completeness} and Bourbaki-complete metric spaces (Subsection 1.3.1), as well as, we will give many examples and counterexamples regarding the behavior of the star-finite modification and the point-finite modification (Subsections 1.3.1 and 1.3.1). Other classical problems, as finding families of real-valued functions characterizing Bourbaki-completeness (see \cite{beer-between1}), and metrization by a Bourbaki-complete metric, are considered, respectively, in Part 3 and Part 2 of the thesis. We must notice here that different versions of the Cantor's theorem for Bourbaki-complete metric spaces can be found in \cite{merono.completeness}. We have not included here these results by practical reasons. Moreover, unexpected results relating Bourbaki-complete metric spaces and metric spaces satisfying that the family of all the real-valued uniformly continuous functions is a ring, are published in \cite{merono.beer}. We have not included these results either because they are somehow out of the topic.

\bigskip

\bigskip

\mysection{-4}{Cofinally Bourbaki-complete uniform spaces}


\hspace{15pt} In parallel, we study another completeness like-property called {\it cofinally Bourbaki-completeness}. The reason  to study this property was motivated by Beer's paper \cite{beer-between1} were a long-standing property  of metric and uniform spaces, called {\it cofinal completeness}, is studied in the frame of metric  spaces.

\begin{dfn} A filter $\mathcal{F}$ of a uniform space $(X,\mu)$ is a {\it  cofinally Cauchy filter} if for every $\mathcal{U}\in \mu$ there is some $U\in \mathcal{U}$ such that $F\cap U\neq \emptyset$, for every $F\in \mathcal{F}$ 
\end{dfn}

\begin{dfn} A uniform space $(X,\mu)$ is {\it cofinal complete} if every cofinal Cauchy filter clusters.
\end{dfn}

Cofinal completeness is a property stronger than usual completeness ({\it every Cauchy filter is cofinal Cauchy}) and was first implicitly considered by Corson in 1958 in \cite{corson}, in order to express the property of paracompactness for uniform spaces. Precisely, he proved that a {\it a Tychonoff space is paracompact if and only if the space $X$ endowed with the fine uniformity ${\tt u}$ is cofinally complete.} Moreover, he gave these cofinal Cauchy filters the name of {\it weakly Cauchy}.  On the other hand, the term ``cofinal completeness" appeared in 1971 in Howes's paper \cite{howes.completeness}, where it is defined, equivalently, by nets. 

Independently, in 1978, Rice (\cite{rice}) introduced the property of {\it uniform paracompactness} which is one of the possible extension of paracompactness to the frame of uniform spaces.

\begin{dfn} A uniform space $(X,\mu)$ is {\it uniformly paracompact} is every open cover $\mathcal{G}$ has an open refinement $\mathcal{A}$ which is {\it uniformly locally finite}, that is, there exists some $\mathcal{U}\in \mu$ such that every $U\in \mu$ meets at most finitely many $A\in \mathcal{A}$.
\end{dfn}

\noindent The reviewer of Rice’s paper (see \cite{smith.r}) observed the following important equivalence (see also \cite{howesbook}).

\begin{thm}  A uniform space $(X,\mu)$ is cofinally complete if and only if it is uniformly paracompact. 
\end{thm}

\noindent Then, it is clear that the previous result by Corson can be obtained as a corollary of the above result. 

The basic bibliography on cofinal completeness is  Howes's book \cite{howesbook} and  Beer's paper \cite{beer-between1}. Other possible uniform extensions of paracompactness which are weaker than Rice's definition can be found in \cite{musaev}. However, here we will only consider 
Rice's notion.
\medskip


\begin{dfn} A filter $\mathcal{F}$  of a uniform space $(X,\mu)$ is {\it cofinally Bourbaki-Cauchy} in $X$ if for every uniform cover $\mathcal{U}\in \mu$ there is some  $m\in \Nset$ and  $U\in \mathcal{U}$ such that $$F\cap St^{m}(U,\mathcal{U})\neq \emptyset \text{ for every }F\in \mathcal{F}.$$ 
\end{dfn}

\begin{dfn} A uniform space $(X,\mu)$ is {\it cofinally Bourbaki-complete} if every cofinally Bourbaki-Cauchy filter clusters.
\end{dfn}

Since {\it every Bourbaki-Cauchy filter is cofinally Bourbaki-Cauchy} and every {\it cofinally Cauchy filter is also cofinally Bourbaki-Cauchy} then {\it every cofinally Bourbaki-complete uniform space is Bourbaki-complete and cofinally complete} (see Subsection 2.2.2 in order to see that {\it cofinal completeness together with Bourbaki-completeness does not implies cofinal Bourbaki-completeness}.)

Cofinal completeness and cofinal Bourbaki-completeness have the characteristic of transforming local properties of uniform spaces  into uniform local properties. More precisely we have the next result.

\begin{thm} {\rm (Theorem \ref{unif.local.compact.th})} For a uniform space $(X,\mu)$ the following properties are equivalent:
\begin{enumerate}

\item $(X,\mu)$ is uniformly locally compact;

\item $(X,\mu)$ is cofinally complete and locally totally bounded;

\item $(X,\mu)$ is cofinally Bourbaki-complete and locally Bourbaki-bounded.

\end{enumerate}

\end{thm}

\noindent Observe that the equivalence between $1$ and $2$ was noticed in \cite{rice} and fully proved in \cite{fletcher}.

\smallskip

As cofinal completeness is equivalent to uniform paracompactness, as we have previously said, we wonder if cofinal Bourbaki-completeness is equivalent to some kind of uniform paracompactness property. The answer will be done by means of the so called {\it strong paracompactness}.

\begin{dfn} A Tychonoff space is {\it strongly paracompact} if every open cover has an open star-finite refinement.
\end{dfn}

It is clear that {\it every strongly paracompact space is paracompact}. To give an example of a paracompact space which is not strongly paracompact we can take into account the result of Morita \cite{morita.star-finite} that states that {\it a connected  paracompact space $X$ is strongly paracompact if and only if it is Lindel\"{o}f}.

The property of strong paracompactness was extended to a uniform property by Hohti \cite[6.1]{hohti-thesis}. He called it {\it uniform hypocompactness}.

\begin{dfn} A uniform space $(X,\mu)$ is {\it uniformly strongly paracompact} if every open cover $X$ has a uniformly star-finite open refinement, where a cover $\mathcal{A}$ is {\it uniformly star-finite} if there exists $\mathcal{U}\in \mu$ such that for every $A\in\mathcal{A}$, $St(A,\mathcal{U})$ meets at most finitely many $A'\in \mathcal{A}$.
\end{dfn}

It turns out that,  {\it uniform strong paracompactness of a uniform space is equivalent to cofinal Bourbaki-completeness} (Theorem \ref{uniform.strongly.paracompact}). This is mainly proved in \cite{merono.completeness} and \cite{merono.nets} for metric spaces and uniform spaces (using nets), respectively. However, here we prove it considering again the star-finite modification of a uniform space. This can be summarized by the following two results.

\begin{lm} {\rm (Theorem \ref{star-finite1})} Let $(X,\mu)$ be a uniform space. Then $\mathcal{F}$ is a cofinally Bourbaki-Cauchy filter of $(X,\mu)$ if and only if it is a cofinally Cauchy filter of $(X,s_f\mu)$.
\end{lm}

\begin{thm} {\rm (Theorem \ref{uniform.strongly.paracompact})} Let $(X,\mu)$ be a uniform space. The following statements are equivalent:
\begin{enumerate}
\item $(X,\mu)$ is cofinally Bourbaki-complete;

\item $(X,s_f\mu)$ is cofinally complete and $\mu=s_f\mu$;

\item $(X,s_f\mu)$ is uniformly paracompact and $\mu=s_f\mu$;

\item $(X,\mu)$ is uniformly strongly paracompact.
\end{enumerate}
\end{thm}

As corollary we have a result parallel to Corson's one. Namely, {\it a Tychonoff space $X$ is strongly paracompact if and only if $(X,{\tt u})$ is cofinally Bourbaki-complete} (Theoreom \ref{uniformly.strong.fine})

Moreover, the above theorem give us the additional information that every uniformly strongly paracompact uniform space has a base of star-finite uniform covers, that is, $\mu=s_f\mu$. This is not necessarily satisfied by Bourbaki-complete uniform spaces (see Example  \ref{point-finite}.)

\smallskip

Similarly to Bourbaki-completeness, we will study subspaces, products, hyperspaces of confinally Bourbaki-uniform spaces (Subsection 1.2.4), {\it sequential cofinal Bourbaki-completeness} and cofinally Bourbaki-complete metric spaces (Subsection 1.3.3). For additional characterizations of cofinally Bourbaki-complete metric spaces in terms of functionals see \cite[Theorem 28]{merono.completeness} or \cite{manisha2}.  In addition, we will also give several examples and counterexamples (Subsections 1.2.2 and 2.2.2). The problem of metrization by a cofinally Bourbaki-complete metric will be studied later, in Part 2. 

\bigskip

\bigskip



\mysection{-3}{Embeddings  of Bourbaki-complete spaces and metrization results}

\hspace{15pt} To find a suitable embedding $\varphi$ of any Bourbaki-complete space into some ``universal" Bourbaki-complete space is of great help. In fact, this is done in \cite{merono.hh} in order to characterize those spaces which are metrizable by a Bourbaki-complete metric. Precisely, it is proved that {\it a metrizable space is metrizable by a Bourbaki-complete metric if and only if it is homeomorphic to a closed subspace of $D^{\omega _0}\times \Rset ^{\omega _0}$ for some discrete space $D$.}

From this result we can deduce that not every {\it completely metrizable space}, that is, a space metrizable by a complete metric, is metrizable by a Bourbaki-complete metric. Indeed, recall the following class of spaces and its characterization.

\begin{dfn} A Tychonoff space is  {\it strongly metrizable} if it has a base of the topology which consists of the members from countably many star-finite open covers. 
\end{dfn}

\begin{thm} {\rm (see \cite[Proposition 3.23 and notes p. 110]{pears})} A Tychonoff space is strongly metrizable if and only if it is homeomorphic to a subspace of $D^{\omega _0}\times \Rset ^{\omega _0}$ for some discrete space $D$.
\end{thm}

By the above result it is clear that {\it every Bourbaki-complete metric space is strongly metrizable}. Moreover, there is a result by Wiscamb \cite{wiscamb} that states that {\it a connected Tychonoff space is strongly metrizable if and only if it is Lindel\"of}. Therefore, {\it not every metrizable space is strongly metrizable} and hence, {\it not every completely metrizable space is Bourbaki-completely metrizable}. Thus, complete metrizability and Bourbaki-complete metrizability are not equivalent topological properties.

The following result summarizes some of the characterizations of Bourbaki-completely metrizable spaces that can be found in Theorem \ref{summary} (see also \cite[Theorem 3.6, Theorem 4.6 and Corollary 4.7]{merono.hh}).

\begin{thm} For a metrizable space $X$ the following statements are equivalent:
\begin{enumerate}
\item $X$ is metrizable by a Bourbaki-complete metric;
\item $X$ is metrizable by complete metric such that the induced metric uniformity has a star-finite base;
\item $X$ is completely metrizable and strongly metrizable;
\item $X$ is homeomorphic to a closed subspace of $D^{\omega _0} \times \Rset ^{\omega _0}$ for some discrete space $D$;
\item $X$ is homeomorphic to a closed subspace of a countable product of locally compact metric spaces.
\end{enumerate}
\end{thm}

From the previous theorem we will provide a new characterization of the strongly metrizable spaces (see Theorem \ref{summary2} and \cite{balogh} for additional characterizations). 

\begin{thm} A metrizable space $X$ is strongly metrizable if and only if it is metrizable by a metric $d$ satisfying that every Bourbaki-bounded subset of $(X,d)$ is totally bounded. 
\end{thm}

The way to prove the above results in \cite{merono.hh} was trough a \v{C}ech-complete like property called {\it strong \v{C}ech-completeness}. This was motivated by  \cite[Theorem 23]{merono.completeness}. Moreover, in \cite{merono.hh}  Bourbaki-complete metrizable spaces are characterized as those spaces which are metrizable by a complete {\it star-finite metric}. By a star-finite metric we mean a metric satisfying that for every $\varepsilon >0$ the family of open balls $\{B_{d}(x,\varepsilon):x\in X\}$ is star-finite (see \cite{balogh}). In \cite[Theorem 4.3]{merono.hh} an easy contruction of such a metric, compatible with the euclidean topology, is given for the real line $\Rset$.

Here, instead of following \cite{merono.hh}, we proceed by defining an embedding $\varphi$ for any Bourbaki-complete metric space, with the additional property that $\varphi$ must be uniformly continuous. More precisely, the following theorem will be proved.

\begin{thm} {\rm (Theorem \ref{embedding.metric})}  Let $(X,d)$ be a Bourbaki-complete metric space. Then, there exists an embedding $$\varphi: (X,d)\rightarrow \Big((\prod_{n\in \Nset}\kappa _n) \times \Rset ^{\omega _0} ,  \pi\Big)$$ where each $\kappa _n$ is a cardinal endowed with the uniformly discrete metric, $\pi$ denotes the usual product metric, $\varphi$ is uniformly continuous and $\varphi (X)$ is a closed subspace of  $(\prod_{n\in \Nset}\kappa _n )\times \Rset ^{\omega _0}.$ 
\end{thm}

We require the uniform continuity in the above theorem since we will need it for several results of Part 3. On the other hand, it is certainly stronger than the topological results in \cite{merono.hh}. In fact, in order to prove it, we need first an embedding result for complete metric spaces having a star-finite base for the metric uniformity (Theorem \ref{embedding.star-finite1}). And then that there is an admissible metric $d'$ on $X$ such that $(X,d')$ is complete, the metric uniformity $\mu _{d'}$ has a star-finite base and the identity map $i: (X,d)\rightarrow (X, d')$ is uniformly continuous (see Theorem \ref{metrizable.star-finite}), even if not every Bourbaki-complete metric space $(X,d)$ has a star-finite base for its uniformity.
\smallskip

We will also similarly prove an embedding theorem for Bourbaki-complete uniform spaces (see, in addition, Theorem \ref{embedding.star-finite2}).

\begin{thm} {\rm (Theorem \ref{universal.uniform})} Let $(X,\mu)$ be Bourbaki-complete uniform space. Then there exist  an embedding $$\varphi: (X,\mu)\rightarrow \Big((\prod_{i\in I, n\in \Nset}\kappa _n ^i)\times \Rset ^{\alpha} ,\pi\Big)$$ where each $\kappa _n ^i$ is a cardinal endowed with the uniformly discrete metric, $\pi$ denotes the usual product uniformity, $\alpha \geq \omega _0$,  $\varphi$ is uniformly continuous and $\varphi (X)$ is a closed subspace of $(\prod_{i\in I, n\in \Nset}\kappa _n ^i)\times \Rset ^{\alpha}.$ 
\end{thm}

From the above result we can characterize the class of all the Tychonoff spaces that are uniformizable by a Bourbaki-complete uniformity (Theorem \ref{summary.uniform}). This coincides exactly with the $\delta$-{\it complete} spaces of Garc\'ia-M\'aynez \cite{garcia.delta-complete}, that is, those Tychonoff spaces $X$ satisfying that $(X, s_f {\tt u})$ is complete where $s_f{\tt u}$ denotes the star-finite modification of the fine uniformity ${\tt u}$.

\bigskip

\bigskip

\mysection{-2}{Metrization by a cofinally complete and Bourbaki-complete metric or by a cofinally Bourbaki-complete metric.}


\hspace{15pt} After having studied the metrization of a space by a Bourbaki-complete metric, we consider the problem of metrization by a cofinally Bourbaki-complete metric. This problem, of course is strictly related to the problem of metrization by a cofinally complete metric that was solved by Romaguera in \cite{romaguera}.

\begin{thm} A metrizable space $X$ is metrizable by a cofinally complete metric if and only if $nlc(X)$ is a compact subset, where $nlc(X)$ denotes the subset of $X$ of all the points of non local compactness, that is, those points no having a locally compact neighborhood in $X$.
\end{thm}

In \cite{merono.completeness} (see Theorem \ref{cof.B.complete.metrizable}, see also \cite[Theorem 34]{merono.completeness} for additional characterizations and \cite{heikki} for further considerations) we characterize those metric spaces that are metrizable by a cofinally Bourbaki-complete metric, applying the above result of Romaguera.

\begin{thm}A space $X$ is metrizable by a cofinally Bourbaki-complete metric if and only if it is strongly paracompact and $nlc(X)$ is compact.
\end{thm}

Related to this problem, we ask which spaces are those being metrizable by a metric which is Bourbaki-complete metric and cofinally complete at the same time. At first, one can think that it is the same than being metrizable by a metric which is cofinally Bourbaki-complete. However, it is topologically weaker. Indeed, we prove the following result (Theorem \ref{sigma.unif3}).

\begin{thm} A metrizable space $X$ is metrizable by a metric $X$ which is at the same time Bourbaki-complete and cofinally complete if and only if $X$ is strongly metrizable and $nlc(X)$ is compact.
\end{thm}

\noindent  Moreover, by the help of these results, we are capable of giving an example of a metric space cofinally complete and Bourbaki-complete at the same time which is not metrizable by a cofinally Bourbaki-complete metric. (Example \ref{ex.strongly.metrizable}). 

\smallskip

The above problem is also related to the uniform problem that we explain next.

Recall that, Zarelua in \cite[Lemma 5]{zarelua} proved that {\it a metrizable is strongly metrizable if and only if it is completely paracompact}, where completely paracompact is defined as follows.

\begin{dfn} A Tychonoff space $X$ is {\it completely paracompact} if every open cover $\mathcal{G}$  has an open refinement $\mathcal{V}$ which is a subcollection of a family of sets $\bigcup _{n\in \Nset} \mathcal{V}_n$ where each $\mathcal{V}_n$ is an open star-finite cover of $X$.
\end{dfn}

It is known that {\it every strongly paracompact space is completely paracompact}, that {\it every completely paracompact space is paracompact}, and that all these three topological properties are different (see, for instance, \cite{pears}). Therefore,  complete paracompactness is a property that lies between strong paracompactness and paracompactness, so it is natural to ask if there exists a uniform version of complete paracompactness lying between uniform strong paracompactness and uniform paracompactness, that satisfies in addition a Corson's like-theorem: {\it a Tychonoff space $X$ is completely paracompact if and only if $(X, {\tt u})$ is uniformly completely paracompact} (Theorem \ref{sigma.unif2}). Therefore, we propose the next definition of uniform complete paracompactness.

\begin{dfn} A uniform space $(X,\mu)$ is {\it uniformly completely paracompact} if it is cofinally complete (equivalently uniformly paracompact) and every open uniform cover has a $\sigma$-star-finite uniform open refinement, where by $\sigma$-star-finite open cover we mean a cover  $\mathcal{V}$ which is a subcollection of a family of sets $\bigcup _{n\in \Nset} \mathcal{V}_n$ where each $\mathcal{V}_n$ is an open star-finite cover of $X$.
\end{dfn}

Observe that the above definition is not stated  only in terms of covers, distinctly to  uniform paracompactness or uniform strong paracompactness. Therefore, it is open the question of giving such a kind of definition. However, we think that the above definition is not so bad, since a uniform space is uniformly strongly paracompact if and only if it is cofinally complete and every uniform open cover has a star-finite uniform (open)  refinement, equivalently the uniformity, has a star-finite base (Theorem \ref{uniform.strongly.paracompact}). The difference is that we do not know if the family of all the uniform $\sigma$-star-finite open covers from a uniform space is a base for some compatible uniformity on the space. However, we know that in the particular case of completely paracompact fine spaces (Theorem \ref{sigma.unif2}) or uniform spaces having a countable base for its uniformity (Theorem \ref{sigma.unif1}), the family of all the uniform $\sigma$-star-finite open covers from a uniform space is a base for a compatible uniformity on theses spaces.


It turns out that {\it a metrizable space $X$ is metrizable  by a metric whose uniformity is uniformly completely paracompact if and only if $X$ is strongly metrizable and $nlc(X)$ is compact, that is, if and only if $X$ is metrizable by a metric which is cofinally complete and Bourbaki-complete at the same time} (see Theorem \ref{sigma.unif3}). However, even if the two problems studied in this result are topologically equivalent, we know that they are not uniformly equivalent. Indeed, Example \ref{ex.strongly.metrizable} shows a uniformly completely paracompact metric space which is not Bourbaki-complete nor cofinally complete.

\bigskip

\bigskip



\mysection{-1}{Samuel realcompact spaces and Kat\v{e}tov-Shirota-type Theorem}

\hspace{15pt} One of the main results of this part of the thesis is the characterization of those uniform spaces $(X,\mu)$ which are complete when they are endowed with the weak uniformity $wU_{\mu}(X)$ induced by all the real-valued uniformly continuous functions $U_\mu(X)$ on $(X,\mu)$ (\cite{isbell.weak}). 

Recall the following definitions.

\begin{dfn} A Tychonoff space $X$ is {\it realcompact} if it is homeomorphic to a closed subspace of a product of real-lines $\Rset^{\alpha}$. A {\it realcompactification} of a Tychonoff space $X$ is a realcompact space $Y$ in which $X$ is densely embedded.
\end{dfn}

We will call the completion of $(X, wU_\mu(X))$, {\it the Samuel realcompactification} of $(X,\mu)$ (Theorem \ref{completion.realcompactification}). 

The name of Samuel comes by several resemblances with the Samuel compactification $s_\mu X$ of $(X,\mu)$. Recall that the Samuel compactification is defined as the completion of $(X, wU^{*}(X))$, where $wU^{*}(X)$ denotes the weak uniformity on $(X,\mu)$ induced by all the bounded real-valued uniformly continuous functions $U^{*}_{\mu}(X)$.

Precisely, as well as the {\it Samuel compactification} is the smallest realcompactification and compactification of $(X,\mu)$, in the natural order of realcompactifiactions (\cite{engelking.realcompact}), such that every bounded real-valeud uniformly continuous function on $(X,\mu)$ can be continuously extended to it, the Samuel realcompactification is the smallest realcompactification of $(X,\mu)$ such that every real-valued uniformly continuous function on $(X,\mu)$ can be continuously extended to it. 

Moreover, as well as every compactification of a Tychonoff space $X$ is the Samuel compatification for some compatible uniformity on $X$, the same works for realcompactifications, that is, every realcompactification of $X$ is the Samuel realcompactification for some compatible uniformity on $X$ (Theorem \ref{all.samuel}). For instance, the well-known Hewitt realcompactification $\upsilon X$ of a Tychonoff space $X$, induced by the family of real-valued continuous functions $C(X)$,  is the Samuel realcompactification of $(X,{\tt u})$, since for the fine unifomity ${\tt u}$, $C(X)=U_{\tt u}(X)$. 

We will denote by $H(U_{\mu}(X))$ the Samuel realcompactification because it is exactly the set of all the real unital lattice homomorphisms on $U_\mu(X)$. 

\begin{dfn} A uniform space $(X,\mu)$ is {\it Samuel realcompact}  if $(X, wU_{\mu}(X))$ is complete.
\end{dfn}

Our result characterizing Samuel realcompact spaces is very near to the classical Kat\v{e}tov-Shirota Theorem, which relates realcompactness to topological completeness. Recall that a Tychonoff space is {\it topologically complete} if and only if $(X,{\tt u})$ is complete.  

\begin{thm}{\rm ({\sc Kat\v{e}tov-Shirota Theorem}, \cite{shirota1}, \cite{katetov})} Let $X$ be a Tychonoff space, then $X$ is realcompact if and only it is topologically complete and there is no closed discrete subspace of Ulam-measurable cardinal.
\end{thm}

The importance of the Kat\v{e}tov-Shirota Theorem lies on the fact that the class of topological complete spaces satisfying that there is no closed discrete subspace having Ulam-measurable cardinal is ``wide enough". This ``wideness'' is given by the ``huge size'' of Ulam-measurable cardinals and by the fact that  the existence of Ulam-measurable cardinals is not provable in ZFC  (equivalently, assuming that ZFC is consistent, ZFC+`` there is no Ulam-measurable cardinal" is consistent). Indeed, there is an old result by Tarski and Ulam (\cite{ulam})  saying that the {\it least Ulam-measurable cardinal is strongly inaccessible}. This implies that  starting from $\omega _0$, any cardinal that we can obtain by usual arithmetic of cardinals, that is, by addition, multiplication, exponentiation, formation of suprema, and by the passage from one of these cardinals obtained to its immediate successor or to any smaller cardinal will be always smaller than the least Ulam-measurable cardinal (\cite{gillman}). 

\begin{dfn} A filter $\mathcal{F}$ on a set $S$ satisfies the {\it countable intersection property} if for every countable family $\{F_n:n\in \Nset\}\subset \mathcal{F}$, $$\bigcap_{n\in \Nset} F_n\neq \emptyset.$$
\end{dfn}


\begin{dfn} An infinite cardinal $\kappa$ is {\it Ulam-measurable} if in any set of cardinal $\kappa$ there is a free (non-principal) ultrafilter satisfying the countable intersection property. 
\end{dfn}

The above definition of Ulam-measurable cardinal is not the original one (\cite{gillman}). In fact, the concept of Ulam-measurable cardinal was motivated by the measure problem (as its name tell us) of knowing  it there exist certain non-trivial measures on any set (\cite{jech}). In particular, what is telling the above definition is that a {\it discrete space is realcompact if and only if it does not have Ulam-measurable cardinal}. This result was not evident at first and it was proved by Mackey in \cite{mackey} (see also \cite{gillman}).

\smallskip

\smallskip

The following result is our characterization of the Samuel realcompact spaces. This is parallel to the classical Kat\v{e}tov-Shirota result, but the roll of topological completeness is now played by Bourbaki-completeness and, instead of  closed discrete subspaces, we consider uniform partitions, that is, partitions of the space that are refined by some uniform cover.

\begin{thm} {\rm (Theorem \ref{measurable2}, {\sc Kat\v{e}tov-Shirota-type Theorem})} A uniform space $(X,\mu)$ is Samuel realcompact if and only if it is Bourbaki-complete and there is no uniform partition having Ulam-measurable cardinal. 
\end{thm}

\noindent The above theorem was originally proved in \cite{merono.Samuel.metric}, \cite{merono.Samuel.uniform} and \cite{husek}. However, the proof that we present here is different and it strongly depends on the embeddings of the previous part of the thesis.
\smallskip

Next, if we apply our Kat\v{e}tov-Shirota-type Theorem to the fine uniformity ${\tt u}$ on a space $X$ then we get a characterization of realcompactness in terms of the property of $\delta$-completeness, which is stronger than topological completeness.

\begin{crl} {\rm (Theorem \ref{garcia})} A Tychonoff space $X$ is realcompact if and only if $X$ is $\delta$-complete and there is no open partition of $X$ having Ulam-measurable cardinal.
\end{crl}

We must notice here that there is a characterization of Samuel realcompact uniform space due to Rice and Reynolds \cite{reynolds}. More precisely, their result is the following.

\begin{thm} {\rm (\cite[Corollary 2.5]{reynolds})} Let $(X,\mu)$ be a uniform space. Then the space $(X,wU_\mu(X))$ is complete if and only if $(X, s_f\mu)$ is complete and there is no uniform cover in $\mu$ of Ulam-measurable cardinal.
\end{thm}  

Since we already know that completeness of $(X,s_f\mu)$ is equivalent to Bour- baki-completeness of $(X,\mu)$, one can think that our result characterizing Samuel realcompactness can be derived from the above result of Rice and Reynolds. However, the proofs that we gave in \cite{merono.Samuel.metric} and \cite{merono.Samuel.uniform} strongly depend on Bourbaki-Cauchy filters, while the proof of Rice and Reynolds strongly depends on a previous result \cite[Theorem 1.1]{reynolds} regarding uniform covers with a point-finite uniform refinement. 

Moreover, we prove our Kat\v{e}tov-Shirota-type Theorem after having proved Theorem \ref{measurable1} (see next result) that characterizes exactly those uniform spaces having no uniform partition of Ulam-measurable cardinal. This family of spaces do not coincide in general with the family of uniform spaces having no uniform cover of Ulam-measurable cardinal. However, for uniform spaces having a star-finite base for the uniformity, it is equivalent to have no uniform partition of Ulam-measurable cardinal than to have no uniform cover of Ulam-measurable cardinal.

\begin{thm} {\rm (Theorem \ref{measurable1})} Let $(X,\mu)$ be a uniform space. The following statements are equivalent:

\begin{enumerate}
\item every uniform partition of $(X,\mu)$ has no Ulam-measurable cardinal;

\item the completion of $(X, s_f\mu)$ is realcompact;

\item every uniform cover of $(X, s_f\mu)$ ha no Ulam-measurable cardinal;

\item every uniform partition of $(X,s_f\mu)$ has no Ulam-measurable cardinal;

\end{enumerate}
\end{thm}

Finally, as an application of the previous result, we characterize the Samuel realcompactification  of a uniform space $(X,\mu)$, and the Samuel realcompact uniform spaces, using the modification $e(s_f\mu)$ and $s_f(e\mu)$ (see Theorem \ref{countable1} and Theorem \ref{countable2}).

\bigskip

\bigskip

\mysection{0}{Relating the Samuel and the Hewitt realcompactification}


\hspace{15pt} The results of the last section of the thesis are motivated by the following problem: {\it to characterize those uniform spaces such that their Samuel realcompactification and their Hewitt realcompactification are equivalent, independently of the cardinality of their uniform partitions}.
\smallskip

Therefore, the first problem that we study is (Theorem \ref{realcompactifcations->Bourbaki-complete}) if  every metric space $(X,d)$ satisfying that $H(U_d(X))$ and $\upsilon X$ are equivalent realcompactifications is Bourbaki-complete.  Observe that this question is clearly false is the frame of uniform spaces. Indeed, the ordinal space of all the countable ordinals $([0, \omega_1),{\tt u})$ is a counterexample. 

In  order to solve this problem we characterize Bourbaki-completeness of a metric space $(X,d)$ by means of the subalgebra of all the Cauchy-continuous functions $CC_{s_f \mu _d}(X)$ over $(X, s_f \mu_d)$. More precisely, we have the following result.

\begin{thm} {\rm  (Theorem \ref{Bourbaki-complete.functions})} A metric space $(X,d)$ is Bourbaki-complete if and only if $C(X)=CC_{s_f\mu _d}(X)$.
\end{thm}

Finally, from the above theorem we get the desired result.

\begin{thm} {\rm (Theorem \ref{realcompactifcations->Bourbaki-complete})} Let $(X,d)$ be a metric space. If  $\upsilon X$ and $H(U_d(X))$ are equivalent realcompactifications, then $(X,d)$ is Bourbaki-complete.
\end{thm}

The next problem that we consider, is the converse to the previous one, that is if every Bourbaki-complete metric space or Bourbaki-complete uniform space satisfies that its Samuel realcompactification and its Hewitt realcompactification are equivalent. One can argue that, since Ulam-measurable cardinals are ``wide enough", as we have previously explained, by the Kat\v{e}tov-Shirota Type Theorem \ref{measurable2}, it is enough to consider the Samuel realcompact spaces. Indeed, it is clear that Samuel realcompact spaces satisfies that both realcompactifications are equivalent as they are homeomorphic to the original space. However, we believe that it is not enough to stop at the level of non Ulam-measurable cardinals, mainly by the following fact: {\it every UC space (in particular, every uniformly discrete space), independently of its cardinal, satisfies that the Samuel realcompactification and the Hewitt realcompactifications are equivalent}. Indeed, recall that the UC spaces are exactly those uniform spaces satisfying that every continuous real-valued function is uniformly continuous function and the in the frame of metric spaces every UC space is Bourbaki-complete (see \cite{merono.completeness}). However, there are UC spaces which are not necessarily complete as, again, the ordinal space of all countable ordinals $([0,\omega _1),\tt{u})$.


By all the foregoing, in Theorem \ref{positive}, we have characterized a class of Bourbaki-complete metric spaces (non necessarily being UC spaces), having uniform partitions with Ulam-measurable cardinal  and satisfying that both realcompactifications are equivalent. The result is as follows.

\smallskip

Given $n\in \Nset$, we can choose representative points $x_i\in X$, $i\in I_n$ such that there is a partition $$\mathcal{P}_n=\{B_{d}^{\infty}(x_i,1/n):i\in I_n\}$$ of $X$ where for each $i\in I_n$, $B_{d}^{\infty}(x_i,1/n)=\bigcup_{m\in \Nset}St^{m}(x_i, \mathcal{B}_{1/n})$ and $\mathcal{B}_{1/n}=\{B_{d}(x,1/n):x\in X\}$. Then $\mathcal{P}_n$ is a uniform partition of $(X,d)$ as the cover $\mathcal{B}_{1/n}$ refines it. The family $\mathcal{P}_n$ is called the family of the {\it chainable components induced by $\mathcal{B}_{1/n}$}.

\begin{thm} {\rm (Theorem \ref{positive})} Let $(X,d)$ be a Bourbaki metric space and $\mathcal{P}_{n}$ be the family of all the chainable components induced by the cover of open balls $\mathcal{B} _{1/n}$. Suppose that for some $n_0\in \Nset$, for every $P\in \mathcal{P}_{n_0}$ and for every  $n> n_0$, the subfamily of chainable components $\{Q\in \mathcal{P}_n:Q\subset P\}$ does not have Ulam-measurable cardinal. Then $\upsilon X$ and $H(U_{d}(X))$ are equivalent realcompactifications.
\end{thm}

Observe that in order to prove the above theorem we have applied the results on embeddings of the previous part of the thesis and we have also proved the following two results on products and subspaces.

\begin{thm} {\rm (Corollary \ref{product.Samuel.discrete})} If $(X,d)$ is a Samuel realcompact metric  space and $(D, \chi)$ is a uniformly discrete metric space, then the Samuel realcompactification of $(X\times D, d+\chi)$ is homeomorphic to $\upsilon (X\times D)$.
\end{thm}

\begin{lm}{\rm (Lemma \ref{mapping.Samuel})} Let $(Y,\rho)$ be a metric space satisfying that $\upsilon Y$ and $H(U_\rho(Y))$ are equivalent realcompactifications. Let $(X,d)$ be a metric space such that there exists an embedding $\varphi:(X,d)\rightarrow (Y,\rho)$ satisfying that $\varphi$ is uniformly continuous and $\varphi(X)$ is a closed subspace of $Y$. Then $\upsilon X$ and $H(U_d(X)))$ are also equivalent realcompactifications. In particular, that $\upsilon Y$ and $H(U_{\rho}(Y))$ are homeomorphic realcompactifications is a property inherited by closed subspaces.
\end{lm}

We would like to know if the above results can be extended to uniform spaces. However, it is not very difficult to find an example of Bourbaki-complete uniform spaces satisfying not so  ``bad conditions'' in the cardinality of its uniform partitions and such that its Samuel realcompactification and its Hewitt realcompactification are not equivalent.  It is exactly, the product space $\beta D \times D$ where $D$ is a uniformly discrete space of Ulam-measurable cardinal (see Example \ref{example.betaD}). The reason is that to compare the Samuel realcompactification of the uniform space with its Hewitt realcompactification, is too restrictive. We should consider the $G_\delta$-closure of the space in its Samuel compactification, instead of the Hewitt realcompactification.

The $G_{\delta}$-closure of a uniform space  in its Samuel compactification is a realcompactification of the space. 	It is exactly the smallest realcompactifiation $H(\mathcal{A}(U_\mu (X))$ of $(X,\mu)$ such that not only every the real-valued uniformly continuous functions, but also all the inverse functions  $g=1/f$ where $f\in U_{\mu}(X)$ and $f(x)\neq 0$ for every $x\in X$, can be continuously extended. Therefore, for a uniform  space $(X,\mu)$ the following two facts are equivalent:

\begin{enumerate}

\item {\it $H(U_\mu(X))$ and $H(\mathcal{A}(U_\mu (X)))$ are equivalent realcompacitifcations;}

\item{\it given a uniformly continuous function $f\in U_{\mu}(X)$ which is in addition non-vanishing, its inverse $1/f$ can be continuously extended to the Samuel realcompactification $H(U_{\mu}(X))$ of $(X,\mu)$.}

\end{enumerate}

Observe that for any metric space $H(\mathcal{A}(U_{\mu}(X))$ coincides with $\upsilon X$ (see Theorem \ref{Gdelta.metric}) so in the case of metric spaces we do not have a new problem.

On the other hand, we think that the above question could help to know which functions belong to the subalgebras $C(H(U_\mu(X)))$ of all the real-valued that can be continuously extended to the Samuel realcompactification of $(X,\mu)$. In fact, we would like to know if, in the frame of Bourbaki-complete uniform spaces, the inverse function $1/f$ of every non-vanishing real-valued uniformly continuous can be continuously extended to its Samuel realcompactification.

Observe that we do know much about the subalgebra $C(H(U_{\mu}(X)))$. Just that it is the subalgebra of all the {\it Cauchy-continuous functions} any of the two modifications $wU_{\mu}(X)$ and $s_f(e\mu)$ (see Theorem \ref{cc} and Theorem \ref{countable1}). Besides, this kind of characterization does not clarify very much our questions above.

\smallskip

Finally, with the help of the following two theorems we prove that the above counterexample $(\beta D\times D, {\tt u} \times \mu_\chi)$, satisfies that $H(\mathcal{A}(U_{{\tt u} \times \mu_\chi}(\beta D\times D)))$ is equivalent to $H(U_{{\tt u} \times \mu_\chi}(\beta D\times D))$.

\begin{thm} {\rm (Theorem \ref{Samuel.product})} Let $(X,\mu)$ be a Samuel realcompact space and $(Y,\nu)$ any uniform space. Then $$H(U_{\mu \times \nu}(X\times Y))=H(U_{\mu}(X))\times H(U_{\nu}(Y))=X\times H(U_{\nu}(Y)).$$
\end{thm}

\begin{thm} {\rm (Theorem \ref{product.Gdelta})} Let $(X,\mu)$ be a Samuel realcompact space and $(Y,\nu)$ any uniform space. Then $$H(\mathcal{A}(U_{\mu \times \nu}(X\times Y)))=H(U_{\mu}(X))\times H(\mathcal{A}(U_{\nu}(Y)))=X\times H(\mathcal{A}(U_{\nu}(Y))).$$
\end{thm}

Therefore, we still do not have a counterexample of Bourbaki-complete uniform or metric space satisfying that the considered realcompactifications are not equivalent. This counterexample would be very useful in order to clarify the problem. A candidate is the uniformly 0-dimensional space $(D^{\omega_0},\rho)$ where $D$ is a uniformly discrete space having the cardinal not only Ulam-measurable, but also $\omega _1${\it -strongly compact}.
\smallskip

\newpage
\mbox{}
\thispagestyle{empty}
\newpage

\chapter{Bourbaki-complete spaces and related properties} 
\pagenumbering{arabic}
\setcounter{page}{1}
\thispagestyle{empty}

\newpage

\section{Preliminaries on uniformities and boundedness in uniform spaces}

\subsection{Completeness and different modifications of a uniform space}

\hspace{15pt}Even though, at first, this thesis was intended to live in the framework of metric spaces, it is clear by its subject, that we need to work in the more general setting of uniform spaces. Here we always consider {\it Hausdorff uniform spaces}, defined by means of families of covers of a set, instead  of entourages or pseudometrics. Thus, we are working with {\it Tychonoff spaces} since these are the {\it uniformizable} ones \cite{willard}.  From now on, by ``space" we mean always ``Tychonoff space".  

For instance, any metric space $(X,d)$ is a uniform space. In fact a base for the uniformity is given by the covers of the open balls for a fixed radius $\varepsilon >0$, that is, $\mathcal{B}_{\varepsilon}=\{B_{d}(x,\varepsilon):x \in X\}$. We will denote by $\mu _d$ the uniformity induced by the metric $d$.

Moreover, the following notation will be used along the whole thesis. Given a set  $X$, $\mathcal{C}$ a cover of $X$ and $A\subset X$, we write:
\begin{itemize}

\item $St^{0}(A,\mathcal{C})=A$

\item $St(A,\mathcal{C})= \bigcup \{C\in \mathcal{C}: C\cap A\neq \emptyset\}$

\item $\mathcal{C}^{*}=\{St(C,\mathcal{C}):C\in \mathcal{C} \}$;

\item $St(x,\mathcal{C})=St(\{x\},\mathcal{C})$

\item $St^{m}(A,\mathcal{C})=St(St^{m-1}(A,\mathcal{C}),\mathcal{C})$, $m\geq 1$

\item $St^{\infty}(A,\mathcal{C})=\bigcup_{m\in \Nset} St^{m}(A,\mathcal{C})$

\end{itemize}

Every cover $\mathcal{C}$ of a set $X$ induces a partition of $X$ in the following way. We can choose representative points $x_i$, $i\in I$ of $X$ such that $\{St^{\infty}(x_i,\mathcal{C}):i\in I\}$ is a cover of $X$ and $St^{\infty}(x_i,\mathcal{C})\cap St^{\infty}(x_j,\mathcal{C})=\emptyset$ whenever $i\neq j$. The family of sets $\{St^{\infty}(x_i,\mathcal{C}):i\in I\}$ is called the family of all the {\it chainable components of $X$ induced by $\mathcal{C}$}. 














\smallskip

The main subjects of this thesis deal with completeness-like properties of uniform spaces. By this reason we assume that  the theory about complete uniform spaces is known. Nevertheless, we will recall some facts about it along the discourse. The main references on uniform spaces that we follow are \cite{willard}, \cite{howesbook}, \cite{bourbaki},  \cite{isbellbook}, and \cite{ginsburg}.

\begin{definition} A filter $\mathcal{F}$ of a uniform space $(X,\mu)$ is a {\it Cauchy filter} if for every $\mathcal{U}\in \mu$ there is $U\in \mathcal{U}$ such that $F\subset U$ for some $F\in \mathcal{F}$.
\end{definition}

\begin{definition} A uniform space $(X,\mu)$ is {\it complete} if every Cauchy filter clusters (equivalently, converges).
\end{definition}

Moreover, we recommend the references \cite{beer-between1}, \cite{howes.completeness} and \cite{howesbook} to learn about cofinally complete uniform spaces.

\begin{definition} A filter $\mathcal{F}$ of a uniform space $(X,\mu)$ is a {\it cofinally Cauchy filter} if for every $\mathcal{U}\in \mu$ there is some $U\in \mathcal{U}$ such that $F\cap U\neq \emptyset$ for every $F\in \mathcal{F}$.
\end{definition}

\begin{definition} A uniform space $(X,\mu)$ is {\it cofinally complete} if every cofinal Cauchy filter clusters. 
\end{definition}

In particular, the following implications are clear:
\medskip

\begin{center}

{\it compact } 

$\Downarrow$ 

{\it cofinally complete } 

$\Downarrow$ 

{\it complete}
\end{center}

\smallskip

\noindent Later, we will see counterexamples of the reverses of the above implications.
\smallskip
 
We say that a uniformity $\mu$ on a Tychonoff space $X$ is {\it compatible} whenever $\mu$ induces the same topology than the topology of $X$. Given two compatible uniformities $\mu_1$ and $\mu _2$ on a space $X$ we will write  $\mu _1 \leq \mu _2$ whenever the identity map $id: (X,\mu _2)\rightarrow (X,\mu _1)$ is uniformly continuous, or equivalently, if $\mu _1\subset \mu _2$. Moreover $\mu _1=\mu _2$ if and only if $\mu _1 \leq\mu _2$ and $\mu _2\leq\mu _1$. 

Next, observe that for two compatible uniformities $\mu$ and $\nu$ on a Tychonoff space $X$, their respective completions $(\widetilde{X}, \widetilde{\mu})$ and $(\widehat{X}, \widehat{\nu})$ can be topologically homeomorphic even if $\mu$ and $\nu$ are not equivalent uniformities. More precisely recall the following definition (\cite{borsik}).
 
\begin{definition} A map $f:(X,\mu)\rightarrow (Y,\nu)$ between two uniform spaces is {\it Cauchy-continuous} if it sends Cauchy filters of $(X,\mu)$ to Cauchy filters of $(Y,\nu)$.
\end{definition}

\noindent It is clear that {\it every uniformly continuous map is Cauchy-continuous} and that {\it every Cauchy-continuous map is continuous} (\cite{dimaio}). In particular, we will use frequently the following fact: 

{\it Let $id_1 ;(X,\mu)\rightarrow (X,\nu)$ and $id_2 :(X,\nu)\rightarrow (X,\mu)$, be the identity maps. Then, if both maps are Cauchy-continuous then the completions $(\widetilde{X}, \widetilde{\mu})$ and $(\widehat{X}, \widehat{\nu})$ are topologically homeomorphic.} 
\noindent Observe that the converse is also true, that is, if there exists an homeomorphic map $h:\widetilde{X}\rightarrow \widehat{X}$, leaving $X$ point-wise fixed, then the restriction map $g=h|_X:(X,\mu)\rightarrow (X,\nu)$ is an homeomorphism satisfying that $g$ and its inverse $g^{-1}$ are Cauchy-continuous.

\smallskip

In the family $M=\{\mu _\alpha :\alpha \in A\}$ of all the compatible uniformities over $X$ the above relation $\leq$ is a partial order, and there exists a supremum of  $(M, \leq)$. This is exactly the  finest uniformity compatible with the topology of $X$.  We call it {\it the fine uniformity} and we denote it by ${\tt u}$. It is exactly the uniformity over $X$ having as a subbase the family of covers belonging to $\bigcup_{\alpha \in A} \mu _{\alpha}$.

\begin{definition} A sequence of open covers $\langle\mathcal{G}_n \rangle _{n\in \Nset}$ of a space $X$ is a {\it normal sequence} if for every $n\in \Nset$, $\mathcal{G}_{n+1} ^{*}$ refines $\mathcal{G}_{n}$.  An open cover $\mathcal{G}$ of a space $X$ is a {\it normal cover} if there is a normal sequence such that $\mathcal{G}_1=\mathcal{G}$.  (Thus every cover in a normal sequence  is a normal cover)  
\end{definition}

\noindent By the axioms of uniformity it is clear that {\it the fine uniformity {\tt u} on a space $X$ is exactly the uniformity having as a base the open normal covers of $X$} (\cite{willard}).

Moreover, it is very useful to known that {\it a base for the fine uniformity ${\tt u}$  is given by all the locally finite cozero covers from $X$} \cite{rice}. A {\it cozero} is a set $S$ such that for some real-valued continuous function $f$, $S=X\backslash f^{-1}(0)$. On the other hand, locally finite covers are defined as follows.

\begin{definition} A cover $\mathcal{U}$ of a topological space $X$ is {\it locally finite} ({\it locally countable}) if for every $x\in X$ there exists a neighborhood $V^x$ of $x$ such that $V^x$ meets only finitely (countably) many elements of $\mathcal{U}$.
\end{definition}

Recall that locally finite covers characterize the property of paracompactness.

\begin{definition} A space $X$ is {\it paracompact} if every open cover has a locally finite open refinement.
\end{definition}

\noindent In particular, {\it the fine uniformity {\tt u} on a paracompact space $X$ has as a base the family of all the open covers of the space} by the following well-known result.

\begin{lemma}\label{shrinking} {\rm (\cite{morita.locally.finite})} Let $\{G_i, i\in I\}$ be an open cover of a paracompact space $X$. Then, there exists a cover $\{Q_i:i\in I\}$ of $X$ such that, for each $i\in I$, the set $Q_i$ is a cozero set contained in $G_i$.
\end{lemma}

Recall that for the fine uniformity ${\tt u}$, paracompactness and the above defined property of cofinal completeness are related by the following result of Corson.

\begin{theorem}\label{corson}{\rm (\cite{corson})} A space $X$ is paracompact if and only if $(X,{\tt u})$ is cofinally complete.
\end{theorem}

\begin{definition} A uniform space $(X,\mu)$ is {\it fine} if $\mu={\tt u}$.
\end{definition}

A nice characterization of fine spaces in terms of maps is given by the following result.

\begin{theorem}\label{map.fine} {\rm (\cite{dimaio})} A uniform space $(X,\mu)$ is fine if and only if every continuous map $f:(X,\mu)\rightarrow (Y,\nu)$ is uniformly continuous for every uniform space $(Y,\nu)$.
\end{theorem}

Examples of fine uniform spaces are compact spaces and, in general, spaces that are uniformizable by a unique uniformity as the the space of all the countable ordinals $[0,\omega _1)$. Moreover, every {\it uniformly discrete metric space} is fine. By a uniformly discrete space we mean a set $D$ endowed with the metric $\chi :D\rightarrow [0, \infty)$ defined by $$\chi(d , e)=\begin{cases}  \,\, \,  0 &     \text { if $d= e$}  \\ \,\,\, 1  & \text{ if $d \neq e$.}\end{cases}$$

\medskip

Recall that by {\it modification} of a uniform space $(X,\mu)$ we mean a uniformity $\nu$ compatible with the topology of $X$ having as a base or subbase  a subfamily of covers from the uniformity $\mu$. Next, we are going to consider the following modifications of a uniformity $\mu$: the {\it finite modification} $f\mu$, the {\it countable modification} $e\mu$, the {\it point-finite modification} $p_f \mu$ and the {\it star-finite modification} $s_f \mu$. By \cite[Theorem 1.1]{ginsburg} and \cite[Proposition 28, Chapter IV]{isbellbook}, these modifications have as a base, respectively, the family of all the finite, countable, point-finite and star-finite open covers from $\mu$.

\begin{definition} A cover $\mathcal{U}$ of a set $X$ is {\it point-finite} ({\it point-countable}) if every point of $X$ lies only on finitely (countably) many elements of the cover. 
\end{definition}

\begin{definition} A cover $\mathcal {U}$ of a set $X$ is {\it star-finite} ({\it star-countable}) if every $U\in \mathcal{U}$ meets at most finitely (countably) many elements of $\mathcal{U}$.
\end{definition}

If $\mu=f\mu$, $\mu=e\mu$, $\mu=p_f\mu$ or $\mu=s_f\mu$, then we will say that the uniformity $\mu$ has a finite, countable, point-finite or star-finite base respectively. Other modifications are possible and we will introduce later some of them (see again \cite{ginsburg} and \cite{isbellbook}). Next, we recall some useful facts about the above modifications.
\medskip

$\bullet$ \underline{{\it The finite modification.}} The finite modification of a uniformity $\mu$ is very important because of the Samuel compactification. Recall that the {\it Samuel compactification} $s_\mu X$ of $(X,\mu)$ is, topologically, the completion of $(X, f\mu)$.

We will come back to the Samuel compactification in the second part on the thesis. However, we include now several results that we will need soon (see \cite{ginsburg}, \cite{isbellbook} and \cite{samuel}).

\begin{lemma} \label{Samuel1} Let $(X,\mu)$ be a uniform space and $(\widetilde{X},\widetilde{\mu})$ denotes its completion. Then, the Samuel compactification $s_{\widetilde{\mu}} \widetilde{X}$ is exactly $s_{\mu}X$.
\end{lemma}

The above result states that the completion $(\widetilde{X}, \widetilde{\mu})$ of $(X,\mu)$ is in particular a topological subspace of $s_{\mu} X$.

It is known, that the completion of $(X, f{\tt u})$ is exactly the Stone-\v{C}ech compactification $\beta X$ of $X$. However, if $(X,\mu)$ is any uniform space then its completion $(\widetilde{X},\widetilde{\mu})$ is not necessarily a (topological) subspace of $\beta X$. For instance, $\Rset$ is not a subspace of $\beta \mathbb{Q}$ \cite{engelkingbook}.

\begin{lemma} \label{Samuel2}Let $(X,\mu)$ be a uniform space and $\nu$ a uniformity of $X$ satisfying that $f\mu \leq\nu \leq\mu$ then the Samuel compactification $s_{\nu}X$ is exactly $s_{\mu}X$.
\end{lemma}

\begin{lemma}\label{Samuel3} Let $(X,\mu)$ be a uniform space and $A\subset X$. Then the Samuel compactification of $(A, \mu |_{A})$ is exactly the closure of $A$ in $s_{\mu}X$, that is, $s_{\mu |_A} A={\rm cl} _{s_\mu X} A$.
\end{lemma}

\begin{lemma} \label{Samuel4} Let $(X,\mu)$ be a uniform space and $A,B\subset X$ subsets. Then $${\rm cl}_{s_{\mu}X}A \cap {\rm cl}_{s_{\mu}X}B\neq \emptyset $$ if and only if $$St(A,\mathcal{U})\cap St(B,\mathcal{U})\neq \emptyset$$ for every $\mathcal{U}\in \mu$.
\end{lemma}

$\bullet$ \underline{{\it The countable modification.}} The countable modification $e\mu$ will be of relevance in Part 3 of this thesis. Just observe that every Lindel\"of uniform space $(X,\mu)$ satisfies trivially that $\mu=e\mu$. Other family of uniform spaces having a countable base for its uniformity are uniformities induced by families of real-valued functions, as we will see later.

\smallskip

$\bullet$ \underline{{\it The point-finite modification.}}  Now consider the point-finite modification $p_f\mu$. It is well-known that every countable uniform cover has a countable uniform point-finite refinement (\cite{smith}).  Therefore, in general, $p_f (e\mu)=e\mu$, where by $p_f(e\mu)$ we denote the point-finite modification of $e\mu$. Moreover, by a result in $\cite{smith}$ also, the point-finite modification $p_f\mu$ is equivalent to the uniformity having as a base all the locally finite uniform covers from $\mu$. Then, it is clear that, for the fine uniformity, ${\tt u}=p_f{\tt u}$.

Recall that in 1960 Stone \cite{stone} asked if every uniform space, in particular  every metric space, has a locally finite, equivalently, point-finite base for its uniformity. His question was motivated by the fact that every metric space is paracompact and hence every open cover has a locally finite refinement, and from the additional fact that ${\tt u}=p_f{\tt u}$. So he wanted to extended this problem from the topological level to the uniform level.
This was answered negatively by E. \v{S}\v{c}epin and J. Pelant in 1975 (see \cite[p. 596]{aull}, \cite{pelant-combinatorial} and see  \cite{avart1} and \cite{avart2} for additional information). We will consider later Pelant's couterexample $\ell _\infty (\omega _1)$.

\smallskip

$\bullet$ \underline{{\it The star-finite modification.}} Finally, let us take into account the star-finite modification. In particular, we study the star-finite modification $s_f {\tt u}$ of the fine uniformity {\tt u}, since the modifications of the fine uniformity will be of special interest. 

It is easy to see that, since a base of ${\tt u}$ is given by the locally finite cozero covers of $X$, then {\it a base for the finite modification $f{\tt u}$ and for the countable modification $e{\tt u}$ is given, respectively, by the finite and countable cozero covers of $X$} (see also \cite{garcia.special}). However, the proof of the same fact but for the star-finite modification $s_f {\tt u}$ is not immediate.

\begin{lemma} {\rm (\cite{morita.star-finite}, or \cite[Lemma 5.2.4]{engelkingbook})} \label{morita} Every countable cozero cover of a space $X$ has a countable star-finite cozero refinement.
\end{lemma}

\begin{theorem} \label{base.star-finite} The family of all the star-finite cozero covers of a space $X$ is a base for $s_f{\tt u}$.
\begin{proof} Observe that in \cite{garcia.special} it is proved that the star-finite cozero covers of $X$ are normal covers. On the other hand, let $\mathcal{C}$ a star-finite open cover belonging to ${\tt u}$.  If $\{St^{\infty}(x_i,\mathcal{C}):i \in  I\}$ is the family of chainable components generated by $\mathcal{C}$ then, for every $i \in I$, $St^{\infty}(x_{i}, \mathcal{C})$ contains at most countably many elements $C\in \mathcal{C}$ (see \cite{engelkingbook}). Now, for every $i\in I$ the cover $$\mathcal{C}_i=\{C\in \mathcal{C}: C\subset St^{\infty}(x_i,\mathcal{C})\}\cup \bigcup \{St^{\infty}(x_j,\mathcal{C}): j\in I, j\neq i\}$$ is a countable open cover belonging to ${\tt u}$ as $\mathcal{C}$ refines $\mathcal{C}_i$. As we have previously noticed, we can take a cozero refinement $\mathcal{A} _i$ of $\mathcal{C}$ which is countable. More precisely, we can take $\mathcal{A}_i$ being countable and star-finite by the result of Morita, Lemma \ref{morita}. 

Next, consider the cover $$\mathcal{A}=\{A\in \mathcal{A}_i: A\subset St^{\infty}(x_i,\mathcal{C}), i\in I\}.$$ Then, it is clear that $\mathcal{A}$ is the desired cozero star-finite refinement of $\mathcal{C}$.
\end{proof}
\end{theorem}

By Lemma \ref{morita} and the above result the next corollary is clear.

\begin{corollary} \label{morita2} Let $X$ be a space and ${\tt u}$ the fine uniformity on it. Then $e{\tt u}=s_f(e{\tt u})$.
\end{corollary} 

\begin{remark}\label{star-countable-cozero} Similarly to Theorem \ref{base.star-finite} one can prove that the uniformity $s_f{\tt u}$, has also a base given by all the star-countable cozero covers of $X$. 
Nevertheless, in general, for a uniform space $(X,\mu)$, if we consider the family of all the star-countable uniform covers of $X$, then this family is not necessarily a base for a uniformity but just for a quasi-uniformity $s_c \mu$ \cite[p. 368]{reynolds}. On the other hand, for uniform spaces satisfying that $\mu =e\mu$ it is clear that $s_c \mu$ is a uniformity since $s_c\mu =e\mu=\mu$.
\end{remark}

The following lemma relates the  star-finite modification $s_f{\tt u}$ with the countable modification $e{\tt u}$ of a connected space. We assume that this result is known, implicitly at least, in the papers of Morita \cite{morita.star-finite} and Wiscamb \cite{wiscamb}.  We will need it several times in this thesis.

\begin{lemma} \label{cozero.connected} Let $X$ be a  connected space. Then $s_f{\tt u}=e{\tt u}$. 
\begin{proof} By  Lemma \ref{morita} it is clear that $e{\tt u}\leq s_f{\tt u}$. So, conversely, let us take $\mathcal{C} \in s_f{\tt u}$ a star-finite cozero cover of $X$. We are going to show that it is countable. Let $\{St^{\infty}(x_i,\mathcal{C}):i \in  I\}$ be the family of chainable components generated by $\mathcal{C}$, then, since $\mathcal{C}$ is star-finite, for every $i \in I$, $St^{\infty}(x_{i}, \mathcal{C})$ contains at most countably many elements $C\in \mathcal{C}$. Moreover, it is clear that each chainable component $St^{\infty}(x_i,\mathcal{C})$ is an open subset of $X$. But as $X$ is connected, we must have that $St^{\infty}(x_i,\mathcal{C})\cap St^{\infty}(x_j,\mathcal{C})\neq \emptyset$ for every $i,j\in I$. But this equivalent to have that $St^{\infty}(x_i,\mathcal{C})= St^{\infty}(x_j,\mathcal{C})$ for every $i,j\in I$, that is, there is only one chainable component in $X$ induced by $\mathcal{C}$. Therefore, $\mathcal{C}$ is a countable cover and then $\mathcal{C}\in e{\tt u}$. 
\end{proof}
\end{lemma}

$\bullet$ \underline{{\it Weak uniformities.}}  Examples of uniformities having a base of star-finite covers can be found in \cite{isbell.weak}. For instance, we are interested, mainly in Part 3, in  the uniformities called 
{\it weak uniformities}. By weak uniformity we mean a uniformity generated by a family of real-valued functions. More precisely, given a family $\mathcal{L}\subset C(X)$ of real-valued continuous functions over a space $X$ ($\Rset$ is always endowed with the euclidean topology), then the covers of the form $$\{\{y\in X: |f(x)-f(y)|<\varepsilon\}: x\in X\}\text{, for every } f\in \mathcal{L} \text{ and every } \varepsilon >0,$$ are a subbase for a uniformity $w\mathcal{L}$ on $X$ which is compatible with the topology of $X$ if and only if the family $\mathcal{L}$ {\it separates points from closed sets} of  $X$, that is, for every $x\in X$ and every $F\subset X$ closed subset such that $x\notin F$ there exists some $f\in \mathcal{L}$ such that $f(x)\notin {\rm cl}_\Rset f(F)$. Observe that  a weak uniformity, generated by a family of continuous real-valued functions ${\mathcal{L}}$, has a star-finite and countable base composed by finite intersections  of covers of the form $\{f^{-1}((n-1 )\cdot\varepsilon,(n+1)\cdot \varepsilon)): n\in \mathbb{Z}\}$ for $f \in \mathcal{L}$ and $\varepsilon >0$. In general, the countable covers $\{C_n:n\in \Nset\}$ of a set $X$ satisfying that $C_i \cap C_j = \emptyset $ whenever $|i-j|>1$ are called {\it linear covers} in \cite{isbell.weak} and 2-{\it finite covers} in \cite{garrido-montalvo}.

In the particular case of the weak uniformity $w{U}_{\mu} (X)$, generated by the family $U_{\mu}(X)$ of all the real-valued uniformly continuous functions on the uniform space $(X,\mu)$, we have that a base of $w{U}_{\mu} (X)$  is given by finite intersections of open linear covers from the uniformity $\mu$ (see \cite[1.6]{isbell.weak}). That is, $w{U}_{\mu} (X)$ can be considered also a modification of $\mu$.

\begin{theorem} \label{real.weak}The uniformity $\mu_{d_u}$ over the real-line $\Rset$ induced by the usual euclidean metric $d_u$ is exactly the weak uniformity $wU_{\mu_{d_u}}(\Rset)$. Moreover the product uniformity $\pi$ over any product $\Rset ^{\alpha}$ coincides also with the weak uniformity $wU_{\pi}(\Rset^{\alpha})$.
\begin{proof} The indentity function $g(x)=x$ is a uniformly continuous real-valued function. Since for every $x\in X$, $$\{y\in X:|g(x)-g(y)|<\varepsilon\}=B_{d}(x,\varepsilon)$$ then $wU_{\mu_{d_u}} (\Rset)=\mu_{d_u}$.

Now, recall that the product uniformity $\pi$ on $\Rset ^{\alpha}$ is exactly the weak uniformity generated by all the projection maps each coordinate $p_\kappa:\Rset^{\alpha}\rightarrow \Rset$, $\kappa <\alpha$. Thus, $\pi$ is the smallest uniformity making the projections $p_{\kappa}$ uniformly continuous, that is $\pi\leq wU_{\pi}(X)$. But, by definition of the weak uniformity $wU_{\pi}(X)$ satisfies that $\pi\geq wU_{\pi}(X)$ and then $\pi= wU_{\pi}(X)$.
\end{proof}
\end{theorem}

We finish this section relating all the modifications appeared until now. For a uniform space $(X,\mu)$ the following relations are satisfied: 
$$\mu\geq p_f\mu \geq s_f\mu\geq wU_{\mu}(X)\geq f\mu$$
$$\mu \geq p_f\mu\geq p_f(e\mu)=e\mu\geq wU_{\mu}(X) \geq f\mu .$$ 
\smallskip

On the other hand, whenever we consider {\it the weak uniformity $wC (X)$, generated by the family $C(X)$ of all the real-valued continuous functions over the space $X$, then it has as a subbase  all the linear cozero covers of $X$} \cite{garrido-montalvo}. Precisely $wC(X)=wU_{{\tt u}}(X)$ as $C(X)=U_{\tt u}(X)$ (Theorem \ref{map.fine}). In particular, $wC(X) \leq e{\tt u}$. However, $wC(X)\neq e{\tt u}$ (see \cite[Exercise 8.1.I b)]{engelkingbook}).

\noindent Therefore, the linear covers are determinant for weak uniformities.

Next, let $U^{*}_\mu(X)$ the family of all bounded real-valued uniformly continuous functions on a uniform space $(X,\mu)$. It is well-known that $f\mu=wU^{*}_{\mu}(X)$ \cite{isbellbook}. Similarly, let $C^{*}(X)$ denote the family of all the bounded real-valued continuous functions over a space $X$, then $wC^{*}(X)=wU^{*}_{{\tt u}}(X)=f{\tt u}$.

\medskip

\subsection{Bourbaki-bounded subsets and other bornologies}

\hspace{15pt} The initial project of this thesis started by studying the next notion of Bourbaki-boundedness. This can be defined in the frame of uniform spaces and in the particular frame of metric spaces, as a generalization of metric boundedness.

\begin{definition} A subset $B$ of a uniform space $(X,\mu)$ is a {\it Bourbaki-bounded subset} in $(X,\mu)$ if for every uniform cover $\mathcal{U}\in \mu$ there exist  $m\in \Nset$ and  finitely many $U_1,...,U_k\in \mathcal{U}$ such that $$B\subset \bigcup _{i=1}^{k}St^{m}(U_i,\mathcal{U}).$$ In particular $(X,\mu)$ is a {\it Bourbaki-bounded space} if it is a Bourbaki-bounded subset in itself.
\end{definition}

For the definition in metric spaces we need the following notation. Let $\varepsilon >0$ and $x\in X$, then we write $$B^{m}_{d}(x,\varepsilon)=St^{m-1}(B_{d}(x,\varepsilon),\mathcal{B}_{\varepsilon})\text{ for every }n \in\Nset$$ $$B^{\infty}(x,\varepsilon)=\bigcup _{m\in \Nset}B^m (x,\varepsilon).$$ Observe that, for $x\in X$, if $y\in B^m _d (x,\varepsilon)$ for some $\varepsilon >0$ and $m\in \Nset$, then we can choose a ``chain" of points $x_i\in X$, $i=0,...,m$, of ``lenght" $m$, such that $x_0=x$, $x_m=y$ and $d(x_{i-1},x_i)< \varepsilon$ for every $i=1,...,m$.

\begin{definition} A subset $B$ of a metric space $(X,d)$ is a {\it Bourbaki-bounded subset} in $X$ if for every $\varepsilon >0$ there exist a natural number $m\in \Nset$ and finitely many points $x_1,...,x_k\in X$ such that $$B\subset \bigcup _{i=1}^{k}B^m _{d}(x_i, \varepsilon).$$ In particular $(X,d)$ is a {\it Bourbaki-bounded space} if it is a Bourbaki-bounded subset in itself.
\end{definition}

We will denote by ${\bf BB}_{d}(X)$ and ${\bf BB}_{\mu}(X)$ the family of all the Bourbaki-bounded subsets of a metric space $(X,d)$ and of a uniform space $(X,\mu)$, respectively. Moreover, we will denote by ${\bf B}_{d}(X)$ the family of all the bounded subsets of a metric space, that is, those subsets $B$ such that for some $x\in X$ and $K>0$, $B\subset B_{d}(x, K)$.

It is clear that {\it every Bourbaki-bounded subset of a metric space $(X,d)$ is bounded by the metric $d$}, that is, ${\bf BB}_d(X)\subset {\bf B}_{d}(X)$, since $B_{d}^{m}(x, \varepsilon)\subset B_{d}(x,m\cdot\varepsilon)$. Moreover, we have the following result

\begin{theorem} \label{boundedness.banach} In a normed space $(E,d_{||\cdot ||})$, where $d_{||\cdot ||}$ is the metric induced by the norm, every bounded subset is Bourbaki-bounded in $X$.
\begin{proof}We need to prove only one inclusion. For every $\varepsilon >0$ and every $K >0$, if we take some $m\in \Nset$, $m> \frac{K}{\varepsilon}$ then $B_{d_{||\cdot ||}}(0,K)\subset B^m _{d_{||\cdot ||}}(0, \varepsilon )$. Since for every $B\in {\bf B}_{d_{||\cdot ||}}(X)$ there exists some $K>0$ such that $B\subset B_{d_{||\cdot ||}}(0,K)$, the result follows.

\end{proof}
\end{theorem}

However, {\it Bourbaki-boundedness is invariant under uniform transformations} while metric boundedness is not. Indeed,  given and unbounded metric space $(X,d)$ the metric $d^{*}=min \{1, d\}$ is bounded and uniformly equivalent to $d$, but not every bounded subset of $(X, d^{*})$ is bounded in $(X,d)$. Therefore, ${\bf BB}_d (X)={\bf BB}_{d^{*}}(X)$ and ${\bf B}_{d}(X)\neq {\bf B}_{d^*}(X)$.

\begin{remark}Hejcman called the Bourbaki-bounded subsets of a uniform space, just {\it bounded}. In fact, it was Hejcman the first that studied deeply boundedness in uniform spaces in his paper \cite{hejcman} of 1959. This paper is always cited as the main reference on the subject. In the same year Atsuji \cite{atsuji1} considered also the Bourbaki-bounded metric spaces that he called {\it finitely-chainable}. The name of Bourbaki-bounded comes from the book of Bourbaki where in an exercise \cite[II Exercise 4.7]{bourbaki} some properties of them are listed. Observe that Bourbaki also called these subsets simply bounded ({\it borne\'{e}} in the French edition). Moreover, these subsets where already called {\it bounded in the sense of Bourbaki} in \cite{vroegrijk1}

\end{remark}

Now, recall the following definition.

\begin{definition} Let $X$ be a (uniform) space and $d$ a pseudometric on $X$. We say that $d$ is a {\it (uniformly) continuous
pseudometric} of X if $d$ is (uniformly) continuous as a mapping from the product space $X\times X$ to $\Rset$.
\end{definition}

\noindent The following result of Hejcman is an interesting characterization of the  Bourbaki-bounded subsets.

\begin{theorem}\label{hejcman}{\rm (\cite{hejcman})}. Let $(X,\mu)$ be a uniform space and $B$ a subset of $X$. The following statements are equivalent:
\begin{enumerate}
\item $B$ is a Bourbaki-bounded subset in $X$;

\item $f(B)$ is a bounded subset of $(\Rset,d_u)$  for every real-valued uniformly continuous function $f\in U_{\mu}(X)$;

\item $B$ is a bounded subset of $(X,\rho)$ for every uniformly continuous pseudometric $\rho$ on $X$.
\end{enumerate}
\end{theorem}

\begin{corollary} Let $(X,d)$ be a metric space and $B$ a subset of $X$. The following statements are equivalent:
\begin{enumerate}
\item $B$ is a Bourbaki-bounded subset in $X$, that is, $B\in {\bf BB}_{d}(X)$;

\item $B$ is a bounded subset of $(X,\rho)$, that is, $B\in {\bf B}_\rho(X)$ for every uniformly equivalent metric $\rho$.
\end{enumerate}
\end{corollary}

\medskip

Now, observe that {\it Bourbaki-bounded subsets are a generalization of the totally-bounded subsets}.

\begin{definition} A subset $B$ of a uniform space $(X,\mu)$ is {\it totally bounded} if for every $\mathcal{U}\in \mu$ there exist finitely many $U_1,...,U_k\in \mathcal{U}$ such that $$B\subset \bigcup _{i=1} ^k U_i .$$ 
\end{definition}

We will denote by ${\bf TB}_{\mu}(X)$ and ${\bf TB}_{d}(X)$, respectively, the family of all the totally bounded subsets of a uniform space $(X,\mu)$ and of a metric space $(X,d)$. Thus, ${\bf TB}_{\mu}(X)\subset {\bf BB}_{\mu}(X)$. More precisely, if $B$ is a Bourbaki-bounded subset of $(X,\mu)$ such that there exists some $M\in \Nset$ satisfying that for every $\mathcal{U} \in \mu$ the value $m$ (from the definition of Bourbaki-bounded subset), depending of $\mathcal{U}$, is always bounded by $M$, then $B$ is a totally bounded subset.
 
However, in general, not every  Bourbaki-bounded subset is totally bounded and, even though total boundedness is invariant under uniform homeomorphisms, it is too restrictive to be chosen as  the uniform  notion of boundedness. In fact, if $(E, ||\cdot ||)$ is an infinite-dimensional Banach space then the closed unit ball of $E$ is not totally bounded because otherwise it would be compact by completeness of $E$. But this is a contradiction of infinite-dimensionality.

There is another big difference between total boundedness and Bourbaki-boundedness. The definition of Bourbaki-bounded subset depends on the space of ambience. For instance, it is clear that the canonical base $\{e_n:n\in \mathbb{N}\}$ of the Hilbert space of square summable sequences of real numbers $\ell_2$ is a Bourbaki-bounded subset of $\ell _2$ but it is not a Bourbaki-bounded space in itself. However, if a subset $B$ is totally bounded in some uniform space $X$ then it will be always totally bounded in any uniform space $Y$ in which $B$ is uniformly embedded.


\medskip

The above families of bounded subsets, metric bounded subsets, Bourbaki-bounded subsets and totally bounded subsets, have in common that all of them form different {\it closed bornologies}, that is, they are covers of $X$, stable under finite unions, subsets and closure \cite{hejcman}. 

There is another closed bornology that will be relevant along this thesis. It is the bornology ${\bf CB}(X)$ of all the subsets $B$ of a topological space $X$ such that the closure of $B$ is compact, that is, the {\it relatively compact subsets}. The relevance of this last bornology is in the following well-known result.

\begin{theorem} \label{complete.metric} A metric space $(X,d)$ is complete if and only if ${\bf CB}(X)={\bf TB}_{d}(X)$
\end{theorem}


For a uniform space $(X,\mu)$ and for a metric space $(X,d)$, all the above bornologies are related, respectively, as follows:

$${\bf BB}_{\mu}(X)\supset{\bf TB}_{\mu}(X)\supset{\bf CB}(X).$$

$${\bf B}_{d}(X)\supset{\bf BB}_{d}(X)\supset{\bf TB}_{d}(X)\supset{\bf CB}(X).$$ Next, we give an example of a metric space such that all the above bornologies are different. To that purpose, we take into account that all these bornologies are {\it finitely productive}.


\begin{example} {\it There exists a metric space $(X,d)$ such that $$ {\bf B}_{d}(X)\supsetneq {\bf BB}_{d}(X)\supsetneq {\bf TB}_{d}(X)\supsetneq {\bf CB}(X). $$}

\begin{proof} Consider the product $X=\ell _2 \times (0,\infty)$ endowed with the product metric given by the sum $d=d_{||\cdot ||_2}+d_u^*$ of the metric $d_{||\cdot ||_2}$ induced by the norm $||\cdot ||_2$ of $\ell _2$ and the bounded metric $d^{*}_{u}=min\{1, d_u\}$ where $d_u$ denoted the usual euclidean metric restricted over the interval $(0,\infty)$. Then, $${\bf B}_{d_{||\cdot ||_2}}(\ell_2)={\bf BB}_{d_{||\cdot ||_2}}(\ell_2)\neq {\bf TB}_{d_{||\cdot ||_2}}(\ell_2)={\bf CB}(\ell _2) \text{ and}$$ $${\bf B}_{d^* _u}((0,\infty))\neq{\bf BB}_{d^* _u}((0,\infty))= {\bf TB}_{d^* _u}((0,\infty))\neq{\bf CB}((0,\infty)).$$ Hence, by productivity,  we have: $$ {\bf B}_{d}(X)\supsetneq {\bf BB}_{d}(X)\supsetneq {\bf TB}_{d}(X)\supsetneq {\bf CB}(X). $$
\end{proof}
\end{example}

\bigskip

\bigskip

\begin{center} \Denarius \Denarius  \Denarius \end{center}
\bigskip

\bigskip

\section{Bourbaki-complete and cofinally Bourbaki-complete uniform spaces}

\subsection{Bourbaki-boundedness by means of filters and star-finite modification of a uniformity}

\hspace{15pt} Totally bounded subsets of a uniform space can be characterized by means of Cauchy filters and cofinally Cauchy filters imitating the well-known characterization of totally bounded subsets by Cauchy sequences and cofinally Cauchy sequences in the metric context that we will recall later. Indeed, the following result is satisfied.

\begin{theorem}\label{totally.bounded.filters} Let $B$ a subset of a uniform space $(X,\mu)$. The following statements are equivalent:
\begin{enumerate}
\item $B$ is totally bounded; 

\item every filter $\mathcal{F}$ in $B$ is cofinally Cauchy;

\item every ultrafilter $\mathcal{F}$ in $B$ is Cauchy;

\item every filter $\mathcal{F}$ in $B$ is contained in some Cauchy filter $\mathcal{F}'$ in $B$.
\end{enumerate}
\begin{proof} The proof of the equivalence of $(1)$, $(3)$ and $(4)$ can be found in \cite[Theorem 5.4]{warren}. The statement $(2)$ is a little bit unusual but it comes from metric spaces (see \cite[Proposition 3.13]{beer-between1}). 

$(1)\Rightarrow(2)$. If ${B}$ is totally bounded and $\mathcal{F}$ is a filter in $B$, then for every $\mathcal{U}\in \mu$ there exists finitely many $U_1,...,U_k\in \mathcal{U}$ such that $B\subset \bigcup _{i=1} ^k U_i$. In particular, there is some $i_0\in \{1,...,k\}$ such that $F\cap U_{i_0}\neq \emptyset$ for every $F\in \mathcal{F}$. Assume by the contrary that for every $i\in \{1,...,k\}$ we can take some $F_i\in \mathcal{F}$ such that $F_i\cap U_i=\emptyset$.  Then $(\bigcap _{i=1}^k F_i)\cap (\bigcup _{i=1}^k U_i)= \emptyset$. Since   $\bigcap_{i=1}^{k}F_i\subset \bigcup _{i=1} ^k U_i$ then $\bigcap_{i=1}^{k}F_i=\emptyset$ which is a contradiction as $\bigcap_{i=1}^{k}F_i$ is an element of the filter $\mathcal{F}$.

$(2)\Rightarrow (3)$ This implication follows from the easy fact that every cofinally Cauchy ultrafilter is a Cauchy ultrafilter by maximality.
\end{proof}
\end{theorem}

Next, we want to characterize also Bourbaki-boundedness by means of filters. This characterization will lead us to the property of Bourbaki-complete- ness and cofinally Bourbaki-completeness. The idea is to imitate the characterization of total boundedness. 

\begin{definition} A filter $\mathcal{F}$  of a uniform space $(X,\mu)$ is {\it Bourbaki-Cauchy} in $X$ if for every uniform cover $\mathcal{U}\in \mu$ there is some $m\in \Nset$ and $U\in \mathcal{U}$ such that $$F\subset St^{m}(U,\mathcal{U}) \text{ for some }F\in \mathcal{F}.$$ 
\end{definition}

\begin{definition} A filter $\mathcal{F}$  of a uniform space $(X,\mu)$ is {\it cofinally Bourbaki-Cauchy} in $X$ if for every uniform cover $\mathcal{U}\in \mu$ there is some  $m\in \Nset$ and $U\in \mathcal{U}$ such that $$F\cap St^{m}(U,\mathcal{U}) \neq \emptyset \text{ for every }F\in \mathcal{F}.$$ 
\end{definition}

More precisely, we want to prove the next result.

\begin{theorem}\label{BB.filters} Let $B$ a subset of a uniform space $(X,\mu)$. The following statements are equivalent:
\begin{enumerate}
\item $B$ is Bourbaki-bounded in $X$; 

\item every filter $\mathcal{F}$ in $B$ is cofinally Bourbaki-Cauchy in $X$;

\item every ultrafilter $\mathcal{F}$ in $B$ is Bourbaki-Cauchy in $X$;

\item every filter $\mathcal{F}$ in $B$ is contained in some filter $\mathcal{F}'$ in $B$ which is Bourbaki-Cauchy in $X$.
\end{enumerate}
\end{theorem}

In order to prove the above result, we first show that Bourbaki-boundedness in $(X,\mu)$ is just total boundedness for uniformities having a star-finite base, generalizing an old result of Nj\aa stad.

\begin{theorem} \label{BB.sf}{\rm (\cite[Theorem 1]{njastad} \cite[Lemma 7]{vroegrijk1})}  Let $B$  a subset of a uniform space $(X,\mu)$. The following statements are equivalent:
\begin{enumerate}
\item $B$ is  a Bourbaki-bounded subset in $(X,\mu)$;

\item $B$ is a totally bounded subset in $(X,s_f\mu)$;

\item for every uniformity $\nu$ on $X$  such that $s_f\mu\geq  \nu \geq wU_\mu(X)$ then $B$ is a totally bounded subset in $(X, \nu)$.
\end{enumerate}
\begin{proof} $(1)\Rightarrow (2)$  Let $\mathcal{U}\in \mu$ being star-finite. Then, by Bourbaki-boundedness there exist finitely many elements $U_1,...,U_k\in \mathcal{U}$ and $m\in \Nset$ such that $B\subset \bigcup_{i=1}^{k} St^{m}(U_i,\mathcal{U})$. Since $\mathcal{U}$ is star-finite then, for every $i=1,...,k$, there are only finitely many $U\in \mu$ such that $U\subset St^{m}(U_i,\mathcal{U})$. Therefore $B$ is totally bounded in $(X, s_f\mu)$. 

\smallskip
$(2)\Rightarrow (3)$ This is trivial.

\smallskip
$(3)\Rightarrow (1)$ Let us take $\nu =wU_{\mu}(X)$. Since $B$ is totally bounded in $(X,wU_{\mu}(X))$, in particular it is Bourbaki-bounded. Moreover, the family of real-valued uniformly continuous functions on $(X,wU_{\mu}(X))$ is exactly $U_{\mu}(X)$, then $B$ is Bourbaki-bounded in $(X,\mu)$ by Theorem \ref{hejcman}.
\end{proof}
\end{theorem}

Next, we introduce a couple of lemmas that we need for the proof of Theorem \ref{BB.filters} and that we will use several times along this thesis.

\begin{lemma} \label{corona} Let $\mathcal{U}$ be a uniform cover of a uniform space $(X,\mu)$. Then by $\mathcal{A}( \mathcal{U})$ we denote the  cover generated by $\mathcal{U}$ in the following way. Let $\{St ^{\infty}(x_i, \mathcal{U}),$ $ i \in I\}$ be the family of all the chainable components induced by $\mathcal{U}$. For every $i\in I$, write $A_1 ^i=St^2 (x_i , \mathcal{U})$,  and for every $n\in \Nset$, $n>1$,   let $$A_n ^i=\bigcup \{V\in \mathcal{V}:V\subset St^{n+1} (x_i, \mathcal{U}), V\cap ( X\backslash St^{n-1} (x_i, \mathcal{U}))\neq \emptyset\}.$$ Put $\mathcal{A}(\mathcal{U})= \{A_n ^i, n\in \Nset, i\in I\}$. Then  $\mathcal{A}(\mathcal{U})$ is a uniform cover satisfying that $A^i _n \cap A^j _m=\emptyset$ if $i\neq j$ and $A^i _n \cap A^i _m =\emptyset$ if $|n-m|>1$. In particular, $\mathcal{U}$ refines $\mathcal{A}(\mathcal{U})$ and $\mathcal{A}(\mathcal{U})\in s_f \mu$.
\begin{proof} It is clear.
\end{proof}
\end{lemma}

\begin{lemma}\label{ultrafilter}  Let $(X,\mu)$ be a uniform space. The following statements are satisfied:
\begin{enumerate}
\item A Cauchy filter of $(X,s_f \mu)$ is a  Bourbaki-Cauchy filter in $(X,\mu)$. 

\item A Bourbaki-Cauchy ultrafilter of $(X, \mu)$ is a Cauchy ultrafilter in $(X, s_f\mu)$.

\item $\mathcal{F}$ is a cofinal Cauchy filter of $(X,s_f \mu)$ if and only if $\mathcal{F}$ is a cofinal Bourbaki-Cauchy filter of $(X,\mu)$.

\end{enumerate}
\begin{proof} $1)$ Let $\mathcal{F}$ be a Cauchy filter of $(X,s_f \mu)$ and $\mathcal{U}\in \mu$. Let $\mathcal{A}(\mathcal{U})=\{A^i _n:i\in I, n\in \Nset\}$ the cover from Lemma \ref{corona} induced by $\mathcal{U}$. 
Since $\mathcal{A}(\mathcal{U})\in s_f\mu$, by hypothesis there exists some $A_n ^i \in \mathcal{A}(\mathcal{U})$ such that $A_n^i\in \mathcal{F}$. But $A_n^i\subset St^{n+1} (x_i,\mathcal{U})$ by construction. Hence, $St^{n+1} (x_i,\mathcal{U})\in \mathcal{F}$  and then  $St^{n+1} (U,\mathcal{U})\in \mathcal{F}$ for some $U\subset St^{n+1} (x_i,\mathcal{U})$. Thus, we can deduce that $\mathcal{F}$ is in fact a  Bourbaki-Cauchy filter of $(X,\mu)$.

$2)$ Let $\mathcal{F}$ be a Bourbaki-Cauchy ultrafilter of $(X,\mu)$. Then it is clear that $\mathcal{F}$ is also a Bourbaki-Cauchy ultrafilter of $(X,s_f \mu)$. Let $\mathcal{U}\in s_f \mu$ and $\mathcal{V}$ a star-finite uniform refinement. Then, for some $m\in \Nset$ and $V\in \mathcal{V}$, $St^{m}(V,\mathcal{V})\in \mathcal{F}$. But $\mathcal{V}$ is star-finite, therefore we can choose finitely many $V_i \in \mathcal{V}$, $i=1,...,k$ such that $St^{m}(V,\mathcal{V})\subset \bigcup_{i=1}^{k}V_i$. Since $\mathcal{F}$ is an ultrafilter then for some $i=1,...,k$, $V_i\in \mathcal{F}$. Finally, since $\mathcal{V}$ is a refinement of $\mathcal{U}$, there is some $U\in \mathcal{U}$ such that $U\in \mathcal{F}$. Then $\mathcal{F}$ is Cauchy in $(X,s_f \mu)$.

$3)$ The proof that every cofinal Cauchy filter of $(X,s_f \mu)$  is a cofinal Bourbaki-Cauchy filter of $(X,\mu)$ is exactly similar to $1)$ .

Next, let $\mathcal{F}$ be a cofinal Bourbaki-Cauchy filter of $(X,\mu)$. We want to prove that it is cofinal Cauchy in $(X, s_f\mu)$. The proof is also similar to $2)$. Let $\mathcal{U}\in s_f \mu$ and $\mathcal{V}$ a star-finite uniform refinement. Then, for some $m\in \Nset$ and $V\in \mathcal{V}$, $St^{m}(V,\mathcal{V})\cap F\neq \emptyset$ for every $F\in \mathcal{F}$. But $\mathcal{V}$ is star-finite, therefore we can choose finitely many $V_i \in \mathcal{V}$, $i=1,...,k$ such that $St^{m}(V,\mathcal{V})\subset \bigcup_{i=1}^{k}V_i$. In particular, for some $i\in\{1,...,k\}$, $V_i\cap F\neq \emptyset$ for every $F\in \mathcal{F}$. Otherwise, for every $i=1,...,k$ there is some $F_i\in \mathcal{F}$ such that $V_i\cap F_i=\emptyset$. But then $F'=\bigcap_{i=1} ^k F_i$ is an element of the filter $\mathcal{F}$ such that $F' \cap St^m (V,\mathcal{V})=\emptyset$, which is a contradiction. Hence $\mathcal{F}$ is a cofinally Cauchy filter of $(X, s_f\mu)$.
\end{proof}
\end{lemma}

Finally, it is clear that applying Theorem \ref{totally.bounded.filters}, and Theorem \ref{BB.sf} and Lemma \ref{ultrafilter}, the proof of Theorem \ref{BB.filters}  follows at once.

\smallskip

On the other hand, observe that, as well as the definition of Bourbaki-bounded subset depends on the space of ambience, the same happens to the definition of Bourbaki-Cauchy filter. This not happens to Cauchy filters. Note that, if we restrict the uniformity $s_f\mu$ to a subset  $B$ of $(X,\mu)$, then it can be strictly weaker than the star-finite modification of the restriction of $\mu$ over $B$. 

\medskip

\subsection{Between compactness and completeness}

\hspace{15pt} If we ask the clustering of the Bourbaki-Cauchy or cofinally Bourbaki-Cauchy filters we have the following two properties.

\begin{definition} A uniform space $(X,\mu)$ is {\it Bourbaki-complete} whenever every Bourbaki-Cauchy filter clusters.
\end{definition}

\begin{definition} A uniform space $(X,\mu)$ is {\it cofinally Bourbaki-complete} if every cofinally Bourbaki-Cauchy filter clusters.
\end{definition}

\smallskip

Observe that {\it every Cauchy filter is Bourbaki-Cauchy}, {\it every cofinally Cauchy filter is cofinally Bourbaki-Cauchy} and {\it every Bourbaki-Cauchy filter is cofinally Bourbaki-Cauchy}. Thus, we get the following implications:
\medskip

\begin{center}

{\it cofinally Bourbaki-complete} 

$\Downarrow$ 

{\it Bourbaki-complete }

$\Downarrow$ 

{\it complete}
\smallskip
\end{center}
and
\begin{center}
{\it cofinally Bourbaki-complete }

$\Downarrow$ 

{\it cofinally complete } 

$\Downarrow$ 

{\it complete}
\smallskip

\end{center}

\noindent In addition, compactness trivially implies all the above properties. By this reason we talk of properties lying between compactness and completeness (see \cite{beer-between1}, \cite{beer-between2} for similar ideas). The reverses of the above implications are in general not true as we will see in future examples.

Next, we are going to consider Bourbaki-completeness and cofinal Bourbaki-completeness for the star-finite modification of a uniformity.

\begin{theorem} \label{star-finite1} Let $(X,\mu)$ be a uniform space. The following statements are satisfied:
\begin{enumerate} 

\item Let $Z$ be the subspace of $s_\mu X$ of all the clusters points of the Bourbaki-Cauchy filters of $(X,\mu)$. Then $Z$ is topologically homeomorphic to the completion of $(X,s_f \mu)$.

\item Let $Z_1$ and $Z_2$ be the subspaces of $s_\mu X$ of all the clusters points of the cofinally Cauchy filters of $(X, s_f\mu)$ and of all the cofinally Bourbaki-Cauchy filters of $(X,\mu)$, respectively. Then $Z_1=Z_2$.

\end{enumerate}
\begin{proof} $1)$ Let $(\widetilde{X},\widetilde{s_ f \mu})$ denote the completion of $(X,s_f \mu)$. Recall that, by Lemma \ref{Samuel1} and Lemma \ref{Samuel2}, $\widetilde{X}$ is a subspace of $s_\mu X$. Therefore, we are going to prove that $Z= \widetilde{X}$. By Lemma \ref{ultrafilter} it is clear that $\widetilde{X}\subset Z$. Conversely, let $\mathcal{F}$ be a Bourbaki-Cauchy filter in $(X,\mu)$ and $\xi \in Z$ a cluster point of $\mathcal{F}$. Consider an ultrafilter $\mathcal{H}$ in X containing the family $\{F\cap V:F\in \mathcal{F},$ V is a neighborhood of $\xi$ in $s_\mu X\}$. Again, by Lemma \ref{ultrafilter},  $\mathcal{H}$ is Cauchy in $(X,s_f\mu)$ because it is Bourbaki-Cauchy in $(X,\mu)$ as it contains $\mathcal{F}$.  Since $\xi$ is the only cluster point of it, $\xi$ must be a point in $\widetilde{X}$.

$2)$ This follows at once by Lemma \ref{ultrafilter}.
\end{proof}
\end{theorem} 

\begin{theorem}\label{star-finite2} Let $(X,\mu)$ be a uniform space. The following statements are satisfied:
\begin{enumerate}

\item $(X,\mu)$ is Bourbaki-complete if and only if $(X, s_f \mu)$ is complete. 

\item $(X,\mu)$ is cofinally Bourbaki-complete if and only if $(X, s_ f \mu)$ is cofinally complete. 

\end{enumerate}

\noindent In particular, every uniform space with a base of star-finite covers for the uniformity satisfies that it is (cofinally) complete if and only if it is (cofinally) Bourbaki-complete.

\end{theorem}

Next, recall the following classical result.

\begin{theorem} A uniform space $(X,\mu)$ is compact if and only if it is totally bounded and complete

\end{theorem}

\noindent In parallel we have the next theorem which follows easily by Theorem \ref{BB.sf} and Theorem \ref{star-finite2}. 

\begin{theorem} \label{Bourbaki-compact} The following statements are equivalent for a uniform space $(X,\mu)$:
\begin{enumerate}
\item $X$ is compact;

\item $(X,\mu)$ is Bourbaki-bounded and Bourbaki-complete;

\item $(X,\mu)$ is totally bounded and complete;

\item $(X,s_f\mu)$ is Bourbaki-bounded and Bourbaki-complete;

\item $(X,s_f\mu)$ is totally bounded and complete.
\end{enumerate}

\end{theorem}

\smallskip

In the next example we show that in general
\begin{center} 

{\it complete} $\nRightarrow$ {\it Bourbaki-complete }
\end{center}

\begin{center}
{\it cofinally complete} $\nRightarrow$ {\it Bourbaki-complete.}

\end{center}

\begin{example} \label{hedgehog}{\it There exists a Bourbaki-bounded cofinally complete metric space which is not totally bounded and not Bourbaki-complete: the metric hedgehog.} 

\begin{proof} Let $A$ be a set of cardinal $\kappa \geq \omega _0$. The metric hedgehog $H(\kappa)$ of  $\kappa$ spininess is defined as follows (\cite{engelkingbook}). Let $I=[0,1]\subset  \Rset$ and consider the product $I\times A$. Next take the equivalence relation $\sim$ on $I\times A$ defined as follows: $$(x_1, \alpha_1)\sim (x_2, \alpha_2)\text{ if and only if } x_1=x_2=0.$$ Then  the quotient $I\times A \big{/ \sim}$ is the set of points of $H(\kappa)$ and we endowed this set with the following metric $\rho:H(\kappa)\times H(\kappa)\rightarrow [0,\infty)$ $$\rho([(x_1,\alpha_1)], [(x_2,\alpha_2)])=\begin{cases}  \,\,  |x_1 -x_2| &     \text {if $\alpha_1= \alpha_2$;}  \\ \,\,\, x_1 +x_2  & \text{ if $ \alpha_1 \neq \alpha_2$.}\end{cases}$$  It is easy to prove that the metric space $(H(\kappa), \rho)$ is cofinally complete (and then complete). For that, recall the metric characterization of cofinally completeness in \cite[Theorem 2.1.1]{hohti-thesis} and in \cite[Theorem 3.2]{beer-between1} (we will see it later).  

Moreover, $(H(\kappa), \rho)$  is a Bourbaki-bounded space since $$H(\kappa)\subset B^{n+1}_{\rho}(0, 1/n)$$ for every $n\in \Nset$, where, by abuse of notation,  $0$ denotes also the point $[(0,\alpha)]\in H(\omega_1)$. But  $H(\omega_1)$ is not compact because it contains a uniformly discrete subspace of infinite cardinality $\kappa$, precisely the set of points $\{[(1, \alpha)]: \alpha\in A\}$. Then it cannot be totally bounded, nor Bourbaki-complete, by Theorem \ref{Bourbaki-compact}.  
\end{proof}

\end{example}

The next example, which is a particular case of the previous one, shows that for a uniform space $(X,\mu)$ 
\begin{center} 

{\it completeness of $(X, e\mu)$} $\nRightarrow$ {\it completeness of $(X, s_f\mu)$}
\end{center}

\begin{center}
{\it completeness of $(X, s_c\mu)$} $\nRightarrow$ {\it completeness of $(X, s_f\mu)$},
\end{center}
\noindent even if the star-countable modification $s_c\mu$ is a compatible uniformity.

\begin{example}{\it There exists a separable complete metric space $(X,d)$ which is not Bourbaki-complete. In particular $(X, s_c \mu _d)$ and $(X, e\mu _d)$ are complete.} 

\begin{proof}
Consider a separable complete metric space which is not Bourbaki-complete, like the separable metric hedgehog $(H(\omega _0), \rho)$. Then the uniformity $\mu _\rho$ generated by the metric has a base of countable covers by the Lindel\"of property, that is, $\mu_\rho=e\mu _\rho$. It is clear then that  $\mu _\rho=s_c \mu_\rho =e\mu _\rho$. By completeness of $(H(\omega _0),\rho)$, then $(H(\omega _0),  s_c \mu_\rho)$ and $(H(\omega _0),  e \mu_\rho)$ are also complete. However, $(H(\omega _0), s_f \mu _\rho)$ is not complete since the space is not Bourbaki-complete. 
\end{proof}
\end{example}

It is clear that in Theorem \ref{Bourbaki-compact} one can change Bourbaki-complete by cofinally Bourbaki-complete or complete by cofinally complete. Now, we are going to refine this result since``cofinal complete-like properties'' have the characteristic of transforming local properties to uniformly local properties.  

\begin{definition} Let $P$ denote the property  of being a compact, a totally bounded or a Bourbaki-bounded subset in a uniform space $(X,\mu)$. A uniform space $(X,\mu)$ is {\it  locally P} if for every $x\in X$ there exists a neighborhood $V^x$ of $x$ such that $V^x$ satisfies the property $P$.  It is said {\it uniformly locally P} if there exists some $\mathcal{U}\in \mu$ such that every $U\in \mathcal{U}$ satisfies $P$.
\end{definition}

From now on, recall that the the notation $\mathcal{U}< \mathcal{V}$ means that the cover $\mathcal{U}$ refines the cover $\mathcal{V}$. In addition, whenever $\mathcal{C}$ is a cover of a set $X$, we will denote by $\mathcal{C}^{*}$ the cover $\{St(C,\mathcal{C}):C\in \mathcal{C} \}$. 

Next, a family of sets $\mathcal{L}$ is {\it directed} provided that, for all $L,M\in \mathcal{L}$, there exists $N\in \mathcal{L}$ such that $L\cup M\subset N$. Note that for every cover $\mathcal{A}$ of $X$, the family $\mathcal{A}^{f}=\{\bigcup \mathcal{E}: \mathcal{E}\subset \mathcal{A} \text{ and } \mathcal{E} \text{ is finite}\}$ is a directed cover. Observe that for a directed cover $\mathcal{C}$, satisfying that $X\notin \mathcal{C}$, the family of sets $\{X\backslash C: C\in \mathcal{C}\}$ is a filter base. Conversely, if $\mathcal{F}$ is a filter base then the family of sets $\{X\backslash F: F\in \mathcal{F}\}$ is a directed cover (\cite{howesbook}).

\begin{theorem} {\rm (\cite{rice})} \label{rice}A uniform space $(X,\mu)$ is uniformly locally compact if and only if it is cofinally complete and locally totally bounded (locally compact).
\begin{proof} $\Rightarrow)$ We just need to prove cofinal completeness. Let $\mathcal{F}$ be a cofinally Cauchy filter and fix $\mathcal{U}\in \mu$ such that ${\rm cl}_{X} U$ is compact for every $U\in \mathcal{U}$. Then for some $U_0\in \mathcal{U}$, $U_0\cap F\neq \emptyset$ for every $F\in \mathcal{F}$. Let $\mathcal{H}$ be an ultrafilter containing  $\mathcal{F}\cup \{U_0\}$. Then, by compactness of ${\rm cl}_{X} U_0$,  $$\emptyset\neq \bigcap_{H\in \mathcal{H}} {\rm cl}_{X} H \subset {\rm cl}_{X} U_0$$ that is, $\mathcal{H}$ converges and in particular $\mathcal{F}$ clusters.

$\Leftarrow)$ Let $(X,\mu)$ be cofinally complete and  locally compact. Without loss of generality we can assume that $X$ is not compact. By local compactness, for every $x\in X$ there exists an open neighborhood $V^x$ of $x$ in  $X$  such that  ${\rm cl}_X V^x$ is compact. Put $\mathcal{V}=\{V^x:x\in X\}$ and consider the directed open cover $\mathcal{V}^f$ given by all the finite unions of elements from $\mathcal{V}$. Then, the family of sets $\{X\backslash W:W\in \mathcal{V}^f\}$ is a filter base for a filter $\mathcal{F}$ in $X$ which does not cluster. Therefore, $\mathcal{F}$ cannot be cofinally Cauchy, that is, there exists some $\mathcal{U}\in \mu$ such that for every $U\in \mathcal{U}$, $U\cap F=\emptyset$ for some $F\in \mathcal{F}$. But this is equivalent to say that $\mathcal{U}<\mathcal{V}^f$. Now, for every $U\in \mathcal{U}$ fix $V^{x_1}\cup V^{x_2}\cup...\cup V^{x_k}\in \mathcal{V}^f$ such that $U\subset V^{x_1}\cup V^{x_2}\cup...\cup V^{x_k}$. Then, $${\rm cl}_X U\subset {\rm cl}_X (V^{x_1}\cup V^{x_2}\cup...\cup V^{x_k})={\rm cl}_X  V^{x_1}\cup {\rm cl}_X  V^{x_2}\cup...\cup {\rm cl}_X  V^{x_k}.$$ Thus, ${\rm cl}_X U$ is compact for every $U\in \mathcal{U}$, that is, $(X,\mu)$ is uniformly locally compact.

\end{proof}
\end{theorem}

\begin{lemma} \label{unif.loc.compact}{\rm  (\cite[Lemma 1.17]{hejcman})} Let $(X,\mu)$ be a uniformly locally compact space. Then there exists $\mathcal{V}\in \mu$ such that ${\rm cl}_X St^m (V,\mathcal{V})$ is compact for every $V\in \mathcal{V}$ and every $m\in \Nset$.
\begin{proof} Let $\mathcal{U}\in \mu$ such that for every $U\in \mathcal{U}$, ${\rm cl}_X U$ is compact. Let $\mathcal{V}\in \mu$ open cover such that $\mathcal{V}^{*}<\mathcal{U}$. Then, ${\rm cl}_X St(V,\mathcal{V})$ is compact for every $V\in \mathcal{V}$, because ${\rm cl}_X St(V,\mathcal{V}) \subset {\rm cl}_X U$. Next, take $K\subset X$ a compact subset, then there exists a finite subfamily $\{V_i:i=1,...,m\}\subset\mathcal{V}$ such that $K\subset \bigcup_{i=1}^{m} V_i$. Thus, ${\rm cl}_{X}St(K,\mathcal{V}) \subset {\rm cl}_{X} \bigcup_{i=1}^{m} St(V_i,\mathcal{V})=  \bigcup_{i=1}^{m} {\rm cl}_{X} St(V_i,\mathcal{V}) $ and then ${\rm cl}_{X}St(K,\mathcal{V}) $ is also compact. In particular, by induction, ${\rm cl}_X St^m (V,\mathcal{V})$ is compact for every $m\in \Nset$.
\end{proof}
\end{lemma}

\begin{theorem} \label{unif.local.compact.th}For a uniform space $(X,\mu)$ the following statements are equivalent:
\begin{enumerate}
\item $(X,\mu)$ is uniformly locally compact;

\item $(X,\mu)$ is  locally Bourbaki-bounded in $X$ and cofinally Bourbaki-complete;

\item $(X,\mu)$ is locally totally bounded and cofinally complete;

\item $(X,s_f\mu)$ is uniformly locally compact;

\item $(X,s_f\mu)$ is  locally Bourbaki-bounded in $X$ and cofinally Bourbaki-complete;

\item $(X,s_f\mu)$ is locally totally bounded and cofinally complete.

\end{enumerate}
\begin{proof} The equivalences $(1)\Leftrightarrow (3)$ and $(4)\Leftrightarrow (6)$ follows from Theorem \ref{rice}. Next, $(2)$, $(5)$ and $(6)$ are equivalent by Theorem \ref{BB.sf} and Theorem \ref{star-finite2}. The implication $(4)\Rightarrow (1)$ is trivial. 

Finally, we prove $(1)\Rightarrow (2)$. We just need to check cofinal  Bourbaki-completeness. Applying Lemma \ref{unif.loc.compact}, let  $\mathcal{V}\in \mu$ such that  ${\rm cl}_X St^m (V,\mathcal{V})$ is compact for every $V\in \mathcal{V}$ and every $m\in \Nset$. Now, let $\mathcal{F}$ be a cofinally Bourbaki-Cauchy filter. Then for some $V\in \mathcal{V}$ and $m\in \Nset$, $St^m (V,\mathcal{V})\cap F\neq \emptyset$ for every $F\in \mathcal{F}$. Let $\mathcal{H}$ be an ultrafilter containing $\mathcal{F}\cup \{St^m (V,\mathcal{V})\}$. Then, by compactness of ${\rm cl}_{X} St^m (V,\mathcal{V})$,  $$\emptyset\neq \bigcap_{H\in \mathcal{H}} {\rm cl}_{X} H \subset {\rm cl}_{X} St^m (V,\mathcal{V}),$$ that is, $\mathcal{H}$ converges and in particular $\mathcal{F}$ clusters.

\end{proof}
\end{theorem}

\begin{remark} In the previous result it is  implicit in the proof that {\it every Bourbaki-complete uniform space satisfies that every closed Bourbaki-bounded subsets of it is compact}. We will come back later to this subject. Moreover, it is easy to check that {\it Bourbaki-completeness and cofinal Bourbaki-completeness are properties inherited by closed subspaces}. 
\end{remark}

The next example shows us that in general

\begin{center} 

{\it completeness } $\nRightarrow$ {\it cofinal completeness}

\end{center}
\begin{center} 

{\it Bourbaki-completeness } $\nRightarrow$ {\it cofinal completeness}

\end{center}

\begin{center} 

{\it Bourbaki-completeness } $\nRightarrow$ {\it cofinal Bourbaki-completeness}.

\end{center}

\begin{example} \label{Example.discrete}(\cite[Example 15]{merono.completeness}) {\it There exists a discrete countable metric space, hence separable and locally compact, which is Bourbaki-complete and not uniformly locally compact, nor cofinally Bourbaki-complete, nor cofinally complete.} 

\begin{proof}[Construction.] This will be a subspace of the Banach space $(\ell _\infty, ||\cdot ||_{\infty})$ of all the bounded sequences of $\Rset$ with the supremum norm. Indeed, let $\{e_n:n\in \Nset\}$ the canonical base of $\ell _\infty$ and put $X=\bigcup _{n\in \Nset}A_n $ where $A_{n}=\{e_n\}\cup\{e_n+\frac{1}{n}e_k:k\in \Nset\}$. Consider on $X$ the metric $d$ inherited from $\ell_\infty$. Then $(X,d)$ is Bourbaki-complete since for every Bourbaki-Cauchy filter there is a singleton belonging to it, but it is not uniformly locally compact since for every $n\in \Nset$ the closed ball of radius $\frac{1}{n}$ and center $e_{n}$ is infinite and discrete, and hence not compact. 
\end{proof}
\end{example}

We close this section with an example of a complete discrete metric space which is not Bourbaki-complete. Observe that, in general, every uniformly discrete metric space is Bourbaki-complete, since they are uniformly locally compact. Therefore, the next example cannot be uniformly discrete.

\begin{example} (\cite[Example 4.1]{merono.hh}) {\it There exists a countable discrete metric space, hence separable and locally compact, which is complete, Bourbaki-bounded and not Bourbaki-complete, nor uniformly locally compact, nor cofinally Bourbaki-complete, nor cofinally complete}.  

\begin{proof}[Construction] As in the previous example, this space will be a subspace of the space $(\ell _\infty,||\cdot||_{\infty})$.  Let again $\{e_n:n\in \Nset\}$ be the canonical basis of $\ell _\infty$ and $d$ the metric on $\ell_{\infty}$ generated by the norm. Then $d(e_n, e_k)=1$ for every $n,k\in \Nset$, $n\neq k$. Let $n,k\in \Nset$, $n<k$. For every $i\in \{0\}\cup \Nset$, $0\leq i \leq 2^k$, set $$x_{n,k,i}=e_n+i\cdot 2^{-k}(e_k -e_n).$$ We have that $x_{n,k,0}=e_n$ and $x_{n,k,2^k}=e_k$ and hence the points $x_{n,k,i}$ form a chain $x_{n,k,0}, x_{n,k,1},...$, $x_{n,k,2^k}$. Note that $d(x_{n,k,i}, x_{n,k,i+1})=2^{-k}$ for each $0\leq i<2^k$. Define $$X=\{e_n:n\in \Nset\}\cup \{x_{n,k,i}: n,k,i \in \Nset, n<k \text{ and } 2^{k-n}\leq i<2^k\}.$$ Informally, $X$ consists on the points $e_1, e_2,...$ and some of the points in the segments connecting $e_n$, and $e_k$, for every $n,k\in \Nset$, $n<k$. More precisely, we have let out some points $x_{n,k,i}$ so that now in our chain from $e_n$ to $e_k$, there is a first ``jump'' of length $2^{-n}$ from $e_n$ to $x_{n,k,2^{n-k}}$, and after that, the distance between two consecutive points of the chain is $2^{-k}$. Let us endowed $X$ with the metric $d$ inherited from $\ell _{\infty}$. It is easy to see that $X$ is closed in $\ell _{\infty}$. As a consequence $(X,d)$ is complete. Moreover, it is also a discrete space because $B_{d}(e_{n}, 2^{-n})=\{e_{n}\}$ and $B_{d}(x_{n,k,i}, 2^{-k})=\{x_{n,k,i}\}$ for all $x_{n,k,i}\in X$. 

Now, we prove that $(X, d)$ is Bourbaki-bounded and hence it cannot be Bourbaki-complete because otherwise, it would be compact by Theorem \ref{Bourbaki-compact}. Let $\varepsilon >0$ and $h\in \Nset$ such that $2^{-h}<\varepsilon$. Put $A=\{x_{n,k,i}\in X: k\leq h\}$. Then $A$ is finite and $e_1,...,e_h \in A$. Let $m= 2^{h+1}$. We will prove that $\bigcup \{St^{m}(a, \varepsilon):a \in A\}=X$ by showing that $X\backslash A\subset St^{m}(e_h,\varepsilon)$. Let $x_{n,k,i}\in X\backslash A$. then we have that $k>h$. The chain $$x_{h,k,0}, x_{h,k,2^{k-h}}, x_{h,k,2\cdot 2^{k-h}}, ..., x_{h,k,2^{h}\cdot 2^{k-h}}$$ joins $e_h$ to $e_k$, and for every $0\leq i \leq2^{h}$, we have $$d(x_{h,k,i\cdot 2^{k-h}}, x_{h,k,(i+1)\cdot 2^{k-h}})=2^{k-h}\cdot 2^{-k}=2^{-h}<\varepsilon.$$ To find a similar chain joining $e_k$ to $x_{n,k,i}$, let $j\in \Nset$ such that $(j-1)\cdot 2^{-h}\leq i \cdot 2^{-k}\leq j\cdot 2^{-h}$. Then, the chain $$x_{n,k,2^{h}\cdot 2^{k-h}},x_{n,k,2^{h-1}\cdot 2^{k-h}},...,x_{n,k,j\cdot 2^{k-h}},x_{n,k,i}$$ joins $e_k$ to $x_{n,k,i}$ and consecutive terms of this chain have distance at most $2^{-h}$. The two chains constructed above verify that $e_k\in B_{d} ^{2^ h}(e_h, \varepsilon)$ and $x_{n,k,i} \in B_{d} ^{2^ h}(e_k, \varepsilon)$. It follows that  $x_{n,k,i} \in B_{d} ^{2^ {h+1}}(e_h, \varepsilon)=B_{d} ^m (e_h, \varepsilon)$ as we advanced. 
\end{proof} 
\end{example}

\medskip

\subsection{Uniform strong paracompactness and related topological properties}

\hspace{15pt} In the introduction we have talked about the following uniform version of paracompactness due to Rice \cite{rice}.

\begin{definition} A uniform space $(X,\mu)$ is {\it uniformly paracompact} is every open cover $\mathcal{G}$ has an open refinement $\mathcal{A}$ which is {\it uniformly locally finite}, that is, there exists some $\mathcal{U}\in \mu$ such that every $U\in \mu$ meets at most finitely many $A\in \mathcal{A}$.
\end{definition}

The main result about uniform paracompactness is its equivalence to cofinal completeness.

\begin{theorem} \label{cof.complete=unif.paracompact} {\rm (\cite{howesbook})} A uniform space $(X,\mu)$ is cofinally complete if and only if it is uniformly paracompact. 
\end{theorem}

\noindent As a corollary we have  Corson's Theorem \ref{corson}: {\it  A space $X$ is paracompact if and only if $(X,{\tt u})$ is uniformly paracompact, or equivalently, $(X,{\tt u})$ is cofinally complete}.

Resembling the above result, in \cite{merono.completeness} and \cite{merono.nets} it is proved that cofinal Bourbaki-completeness is equivalent to {\it uniform strong paracompactness}, a uniform property  introduced by Hohti in \cite{hohti-thesis}. It is a uniform extension of the topological property of {\it strong paracompactness}.

\begin{definition} A space is {\it strongly paracompact} if every open cover has an open star-finite refinement.
\end{definition}

\begin{definition}{\rm (\cite{hohti-thesis})} A uniform space $(X,\mu)$ is {\it uniformly strongly paracompact} if every open cover $X$ has a uniformly star-finite open refinement, where a cover $\mathcal{A}$ is {\it uniformly star-finite} if there exists $\mathcal{U}\in \mu$ such that for every $A\in\mathcal{A}$, $St(A,\mathcal{U})$ meets at most finitely many $A'\in \mathcal{A}$.
\end{definition}

\noindent Observe that in the definition of strong paracompactness we can change star-finite by star-countable. Indeed, we have previously noticed in Remark \ref{star-countable-cozero}  that every star-countable cozero cover has a star-finite cozero refinement, and then applying Lemma \ref{shrinking} the fact follows.

\smallskip

Here we change the proofs from \cite{merono.completeness} and \cite{merono.nets} and we start by giving a technical theorem which provide us of a characterization of uniformitites with a star-finite base. This is a more tangible way to express it as (cofinal) Bourbaki-Cauchy filters of a uniform space $(X,\mu)$ is more tangible than (cofinal) Cauchy filters of $(X,s_f\mu)$.

\smallskip

\begin{theorem} \label{star} A uniform space $(X,\mu)$ has a star-finite base for its uniformity, that is, $\mu=s_f\mu$, if and only if it satisfies the following property:
\smallskip

$(\star)$  for every $\mathcal{U}\in \mu$ there is some $\mathcal{V}\in \mu$ satisfying that for every $V\in \mathcal{V}$ and every $n\in \Nset$ there exist finitely many $U_1,...,U_k \in \mathcal{U}$ such that $$St^{n}(V,\mathcal{V})\subset \bigcup _{i=1}^{k}U_i$$
\begin{proof} $\Rightarrow)$ Let $\mathcal{U}\in \mu$ and $\mathcal{V}\in \mu$ star-finite such that $\mathcal{V}<\mathcal{U}$. By the star-finite property, for every $V\in \mathcal{V}$ and every $n\in \Nset$ there exist at most finitely many $V'\in \mathcal{V}$ such that $V'\cap St^{n}(V,\mathcal{V})\neq \emptyset$. Since $\mathcal{V}$ is a refinement of $\mathcal{U}$ the $(\star)$-property follows. 

$\Leftarrow)$ Conversely,  let $\mathcal{U}\in \mu$ and select $\mathcal{V}\in \mu$ such that $\mathcal{V}^{*}<\mathcal{U}$. By hypothesis there is some $\mathcal{W}\in \mu$ such that for every $n\in \Nset$ and every $W\in \mathcal{W}$ there exists finitely many $V_i\in \mathcal{V}$, $i=1,...,k$ such that $St^{n}(W,\mathcal{W})\subset \bigcup _{i=1}^k V_i$. Without loss of generality we can take $\mathcal{W}$ refining $\mathcal{V}$.  
Let $\mathcal{A}(\mathcal{W})=\{A_n ^i, n\in \Nset, i\in I\}$ from Lemma \ref{corona} induced by $\mathcal{W}$. Clearly $\mathcal{A}(\mathcal{W})\in s_f \mu$. By hypothesis, for every $i\in I$ and every $n\in \Nset $ we can fix a finite family $\mathcal{V}_{i,n}\subset \mathcal{V}$ such that $A_n ^i\subset St^{n+1}(x_i,\mathcal{W})\subset St^{n+1} (W_i,\mathcal{W})\subset\bigcup \{V:V\in\mathcal{V}_{i,n}\}$ where $W_i\in \mathcal{W}$ is some set such that $W_i \subset St(x_i,\mathcal{W})$.  Define $$\mathcal{G}= \{A_n^i\cap St(V,\mathcal{V}): V\in \mathcal{V}_{i,n}, i\in I, n\in \Nset\} .$$ Then $\mathcal{A}(\mathcal{W})\wedge \mathcal{V}=\{A\cap V: A\in \mathcal{A}(\mathcal{W}), V\in \mathcal{V}\}<\mathcal{G}<\mathcal{V}^{*}<\mathcal{U}$, and it is easy to check that $\mathcal{G}$ is also star-finite.
\end{proof}

\end{theorem}

\begin{theorem}\label{uniform.strongly.paracompact} Let $(X,\mu)$ be a uniform space. The following statements are equivalent:
\begin{enumerate} 
\item $(X,\mu)$ is uniformly strongly paracompact.

\item $(X,\mu)$ uniformly paracompact and $\mu=s_f\mu$.

\item $(X,\mu)$ is cofinally complete and $\mu=s_f\mu$.

\item $(X,\mu)$ cofinally Bourbaki-complete.
\end{enumerate}
\begin{proof}$(1)\Rightarrow(2)$ That uniform strong paracompactness implies uniform paracompactness is clear from the definitions. In order to see that $\mu=s_f\mu$ we are going to prove that $(X,\mu)$ satisfies the $(\star)$-property from Theorem \ref{star}.

Take $\mathcal{U}\in \mu$ and $\mathcal{W}\in \mu$ an open cover such that $\mathcal{W}<\mathcal{U}$. By uniform strong paracompactness there exists an open refinement $\mathcal{A}$ of $\mathcal{W}$ and there exists $\mathcal{V}\in \mu$ such that for every $A\in \mathcal{A}$, $St(A,\mathcal{V})$ meets at most finitely many $A'\in \mathcal{A}$. By induction it is not difficult to see that for every $V\in \mathcal{V}$ and every $n\in \Nset$, $St^{n}(V,\mathcal{V})$ is covered by finitely many elements from $\mathcal{A}$. Indeed, let $V\in \mathcal{V}$ and $A\in \mathcal{A}$ such that $A\cap V\neq \emptyset$. Then $V\subset St(A,\mathcal{V})$ and  there exists finitely many $A_1,...,A_k\in \mathcal{A}$ such that $V\subset St(A,\mathcal{V})\subset \bigcup_{i=1}^{k} A_i$. Suppose that for some $n\in \Nset$, $St^{n}(V,\mathcal{V})$ is covered by finitely $A_1,...,A_k\in \mathcal{A}$ then $$St^{n+1}(V,\mathcal{V})\subset \bigcup _{i=1}^{k}St(A_i,\mathcal{V}).$$ 	Since each set $St(A_i,\mathcal{V})$ is covered by finitely many elements from $\mathcal{A}$ the result follows. Moreover, as $\mathcal{A}$ is a refinement of $\mathcal{W}$ and $\mathcal{W}$ is a refinement of $\mathcal{U}$ the $(\star)$-property follows.

$(2)\Rightarrow (1)$ Suppose that $\mu=s_f\mu$ and that $(X,\mu)$ is uniformly paracompact. Let $\mathcal{G}$ be an open cover of $X$ and $\mathcal{A}$ a uniformly locally finite open refinement of $\mathcal{G}$. Then there exists some $\mathcal{U}\in \mu$ such that every $U\in \mathcal{U}$ meets at most finitely many $A\in \mathcal{A}$. Let $\mathcal{V}\in \mu$ being star-finite such that $\mathcal{V}<\mathcal{U}$. Observe that we can take $\mathcal{V}$ being open. Next, define $\mathcal{W}=\{V\cap A: V\in \mathcal{V},A\in \mathcal{A}\}$. Then $\mathcal{W}<\mathcal{G}$ and it is an open cover. Moreover, for every $W\in \mathcal{W}$, $St(W,\mathcal{V})$ meets only fintely many $W'\in \mathcal{W}$. Indeed, $\mathcal{V}$ is star-finite and each $V\in \mathcal{V}$ meets at most finitely many $A\in \mathcal{A}$ as $\mathcal{V}<\mathcal{U}$.

$(2)\Leftrightarrow (3)$ This equivalence follows from the equivalence of cofinal completeness and strong paracompactness (Theorem \ref{cof.complete=unif.paracompact}).

$(3)\Rightarrow (4)$ This implication follows from Theorem \ref{star-finite2}.

$(4)\Rightarrow (3)$ Since $(X,\mu)$ is cofinally Bourbaki-complete then it is cofinally complete. On the other hand, suppose by contradiction that it does not satisfies the $(\star)$-property (Theorem \ref{star}). Then there exists some $\mathcal{U}_0\in \mu$, that we can take open, such that for every $\mathcal{V}\in \mu$ there exists $V_0 \in \mathcal{V}$ and $m_0\in \Nset$ for which there is no finite subfamily in $\mathcal{U}_0$ covering $St^{m_0}(V_0, \mathcal{V})$. Let $\mathcal{U} _0 ^{f}$ the cover obtained by taking finite unions of elements of $\mathcal{U}_0$. Then $\mathcal{U} _0 ^{f}$ is a directed open cover of $X$ and $\{X\backslash A:A\in \mathcal{U}^{f}\}$ is a filter base of a filter $\mathcal{F}$ in $X$ (note that $X\notin \mathcal{U}_0 ^{f}$). In particular, $\mathcal{F}$ is cofinally Bourbaki-Cauchy since for every $\mathcal{V}\in \mu$ there exists $V_0\in \mathcal{V}$ such that $F\cap St^{m_0}(V_0,\mathcal{V})\neq \emptyset$ for every $F\in \mathcal{F}$. Therefore, $\mathcal{F}$ clusters contradicting that $\mathcal{U}_0 ^{f}$ is a cover.
\end{proof}
\end{theorem}

The above theorem states that every uniformly strongly paracompact uniform space (equivalently, cofinally Bourbaki-complete uniform space) has a star-finite base for its uniformity. Observe that similarly, Hohti  \cite[pp. 31-32]{hohti-thesis} proved that every cofinally complete uniform space has a point-finite base.

\begin{theorem} {\rm (\cite{hohti-thesis})} \label{point.finite} If a uniform space $(X,\mu)$  is cofinally complete, then $\mu=p_f\mu$. In particular, $(X,\mu)$  is cofinally complete if and only if $(X,p_f\mu)$ is cofinally complete.
\end{theorem}

Next, we apply Theorem \ref{uniform.strongly.paracompact} to the fine uniformity ${\tt u}$ obtaining a result in the line of Corson's Theorem \ref{corson}.

\begin{corollary} \label{uniformly.strong.fine} For a  space $X$ the following statements are equivalent:
\begin{enumerate}
\item $X$ is strongly paracompact;

\item $X$ is paracompact and ${\tt u}=s_f{\tt u}$

\item $(X,{s_f\tt u})$ is uniformly paracompact (equivalently, cofinally complete) and ${\tt u}=s_f{\tt u}$

\item $(X,{\tt u})$ is uniformly strongly paracompact (equivalently, cofinally Bourbaki-complete).
\end{enumerate}
\end{corollary} 

By the above result cofinal Bourbaki-completeness of $(X,{\tt u})$ is just strong paracompactness of the topological space $X$. Thus, cofinal completeness of $(X,{\tt u})$, that is, paracompactness of $X$, is a weaker property. We ask now, if Bourbaki-completeness of  $(X,{\tt u})$ is also stronger than completeness of $(X,{\tt u})$. To that purpose recall the following facts.

\begin{definition} A space $X$ is {\it topological complete} if it is uniformizable by a uniformity $\mu$ such that $(X,\mu)$ is complete.
\end{definition}

\noindent The following result is immediate \cite{engelkingbook}.

\begin{theorem} A space $X$ is topologically complete if and only if $(X, {\tt u})$ is complete.
\end{theorem}

\begin{definition} {\rm (\cite{garcia.delta-complete})} A  space $X$ is $\delta$-{\it complete} if $(X, s_f{\tt u})$ is complete.
\end{definition}

\begin{theorem} \label{delta}{\rm (\cite{merono.completeness})} Let $X$ be a  space. The following statements are equivalent:
\begin{enumerate}
\item $X$ is $\delta$-complete, that is, $(X, s_f{\tt u})$ is complete

\item $(X, s_f{\tt u})$ is Bourbaki-complete;

\item $(X, {\tt u})$ is Bourbaki-complete.
\end{enumerate} 

\begin{proof} The proof follows from Theorem \ref{star-finite2}. 
\end{proof}
\end{theorem}

Now, we apply the above results to the countable modification  $e{\tt u}$. But first, we recall the definition of realcompact space.

\begin{definition} A space $X$ is realcompact if and only if $X$ is homeomorphic to a closed subspace of a product of real-lines.
\end{definition}

\noindent We will come back to realcompactness in the next part of thesis. Nevertheless we will recall now Shirota's Theorem which is well-known.

\begin{theorem}\label{shirota}{\rm (\cite{shirota1})} A space $X$ is  realcompact if and only if $(X,e{\tt u})$ is complete.
\end{theorem}

\begin{corollary} \label{delta1} A space $X$ is realcompact if and only if $(X,e{\tt u})$ is Bourbaki-complete.
\begin{proof}  By Corollary \ref{morita2},  $s_f (e{\tt u})=e{\tt u}$. Therefore the proof is immediate from Theorem \ref{star-finite2}.
\end{proof}
\end{corollary}

It is clear that {\it every realcompact space is $\delta$-complete} because $e{\tt u}\leq s_f\tt{u}$ (Lemma \ref{morita}). However, not every $\delta$-complete space is realcompact as it is show in the next example.

\begin{example}\label{discrete.measurable} {\it There exists a uniform space which is $\delta$-complete but not realcompact.}
\begin{proof} Let $(D,\chi)$ be a uniformly discrete metric space of Ulam-measurable cardinal. As we have said in the introduction, this is equivalent to have $(D,\chi)$ failing realcompactness, that is, $(D, e{\tt u})$ is not complete. However, $(D,\chi)$ is Bourbaki-complete as every Bourbaki-bounded subset in $D$ is a finite set. Since $\mu _\chi=s_f\mu _\chi \leq s_f {\tt u}$, $(D, s_f{\tt u})$ is complete, that is, $D$ is $\delta$-complete.
\end{proof}
\end{example}

On the other hand, it is immediate that {\it every $\delta$-complete space is topologically complete}. However, to give a counterexample that the reverse implication is not true, is not so clear. Observe that by Lemma \ref{cozero.connected} the following result is immediate.

\begin{theorem}\label{connected.realcompact} A connected space is realcompact if and only if it is $\delta$-complete.
\end{theorem}

\noindent Therefore, by the previous result, a connected topological complete space which is not realcompact is an example of topological complete space space which is not $\delta$-complete. 

\begin{example}(\cite{garcia.equality}, \cite{comfort})\label{hedgehog.measurable}{\it There exists a topological complete space which is not $\delta$-complete.}
\begin{proof}The metric hedgehog $H(\kappa)$ (Example \ref{hedgehog}) where $\kappa$ is an Ulam-measurable cardinal is a complete metric space and hence it is topological complete. That $H(\kappa)$ is not realcompact follows from the fact that realcompactness is inherited by closed subspaces. Nevertheless, observe that the uniformly discrete subspace $\{[(1,\alpha)]:\alpha <\kappa\}$ is not realcompact as its cardinality is exactly $\kappa$. By Theorem \ref{connected.realcompact}, it is not $\delta$-complete either.
\end{proof}
\end{example}

Finally,  recall that in \cite{howes.completeness} Howes proved the following result.

\begin{theorem} \label{howes} A space $X$ is Lindel\" of if and only if $(X,e{\tt u})$ is cofinally complete.
\end{theorem}

In particular, similarly to Corollary \ref{delta1}, we have the next result.

\begin{corollary}A space $X$ is Lindel\" of if and only if $(X,e{\tt u})$ is cofinally Bourbaki-complete.
\end{corollary}

\noindent Recall that {\it every Lindel\"of space is strongly paracompact as $s_f{\tt u}\geq e{\tt u}$} (Lemma \ref{morita}). In particular, from Lemma \ref{cozero.connected} we can deduce the following result of Morita (see \cite{morita.star-finite}). 

\begin{theorem}  \label{strongly.paracompact.morita} A connected space is strongly paracompact if and only if it is Lindel\"of.
\end{theorem}

Not every paracompact space is strongly paracompact and not every strongly paracompact space is Lindel\"of as we show in the next examples.

\begin{example}{ \it There exists a paracompact space which is not strongly paracompact.}
\begin{proof} By theorem  \ref{strongly.paracompact.morita}, the metric hedgehog $H(\omega _1)$ (Example \ref{hedgehog}) is not strongly paracompact because it is not Lindel\"of. However it is paracompact because it is a metric space.
\end{proof}
\end{example}

\begin{example}\label{baire}{ \it There exists a strongly paracompact space which is not Lindel\"of.}
\begin{proof} Every strongly-zero dimensional is strongly paracompact since every open cover has an open refinement which is a clopen partition of $X$ (see \cite{engelkingbook}). Therefore any strongly zero-dimensional which is not Lindel\"of is an example: for instance the countable product $D^{\omega _0}$ where $D$ is a discrete space of uncountable cardinal.
\end{proof}
\end{example}

We close this section with the next example that shows that, in general, for a uniform space $(X,\mu)$
\begin{center}
 {\it $\delta$-complete} $\nRightarrow$ {\it Bourbaki-complete}
\end{center}

\begin{center}
 {\it strongly paracompact} $\nRightarrow$ {\it cofinally Bourbaki-complete} 
\end{center}
as it was expected.

\begin{example}{\it There exists a strongly paracompact, in particular $\delta$-complete, metric space which is not cofinally Bourbaki-complete, nor Bourbaki-complete, for the metric uniformity. } 

\begin{proof} Let $(X,d)$ be a separable non-complete metric spaces, for instance, the open interval $(0,1)$ endowed with the euclidean metric. Then $X$ is strongly paracompact. However it is not Bourbaki-complete for the metric uniformity, nor cofinally Bourbaki-complete
\end{proof}
\end{example}

\medskip

\subsection{Products and hyperspaces}

\hspace{15pt} First we consider products of Bourbaki-complete and cofinally Bourbaki-complete uniform spaces. 




\begin{theorem} \label{product} Any nonempty product of uniform spaces is Bourbaki-complete if and only if each factor is Bourbaki-complete.

\begin{proof} $\Rightarrow )$ Suppose  $(\prod_{i\in I} X_i,\prod _{i\in I} \mu _i)$ is Bourbaki-complete. Since each factor  $(X_i,\mu _i)$ is uniformly homeomorphic to a closed subspace of this product, then it must be Bourbaki-complete, as Bourbaki-completeness is inherited by closed subspaces.

$\Leftarrow)$ Conversely, suppose that $(X_i,\mu _i)$ is Bourbaki-complete for every $i\in I$ and let $\mathcal{F}$ be a Bourbaki-Cauchy filter in the product. Take  $\mathcal{H}$  an ultrafilter containing $\mathcal{F}$. Clearly, $\mathcal{H}$ is also Bourbaki-Cauchy and then its projection into $x_i$ will be a Bourbaki-Cauchy ultrafilter, for every $i\in I$ (by uniform continuity of the projections). Now, from the Bourbaki-completeness of every factor, this projection must converges to a point in $X_i$. Therefore,  $\mathcal{H}$ also converges to a point in the product, and this means, in particular, that the initial filter $\mathcal{F}$  clusters, as we wanted.
\end{proof}
\end{theorem}

\begin{remark}Consider the euclidean real-line $(\Rset,d_u)$ and a uniformly discrete space $(D, \chi)$. Both metric spaces are Bourbaki-complete by uniform local compactness. Then, any closed subspace of $(D^{\alpha}\times \Rset^{\alpha},\pi)$, where $\pi$ denotes the product uniformity and $\alpha$ is any cardinal, is a Bourbaki-complete space. We will see soon how the above space can be considered universal for all the Bourbaki-complete spaces. On the other hand, it is well-known that {\it a space is $\delta$-complete if and only if it is homeomorphic to a closed subspace of $D ^{\alpha}\times \Rset^{\alpha}$} for some discrete space $D$ and some cardinal $\alpha$ (see \cite{garcia.equality}). Thus, by Theorem \ref{delta}, the space $D ^{\alpha}\times \Rset^{\alpha}$ is universal for the spaces which are Bourbaki-complete when they are endowed with the fine uniformity {\tt u}. 

\end{remark}

\smallskip

Products of cofinally Bourbaki-complete uniform spaces do not have such a good behavior as product of Bourbaki-complete uniform spaces. In fact, products of cofinally complete uniform spaces are not good either.

\begin{theorem} {\rm (\cite{hohti-thesis})} Let $(X, d)$ and $(Y,\rho)$ be cofinally complete metric spaces. Then  $(X\times Y, d+\rho)$ is cofinally complete if and only if $X$ and $Y$ are (uniformly) locally compact or at least one of them is compact.
\end{theorem}

\begin{theorem}{\rm (\cite{marconi})} Let $(X, \mu)$ and $(Y,\nu)$ be cofinally complete uniform spaces. Then if both are (uniformly) locally compact or at least one of them is compact then $(X\times Y , \mu\times \nu)$ is cofinally complete. \end{theorem}

\begin{theorem} {\rm (\cite{hohti-thesis})} \label{product.cof.Bourbaki} If $\{(X_i,\mu_i): i\in I\}$ is an infinite family of uniform spaces such that their product $(\prod_{i\in I}X_i,\prod_{i\in I}\mu_i)$ is uniformly paracompact (cofinally complete) then all but at most finitely many factors  are compact.
\end{theorem}

\begin{remark} From Theorem \ref{uniform.strongly.paracompact}  and Theorem \ref{unif.local.compact.th} we get that the above results are also satisfied when we change cofinal completeness by cofinal Bourbaki-completeness. Indeed, uniform local compactness implies cofinally Bourbaki-completeness and cofinal Bourbaki-completeness implies cofinally completeness (see Theorem \ref{unif.local.compact.th})
\end{remark}



Next, we give an example of the bad behavior of products of cofinally Bourbaki-complete uniform spaces. In it, and from now on, we endow  $D^{\omega _0}$, the product of countably many copies of a discrete space $D$, with the ``{\it first difference metric}" $\rho$ defined by
$$\rho(\langle d_n\rangle , \langle e_n\rangle)=\begin{cases}  \,\,  0 &     \text {if $d_n= e_n$ for every $n\in \Nset$;}  \\ \,\,\, 1/n  & \text{ if $d_j = e_j$ for every $j=1,...,n-1$ and $e_n\neq d_n$}\end{cases}.$$ This metric is compatible with the product uniformity on it.

\begin{example}\label{nagata}{\it There is a countable product of cofinally Bourbaki-complete metric spaces which is not strongly paracompact, nor cofinally Bourbaki-com- plete.}

\begin{proof} Let $D$ be a discrete space of cardinality $\kappa \geq \omega _1$. Then the countable product $( D^{\omega _0}\times \Rset, \rho+d_u)$, is Bourbaki-complete by Theorem \ref{product} as each factor is uniformly locally compact. However, it is not strongly paracompact. Indeed, Nagata proved in \cite[Remark p. 169]{nagata} (see also \cite{pears}) that the space $ D^{\omega _0}\times (0,1)$ is not strongly paracompact whenever $D$ is an uncountable discrete space. Since $ D^{\omega _0}\times \Rset$ and $D^{\omega_0}\times (0,1)$ are homeomorphic, the result follows. In particular $( D^{\omega _0}\times \Rset, \rho+d_u)$ cannot be cofinally Bourbaki-complete because every cofinally Bourbaki-complete uniform space is strongly paracompact as ${\tt u} \geq \mu$.
\end{proof}
\end{example}



\medskip

Now, we proceed with the hyperspaces. The results obtained here are, however, partial.

Let $\mathcal{H}(X)$ and $\mathcal{K}(X)$ the set of all the non-empty closed, respectively compact, subsets of  a topological space $X$. For subsets $A_1,...,A_n$ of $X$ we denote by $\langle A_1,....,A_n\rangle$ the family of elements $B$ in $\mathcal{H}(X)$, respectively $\mathcal{K}(X)$, such that  $B\subset A_1 \cup... \cup A_n$ and $B\cap A_i \neq \emptyset$ for every $i\in \{1,...,n\}$. Whenever we deal with uniform spaces $(X,\mu)$ we can endowed the hyperspaces   $\mathcal{H}(X)$ and $\mathcal{K}(X)$ with the respective uniformity $\langle \mu \rangle$ generated by the family of covers $\langle \mathcal{U}\rangle$, $\mathcal{U}\in \mu$, where the elements of  $\langle \mathcal{U}\rangle$ are the sets of the form $\langle U_1,...,U_n\rangle$ for $U_1,...,U_n \in \mathcal{U}$ (see for instance \cite{hohti-hyper} for bibliography and more information). If $(X,\mu)$ is metrizable then the unifomity $\langle \mu \rangle$ is  metrizable by the {\it Hausdorff metric} $d_{H}$ defined as follows $$d_H(A,B)=inf\{\varepsilon >0: A\subset B^{\varepsilon} \text{ and } B\subset A^{\varepsilon}\}.$$  Observe that the metric $d_H$ of $\mathcal{H}(X)$ is finite only when $(X,d)$ is a bounded metric space.

Next, recall the following result.

\begin{theorem}\label{star.hyper}{\rm (\cite[Proposition 3.1]{hohti-hyper})} If $(X,\mu)$ is a uniform space such that $\mu =s_f\mu$ then the uniformity $\langle \mu \rangle$ of $\mathcal{K}(X)$ has also a star-finite base.  Moreover, for any uniform space $(X,\mu)$ the uniformities  $\langle s_f\mu \rangle$ and $ s_f  \langle \mu \rangle$ are equivalent. 
\begin{proof} We just prove the second statement. That $\langle s_f\mu \rangle \leq  s_f  \langle \mu \rangle$ is satisfied follows from Hohti's result. Conversely, let $\langle \mathcal{U} \rangle \in s_f \langle \mu \rangle $ being star-finite and suppose by contradiction that for some $U\in \mathcal{U}$, $U$ meets infinitely many distinct elements $V_1,..,V_n,...$ in $\mathcal{U}$. But then $\langle U\rangle \cap \langle V_n \rangle \neq \emptyset$ for every $n\in \Nset$ which is a contradiction since $\langle \mathcal{U} \rangle$ is star-finite.
\end{proof}
\end{theorem}

\noindent However, the above statement is not in general true for $\mathcal{H}(X)$ as it was shown by Hohti in \cite{hohti-hyper}.

In \cite[Theorem 1.5]{morita.hyper}, Morita proved that {\it  $(X,\mu)$ is complete if and only if $(\mathcal{K}(X),\langle \mu \rangle)$ is also complete} (see also \cite[Theorem 3.5]{cao}). We prove  now the same for Bourbaki-completeness. Observe that this fact was partially proved in \cite[Proposition 3.3]{merono.hh} for metric spaces.

\begin{theorem} Let $(X,\mu)$ be a uniform space. Then $(\mathcal{K}(X),\langle \mu \rangle)$ is Bourba-ki-complete if and only if  $(X,\mu)$ is Bourbaki-complete.
\begin{proof} 
 The results follows at once by Theorem \ref{star.hyper}, Theorem  \ref{star-finite2} and the above result of Morita.

\end{proof}
\end{theorem}

Recall that {\it for a metric space $(X,d)$,  $(\mathcal{H}(X),d_H)$ is complete if and only if $(X,d)$ is complete} (\cite[Theorem 3.2.4]{beer-book}). However, this is not always true in the frame of uniform spaces (\cite[Chapter II, 46 and Theorem 48]{isbellbook}. On the other hand, for a Bourbaki-complete metric space,  $(\mathcal{H}(X),d_H)$ is not necessarily Bourbaki-complete, as it is shown in the next example.

\begin{example} \label{hyper-example}\cite[Example 3.4]{merono.hh} {\it There exists a Bourbaki-complete metric space $(X,d)$ such that the hyperspace $(\mathcal{H}(X), d_{H})$ is not Bourbaki- complete.} 
\begin{proof}
We are going to prove that the hyperspace $(\mathcal{H}(\Rset), d_{H})$ of the real-line $(\Rset, d^{*})$, where $d^{*}(x,y)=min \{1,d_u(x,y)\}$, is not Bourbaki-complete since it contains a closed isometric copy of the metric hedgehog  $H(\mathfrak{c})$ having continuum-many spines. Thus, $(\mathcal{H}(\Rset), d_{\mathcal{H}})$ cannot be Bourbaki-complete because $(H(\mathfrak{c}),\rho)$ is a complete metric space, in particular, a closed subspace of $(\mathcal{H}(\Rset), d_{\mathcal{H}})$, but it is not Bourbaki-complete.   

Denote by $\mathcal{F}$ the set of all functions $\Nset \rightarrow \{-1,1\}$. For each $f\in \mathcal{F}$, let $E_f=\{4n+f(n):n\in \Nset\}$. Let $g\in \mathcal{F}$. For every $t\in [0,1]$, we set $E_{g}(t)=\{4n+t\cdot g(n):n\in \Nset\}$. The set $\mathcal{L}_g$ is a ``line segment'' in $\mathcal{H}(\Rset)$  joining the ``endpoint'' $E_{g}(1)=E_g$ to the ``origin'' $E_{g}(0)=C$ where $C=\{4n:n\in \Nset\}$. Note the that for all $f,g\in \mathcal{F}$ and $s,t\in [0,1]$, we have that $d_H (E_f(s), E_g(t))=s+t$ if $f\neq g$ and $d_H (E_f(s), E_g(t))=|s-t|$ if $f=g$. Consider the subset $H=\bigcup _{f\in \mathcal{F}} \mathcal{L}_f$ of $\mathcal{H}(\Rset)$ endowed with metric $d_H$ restricted on it. By all the foregoing $(H, d_{H})$ is the isometric copy of $(H(\mathfrak{c}),\rho)$.
\end{proof}
\end{example}

Finally, we study cofinal Bourbaki-completeness of $\mathcal{K}(X)$ and $\mathcal{H}(X)$. As happens to products, the work is already done for cofinal completeness, at least for metric spaces.

\begin{theorem} Let $(X,d)$ be a metric space. The following are equivalent:
\begin{enumerate}

\item $(\mathcal{K}(X), d_{H})$ is uniformly locally compact;

\item $(\mathcal{K}(X), d_{H})$ is cofinally Bourbaki-complete;

\item $(\mathcal{K}(X), d_{H})$ is cofinally complete;

\item $(X,d)$ is uniformly locally compact.

\end{enumerate}
\begin{proof} The proof follows at once applying  \cite[Corollary 4.2]{hohti-hyper}, Theorem \ref{unif.local.compact.th} and Theorem \ref{uniform.strongly.paracompact} .
\end{proof}
\end{theorem}

\begin{theorem}  Let $(X,d)$ be a metric space. The following are equivalent:
\begin{enumerate}

\item $(\mathcal{H}(X), d_H)$ is uniformly locally compact;

\item $(\mathcal{H}(X), d_H)$ is cofinally Bourbaki-complete;

\item $(\mathcal{H}(X), d_H)$ is cofinally complete;

\item $X$ is a point of local compactness of $(\mathcal{H}(X), d_H)$.
\end{enumerate}
\begin{proof} Again, the proof follows at once applying  \cite{burdick},  \cite[Theorem 3.9]{beer-dimaio}, Theorem \ref{unif.local.compact.th} and Theorem \ref{uniform.strongly.paracompact}.
\end{proof}
\end{theorem}

We don't know whether or not the above result are true in the frame of uniform spaces. On the other hand, Burdick proved in \cite{burdick} that the hyperspace $(\mathcal{H}(\Rset), d_{\mathcal{H}})$ (Example \ref{hyper-example}) is not cofinally complete in spite of the fact that $(\Rset, d^{*})$ is uniformly locally compact.

\bigskip

\bigskip

\begin{center} \Denarius \Denarius  \Denarius \end{center}
\bigskip

\bigskip

\section{The particular case of metric spaces}

\subsection{Bourbaki-complete metric spaces}

{\hspace{15pt}} Before studying Bourbaki-completeness of metric spaces, we first show that Bourbaki-Cauchy sequences and cofinally Bourbaki-Cauchy sequences determine Bourbaki-boundedness in the frame of metric spaces.

\begin{definition} A sequence $(x_n)_n$ of a uniform space $(X,\mu)$ is {\it Bourbaki-Cauchy} in $X$ if for every $\mathcal{U}\in \mu$ there exist $m,n_0\in \Nset$ such that for some $U\in \mathcal{U}$, $x_n \in St^{m}(U,\mathcal{U})$ for every $n\geq n_0$.
\end{definition}

\noindent If $(x_n)$ is a Bourbaki-Cauchy sequence of a metric space $(X,d)$, then for every $\varepsilon >0$ there exist, $m, n_0\in \Nset$ such that for every $n\geq n_0$ we can join the points $x_n$ and $x_{n_0}$ by a chain of points $u_0,....,u_m$ in $X$, where $u_0=x_{n_0}$, $u_m=x_n$ and $d(u_{i-1}, u_i)<\varepsilon$ for every $i=1,...,m$.

\begin{definition} A sequence $(x_n)_n$ of a uniform space $(X,\mu)$ is {\it cofinally Bourbaki-Cauchy} in $X$ if for every $\mathcal{U}\in \mu$ there exists $m\in \Nset$ such that for some $U\in \mathcal{U}$ and for every $n\in \Nset $ there exists $k\geq n$ such that $x_k \in St^{m}(U,\mathcal{U})$.
\end{definition}

\noindent Similarly, we have that, if $(x_n)$ is a cofinal Bourbaki-Cauchy sequence of a metric space $(X,d)$, then for every $\varepsilon >0$ there exist, $m,\in \Nset$ and an infinite set $\Nset ^{*}\subset \Nset$ such that for every $n\in \Nset ^{*}$ we can join the points $x_n$ and with some fixed $x_{n_0}$, $n_0\in \Nset ^{*}$ by a chain of points in $X$, satisfying the same properties than above.

\begin{theorem} \label{metric.BB} {\sc (\cite{merono.completeness})}  Let $(X,d)$ be a metric space and $B$ a subset of $X$. The following statements are equivalent:
\begin{enumerate}
\item $B$ is a Bourbaki-bounded subset of $X$;

\item every countable subset of $B$ is a Bourbaki-bounded subset of $X$;

\item every sequence in $B$ has a Bourbaki-Cauchy subsequence in $X$;

\item every sequence in $B$ is cofinally Bourbaki-Cauchy in $X$.

\end{enumerate}
\begin{proof} $(1) \Rightarrow (2)$ This implication is clear since every subset of a Bourbaki-bounded subset in $X$ is a Bourbaki-bounded subset in $X$.
\smallskip

$(2)\Rightarrow (3)$ Let $(x_n)_n$ be a sequence in $B$.  By hypothesis, the set $\{x_n:n\in \Nset\}$ is Bourbaki-bounded in $X$, and then for $\varepsilon =1$ there exist $m_1 \in \Nset$ and finitely many points $p_1 ^1, , ... , p_{j_1}^1 \in X$ such that, $$\{x_n :n \in \Nset\}\subset \bigcup \{B_{d}^{m_1}(p_i ^1, 1): i=1,...,j_1\}.$$ Since the family $\{B_d ^{m_1}(p_i ^1, 1): i=1,...,j_1\}$ is finite, some $B_d ^{m_1}(p_i ^1,1)$ contains infinite terms of the sequence. Therefore, there exists a subsequence $(x_n ^1)_n$ of $(x_n)_n$ inside $B_d ^{m_1}(p_i ^1,1)$.

By induction, we have that, for every $k\geq 2$ and $\varepsilon = 1/k$, there exist some $m_k\in \Nset$ and finitely many points $p_1 ^k, , ... , p_{j_k}^k \in X$ such that, $\{B_d ^{m_k}(p_i ^k, 1/k): i=1,...,j_k\}$  is a finite cover of the set $\{x_n ^{k-1}:n\in \Nset\}$. Then there exists some $B_d ^{m_k}(p_i ^k, 1/k)$  containing some subsequence $(x_n ^k)_n$ of $(x_{n}^{k-1})_n$.

Finally, choosing the standard diagonal subsequence $(x_n^n)_n$ we can easily check that it is the required Bourbaki-Cauchy subsequence of $(x_n)_n$.
\smallskip

$(3)\Rightarrow (4)$ It is clear.

$(4)\Rightarrow (1)$ Assume that $B$ is not a Bourbaki-bounded subset of $X$. Then there exists some $\varepsilon_0 >0$ such that, for every $m\in \Nset$, the family $\{B_{d}^m(x, \varepsilon _0): x\in X\}$ does not contain any finite subcover of $B$. Fix $x_0 \in X$ and for every $m\in \Nset$ choose $x_m\in B$ such that $x_m\notin B^{m}_{d}(x_i,\varepsilon _0)$, for every $i=0,...,m-1$. Then, the sequence $(x_m)_m$ constructed in this way is not cofinally Bourbaki-Cauchy in $x$. Otherwise, for this $\varepsilon _0$ there must exists $m_0\in \Nset$ and an infinite subset $N_0 \subset \Nset$ such that for some $p_0\in X$ we have that $x_n \in B^{m}_{d}(p_0, \varepsilon _0)$ for every $n\in N _0$. Then taking $n_0\in N_0$, we have that there are infinitely many terms of the sequence $(x_m)_m$ in $B_{d}^{2\cdot m_0}(x_{n_0},\varepsilon _0)$, which is a contradiction.
\end{proof}
\end{theorem}

Observe that, in the proof of the above result, we have not use the following result that characterizes total-boundedness by sequences, as we have use Theorem \ref{totally.bounded.filters} to prove Theorem \ref{BB.filters}.

\begin{theorem} \label{TB.sequences}{\rm (\cite{beer-between1})} Let $(X,d)$ be a metric space and $B$ a subset of $X$. The following statements are equivalent:
\begin{enumerate}
\item $B$ is a totally-bounded subset;

\item every sequence of $B$ has a Cauchy subsequence in $X$;

\item every sequence of $B$ is cofinally Cauchy in $X$.

\end{enumerate}
\end{theorem}

\noindent In fact, for a metric space $(X,d)$ the star-finite modification $s_{f}\mu _d$ is not necessarily metrizable and the  Cauchy sequences of $(X,s_f\mu _d)$ are not enough to characterize Bourbaki-boundedness as we will see in Example \ref{BB.non.seq.characterized}.

Next, recall  the following facts.

\begin{definition} A uniform space $(X,\mu)$ is {\it sequentially complete} if every Cauchy sequence clusters (equivalently, converges). 
\end{definition}

\begin{theorem} [\cite{willard}] A metric space $(X,d)$ is complete if and only if it is sequentially complete.
\end{theorem}

\begin{definition} A uniform space $(X,\mu)$ is {\it sequentially Bourbaki-complete} if every Bourbaki-Cauchy sequence clusters. 
\end{definition}

In the following result we check that, for a metric space, Bourbaki-comple- teness and sequentially Bourbaki-completeness coincide, as it happens to usual completeness. Moreover, we also give a characterization by means of bornologies similarly to Theorem \ref{complete.metric}. Nevertheless, we will see that these results are not evident.

\begin{theorem} \label{metric.BC} {\sc (\cite{merono.completeness})} Let $(X,d)$ be a metric space. The following statements are equivalent:
\begin{enumerate}
\item $(X,d)$ is Bourbaki-complete;

\item $(X,d)$ is sequentially Bourbaki-complete;

\item the closure of every $B\in {\bf BB}_{d}(X)$ is compact, that is, ${\bf BB}_{d}(X)={\bf CB}(X);$

\item $(X,d)$ is complete and ${\bf BB}_{d}(X)={\bf TB}_{d}(X).$
\end{enumerate}
\begin{proof} $(1)\Rightarrow (2)$ This implication is easy. In fact, the family of sets $\{\{x_n: n\geq k\}:k\in \Nset\}$ is a filter base which is Bourbaki-Cauchy in $X$ whenever $(x_n)_n$ is Bourbaki-Cauchy in $X$. 
\smallskip

$(2)\Rightarrow (1)$ Let $(X,d)$ be a metric space being sequentially Bourbaki-complete and let $\mathcal{F}$ be a Bourbaki-Cauchy filter in $X$. Suppose by contradiction that $\mathcal{F}$ does not cluster, then for every $x\in X$ there is an open neighborhood $V^{x}$ of $x$ such that $F\cap V^{x}=\emptyset $ for some $F\in \mathcal{F}$. By paracompactness of $X$ the open cover $\mathcal{V}=\{V^{x}:x\in X\}$ has a locally finite open refinement $\mathcal{A}$, which is, in particular, point-finite. We are going to construct inductively a Bourbaki-Cauchy sequence $(x_n)_n$  which does not have any cluster point.
 
First, for every $n\in \Nset$, we can fix $y_n\in X$, and $m_n\in \Nset$ such that $B_{d}^{m_n}(y_n, 1/2^n) \in \mathcal{F}$ because $\mathcal{F}$ is Bourbaki-Cauchy in $X$. Next, fix $x_0\in X$ an arbitrary point. As $\mathcal{A}$ is point-finite, there are only finitely many $A_1,...,A_{k_1} \in \mathcal{A}$ such that $ x_0\in A_i$, $i=1,...,k_1$. Let $\mathcal{A}_0$ denote the finite subfamily of $\mathcal{A}$ consisting of all $A_i$, $i = 1,...k_1$. For each $A_i$, there is some $F_i ^1 \in \mathcal{F}$ such that $F_i ^1 \cap A_i =\emptyset$, so $\bigcap_{1}^{k_1} F_i ^1 \subset X\backslash\bigcup \mathcal{A}_0$ . Then we can take some $$x_1\in \bigcap_{1}^{k_1} F_i ^1 \cap B_{d}^{m_1}(y_1, 1/2^1) \subset (X\backslash\bigcup \mathcal{A}_0 )\cap B_{d}^{m_1}(y_1, 1/2^1). $$ Now, take the finite subfamily $\mathcal{A}_1 = \mathcal{A}_0 \bigcup \{A\in \mathcal{A}: x_{1}\in \mathcal{A}\}$. Similarly to the previous step, there are finitely many $F_i ^2\in \mathcal{F}$, $i=1,...,k_2$ such that  we can take some $x_2\in X$ satisfying that $$ x_2\in \bigcap_{1}^{k_2} F_i ^2 \cap B_{d}^{m_2}(y_2, 1/2^2) \cap B_{d}^{m_1}(y_1, 1/2^1) $$ $$\subset (X\backslash\bigcup \mathcal{A}_1) \cap B_{d}^{m_2}(y_2, 1/2^2)\cap B_{d}^{m_1}(y_1, 1/2^1). $$

In general, for every $n \geq 3$, there  are finitely many $F_i ^n\in \mathcal{F}$, $i=1,...,k_n$ such that  we can take some $x_n\in X$ satisfying that $$ x_n\in \bigcap_{1}^{k_n} F_i ^n \cap \bigcap_{j=1}^{n} B_{d}^{m_j}(y_j, 1/2^j)  \subset (X\backslash\bigcup \mathcal{A}_{n-1}) \cap \bigcap_{j=1}^{n} B_{d}^{m_j}(y_j, 1/2^j), $$ where
$\mathcal{A}_{n-1} = \mathcal{A}_{n-2}\bigcup \{A \in \mathcal{A}: x_{n-1}\in A\}$ is a finite subfamily of $\mathcal{A}$.  Note that the sequence $(x_n)_n$ obtained in this way is Bourbaki–Cauchy. However, $(x_n)_n$ does not have a cluster point. Indeed, for any $y\in X$, there is some $A\in \mathcal{A}$ such that $y\in A$. Since $A$ is open then it is a neighborhood of $y$ and by the construction of the sequence, if $x_k \in A$ implies $x_n\in X\backslash A$, for all $n>k$.
\smallskip

$(2)\Rightarrow (3)$ Let $B\in {\bf BB}_{d}(X)$ and $(x_n)_n$ a sequence in its closure ${\rm cl}_{X} B$. As ${\rm cl}_{X} B \in {\bf BB}_{d}(X)$, then, by Theorem \ref{metric.BB} $(x_n)_n$ has a Bourbaki-Cauchy subsequence which  clusters in ${\rm cl}_{X} B$ by hypothesis. Thus, $(x_n)_n$ clusters too. 

\smallskip $(3)\Rightarrow (2)$ It is not difficult see that every Bourbaki-Cauchy sequence is a Bourbaki-bounded subset of $(X,d)$. Hence, the result follows.
\smallskip

$(3)\Leftrightarrow (4)$ Finally, the equivalence of $(3)$ and $(4)$ is immediate since $${\bf CB}_{d}(X)\subset {\bf TB}_{d}(X)\subset \bf {BB}_{d}(X)$$ and by Theorem \ref{complete.metric} a metric space is complete if and only if ${\bf CB}_{d}(X)= {\bf TB}_{d}(X)$.

\end{proof}
\end{theorem}

\begin{corollary} \label{banach.BC} A Banach space $(E, ||\cdot ||)$ is Bourbaki-complete if and only if it is finite dimensional.
\begin{proof} The proof follows at once from Theorem \ref{metric.BC}, Theorem \ref{boundedness.banach} and the well-known fact that a Banach space is finite-dimensional if and only if it satisfies the {\it Heine-Borel property}, that is, ${\bf B}_{d_{||\cdot ||}}(X)={\bf CB} (X)$.

\end{proof}
\end{corollary}

Now, observe that for a uniform space $(X,\mu)$, the following implications which are a uniform extension of the partial implications of Theorem \ref{metric.BC}, are true.

\begin{center}

$(X,\mu)$ {\it is Bourbaki-complete}

$\Downarrow$

$(X,\mu)$ {\it is complete and} ${\bf TB}_{\mu}(X)={\bf BB}_{\mu}(X)$

$\Downarrow$ 

${\bf CB}(X)={\bf BB}_{\mu}(X)$
 
$\Downarrow$
 
$(X,\mu)$ {\it is sequentially Bourbaki-complete}

\end{center}
 
\noindent The first and the second implications follows are clear from the definitions and the third one from the fact that every Bourbaki-Cauchy sequence is a Bourbaki-bounded subset. However, the reverse implications are not true as it is expected. Indeed, next we give the counterexamples to these.

\begin{example}{\it There exists a  complete uniform space $(X,\mu)$ such that ${\bf BB}_{\mu}(X)={\bf CB}(X)$ but which is not Bourbaki-complete.}
\begin{proof} To give such an example we are going to use Ulam-measurable cardinals. Let us consider the metric hedgehog of $\kappa$ spines where $\kappa$ is an Ulam-measurable cardinal (see Example \ref{hedgehog.measurable}), endowed with the fine uniformity ${\tt u}$. Recall that $(H(\kappa), s_f {\tt u})$ is not complete. Then, by Theorem  \ref{delta}, $(H(\kappa), {\tt u})$ is not Bourbaki-complete. However, $(H(\kappa), {\tt u})$ is complete because, as a metric space, $(H(\kappa), \rho)$ is complete. Moreover, $(H(\kappa), {\tt u})$ satisfies that ${\bf BB}_{{\tt u}}(H(\kappa))={\bf CB}(H(\kappa))$. Indeed, recall that the family of all the real-valued continuous functions on a fine space is exactly $C(X)$ (see Theorem \ref{map.fine}). Therefore, by Theorem \ref{hejcman}, $B\in {\bf BB}_{{\tt u}}(H(\kappa))$ if and only if $B$ is relatively pseudocompact in $H(\kappa)$, where by a {\it relatively pseudocompact subset} of a space $X$ we mean a subset $B$ that satisfies that $f(B)$ is a bounded subset of $\Rset$ for every $f\in C(X)$. Hence, by normality, $cl_{H(\kappa)} B$ is a pseudocompact subspace of $H(\kappa)$, and by metrizabitity, it is, in particular, a compact subspace. Thus ${\bf BB}_{{\tt u}}(H(\kappa))={\bf CB}(H(\kappa))$.
\end{proof}
\end{example}

\begin{example} \label{BB=CB.notcomplete} {\it There exists a uniform space $(X,\mu)$ such that ${\bf BB}_{\mu}(X)={\bf CB}_{\mu}(X)$ but which is not Bourbaki-complete, nor complete}.
\begin{proof} In this example we are going to use again realcompactness. Let $(D,\chi)$ be a uniformly discrete metric of Ulam-measurable cardinal. Then $D$ is not realcompact, that is, the uniform space $(D,e{\tt u})$ is not complete (see Example \ref{discrete.measurable}). In particular, as $e{\tt u}=e\mu _\chi$, because $(D,\chi)$ is a fine space, then $(D,e\mu_\chi)$ is not complete either. However every closed and Bourbaki-bounded subset in $(D, e\mu_\chi)$ is compact. Indeed, the two uniformities $e\mu _\chi$ and $\mu _\chi$ share Bourbaki-bounded subsets as $wU_\chi(D)\leq e\mu _\chi \leq \mu _\chi$ (see Theorem \ref{hejcman}). Therefore, since every Bourbaki-bounded subset of $(D,\chi)$ is a finite set, then every closed and Bourbaki-bounded subset of $(D, e\mu_\chi)$ is compact.
\end{proof}
\end{example}

\begin{example}  {\it There exists a uniform space $(X,\mu)$ which is sequentially Bourbaki-complete, Bourbaki-bounded and not compact. In particular, $${\bf BB}_{\mu}(X)\supsetneq{\bf CB}_{\mu}(X).$$}

\begin{proof}[Construction.] Consider the metric hedgehog (Example \ref{hedgehog}) of uncountable many spines $H(\omega _1)$. Then, $(H(\omega _1), \rho)$ is Bourbaki-bounded and not Bourbaki-complete. In particular, by Theorem \ref{BB.sf}, $s_{f}\mu _\rho=f\mu _{\rho}$ and hence, the completion of $(H(\omega _1), s_f \rho)$ is exactly the Samuel compactification $s_{\rho} H(\omega _1)$. Next, let $Y$ be the subspace of $s_\rho H(\omega _1)$ given by all the cluster points of the Bourbaki-Cauchy sequences of $(H(\omega _1), \rho)$. We are going to check that the space $Y$, endowed with the uniformity  inherited  from $s_{\rho}H(\omega _1)$, is sequentially Bourbaki-complete. 

Take $(\xi _n)_n$ a Bourbaki-Cauchy sequence in $Y$ and $\xi$ cluster point of it in $ s_{\mu _\rho}H(\omega _1)$. We have to prove that $\xi \in Y$. For every $n\in \Nset$, let $(x_k ^n)_k$ be a Bourbaki-Cauchy sequence of $(X,d)$ such that $\xi _n$ is a cluster point of it. Let $f :\Nset \times \Nset\rightarrow  \Nset$ any bijection and define the sequence $(y_{f(n,k)})_{n,k}$, where $y_{f(n,k)}:=x_k ^n$ for every $n,k \in \Nset$. Then, $(y_{f(n,k)})_{n,k}$ is a Bourbaki-Cauchy sequence of $(H(\omega _1), \rho)$ since for every $\varepsilon >0$, $$y_{f(n,k)}\in B_{\rho}^{m}(0, \varepsilon)\text{ for every }n,k\in \Nset$$ where $m$ is any natural number bigger than $\frac{1}{\varepsilon}$

Moreover, $\xi$ is a cluster point of it. Indeed, since $\xi$ is a cluster point of $(\xi _n)_n$ then for every neighborhood $V^\xi$ of $\xi$ in $s_{\mu _\rho} H(\omega _1)$, there exists an infinite set $N\subset \Nset$ such that $\xi _n \in V^\xi$ for every $n\in N$. In addition, for every $\xi_n$ lying in $V^\xi$ there exists an open neighborhood $V^{\xi _{n}}$ of $\xi _n$ in $s_{\mu _\rho} H(\omega _1)$ such that $\xi _n \in V^{\xi _n}\subset V^\xi$. Thus, fixed $n\in N$, since $\xi _n$ is a cluster point of $(x_{k}^n)_k$, there exists some infinite set $M_n\subset \Nset$ such that $$x_{k}^n \in V^{\xi _n} \text{ for every }k\in M_n.$$ Therefore, $$y_{f(n,k)}\in V^\xi \text{ for every }(n,k)\in \bigcup_{n\in N}\{n\}\times M_n .$$ Clearly, $\bigcup_{n\in N}\{n\}\times M_n $ is an infinite subset of $\Nset \times \Nset$, so $f\big(\bigcup_{n\in N}\{n\}\times M_n\big) $ is also an infinite subset of $\Nset$ because $f$ is bijective. Hence $\xi $ is a cluster point of $(y_{f(n,k)})_{n,k}$ and $\xi \in Y$.

Finally observe that  $Y$ is a totally bounded subspace of $s_{\mu _\rho} H(\omega _0)$ and hence it is also Bourbaki-bounded space (in itself). However, it is not compact. Indeed, consider the filter $\mathcal{F}$ of $(H(\omega _1), \rho)$ generated by the family of sets $\{F_{\alpha}: \alpha <\omega _1\}$ where $F_{\alpha}=\bigcup\{[(1, \beta)]: \alpha \leq \beta, \beta <\omega _1\}$. Let $\xi \in s_{\rho}H(\omega _1)$ a cluster point of $\mathcal{F}$. We are going to prove that $\xi \notin Y$ by showing that there is no Bourbaki-Cauchy sequence $(x_n)_n$ of $(H(\omega _1), \rho)$ such that $\xi$ is a cluster point of $(x_n)_n$ .
 
Suppose by contradiction that there exists a Bourbaki-Cauchy sequence $(x_n)_n$ in $(H(\omega _1), \rho)$ such that $(x_n)_n$ clusters to $\xi$. Then, it is clear that $$\xi\in {\rm cl}_{s_\rho H(\omega _1)} \{[(1,\alpha)]:\alpha <\omega _1\} \cap  {\rm cl}_{s_\rho H(\omega _1)} \{x_n:n\in \Nset\}$$ and then $\rho(\{[(1,\alpha)]:\alpha <\omega _1\}, \{x_n:n\in \Nset\})=0$ by Lemma \ref{Samuel4}. This implies that there exists a countable subfamily $\{[(1, \alpha _n)]:n\in \Nset\}$ of $\{[(1,\alpha)]:\alpha <\omega _1\}$ such that $\xi$ is also a cluster point of the sequence $([(1,\alpha_n)])$. But this is a contradiction. Indeed, since the cofinality of $\omega _1$ is exactly $\omega _1$, then there is some $\alpha^* <\omega _1$ such that $\alpha _n <\alpha^*$ for every $n\in \Nset$. Put $F=\{[(1,\beta)]:\alpha^*\leq \beta, \beta<\omega _1 \}$, then $F\in \mathcal{F}$. Since $\rho(\{[(1, \alpha _n)]:n\in \Nset\},F)=2$, by Lemma \ref{Samuel4}, $\xi$ cannot belong to the closure in $s_{\rho }H(\omega _1)$ of both sets, contradicting the fact that $\xi$ is a cluster point of $([(1, \alpha _n)])_n$ and of $\mathcal{F}$.
\end{proof}

\end{example}

In spite of the above results we have the following theorem.

\begin{theorem} \label{completion.BC} For a metric space $(X,d)$ the following statements are equivalent:
\begin{enumerate} 
\item The completion of $(X,d)$ is Bourbaki-complete;

\item  ${\bf TB}_d (X)={\bf BB}_d(X)$;

\item every Bourbaki-Cauchy sequence of $(X,d)$ has a Cauchy subsequence.
\end{enumerate}
\begin{proof} $1)\Rightarrow 2)$ Let $(X,d)$ be a metric space such that its completion $(\widetilde{X}, \widetilde{d})$ is Bourbaki-complete and let $B$ a Bourbaki-bounded subset of $(X,d)$. Then $B$ is also a Bourbaki-bounded subset of $(\widetilde{X}, \widetilde{d})$ since $(X,d)$ is isometrically embedded in its completion. Then $B$ is a totally Bounded subset of $(\widetilde{X}, \widetilde{d})$ because its closure in $\widetilde{X}$ is compact. But $B\subset X$ so it is a totally bounded subsets of $X$.
\smallskip

$2)\Rightarrow 1)$ We are going to prove that every Bourbaki-Cauchy sequence  in the completion $(\widetilde{X},\widetilde{d})$, clusters. 

Take $(y_n)_n$ a Bourbaki-Cauchy sequence in $\widetilde{X}$, then for every $k\in \Nset$ there exists $m_k\in \Nset$ and $n_k\in \Nset$ such that $y_n \in B_{\widetilde{d}}^{m_k}(y_{n_k}, 1/3k)$ for every $n\geq n_k$. Next, for every $k\in \Nset$ let $i_k=sup\{k, n_k\}$. Then, for every $n\leq i_k$, by density of $X$ in $\widetilde{X}$, we can take some $z^{k} _n \in X$, such that $\tilde{d}(z^k _{n}, y_n)<1/3k $. We are going to prove that  the set $B=\{z^{n} _k: k\in \Nset, n\leq i_k\}$ is a Bourbaki-bounded subset in $X$. In order to show that this statement is true we prove first the next claim.
\smallskip

{\bf CLAIM.} Let $(X,d)$ a metric space and let $(\widetilde{X}, \tilde{d})$ denote its completion. Then for every $\varepsilon >0$, $y\in \widetilde{X}$ and $m\in \Nset$, $$ \emptyset \neq B^{m} _{\tilde{d}}(y, \varepsilon /3)\cap X\subset B_{d}(x,\varepsilon), \text{ for some }x_0\in X.$$ 

\begin{proof}[Proof of the claim] By density of $X$ in $\widetilde{X}$, we can take some $z\in B^{m} _{\tilde{d}}(y, \varepsilon /3)\cap X$. Then, we can fix some chain of points $u_0, u_1,...,u_m$ in $\widetilde{X}$ such that $u_0=y$, $u_m =z$ and $\tilde{d}(u_{i-1}, u_i)<1/3\varepsilon$, for every $i=1,...,m$. Again, by density of  $X$ in $\widetilde{X}$, we can take $x_0,x_1,...,x_{m-1} \in X$ such that $\tilde{d}(x_{i-1}, u_{i-1})<1/3\varepsilon$, for every $i=1,...,m$. Put $x_m=z$. In particular, $$d(x_{i-1}, x_i)=\tilde{d}(x_{i-1}, x_i)\leq \tilde{d}(x_{i-1}, u_{i-1})+\tilde{d}(u_{i-1}, u_i)+\tilde{d}(u_{i}, x_i)<1/\varepsilon$$ for every $i=1,...,m$. Therefore, $z\in B^{m}_d(x_0, \varepsilon)$ and the claim follows.

\end{proof}

Next, observe that for every $k\in \Nset$, $$B\subset B^{m_k}_{\tilde{d}}(y_{i_k},1/3k)\cup \bigcup_{n< i_k} B_{\tilde{d}}(y_{n}, 1/3k) \cup \{z^j _{n}:n<i_k, j<k\}\subset \widetilde{X}.$$ Therefore, by the above claim, for every $k\in \Nset$, there exists finitely many points $x_i\in X$, $i=1,...,i_k$ such that

$$B\subset B^{m_k}_{d}(x_{i_k},1/k)\cup \bigcup_{n< i_k} B_{d}(x_{n}, 1/k) \cup \{z^j _{n}:n<i_k, j<k\}\subset X.$$ Since the set $\{z^j _{n}:n<i_k, j<k\}$ is finite, it follows that $B$ is Bourbaki-bounded in $X$. More precisely, by the hypothesis $B$ is a totally bounded subset of $X$.

Finally, notice that, by the choose of the points of $B$, $\{y_n:n\in \Nset\}\subset {\rm cl}_{\widetilde{X}}B$. Moreover, by total boundedness of $B$, the closure $\subset {\rm cl}_{\widetilde{X}}B$ is a totally bounded subset of the completion $(\widetilde{X}, \tilde{d})$. Therefore, by completeness, $\subset {\rm cl}_{\widetilde{X}}B$ is compact and the the sequence $(y_n)_n$ clusters in $\widetilde{X}$, that is, $(\widetilde{X}, \tilde{d})$ is Bourbaki-complete.

\smallskip

$2)\Rightarrow 3)$ Let $(x_n)_n$ be a Bourbaki-Cauchy sequence of $(X,d)$. Then $\{x_n:n\in \Nset\}$ is a Bourbaki-bounded subset of $(X,d)$, and by hypothesis it is a totally bounded subset. Then, by Theorem \ref{TB.sequences} the result follows.

\smallskip

$3)\Rightarrow 2)$ Conversely, let $B$ a Bourbaki-bounded subset in $X$. By Theorem \ref{metric.BB} $B$, every sequence contained in $B$ has a Bourbaki-Cauchy subsequence in $X$. Therefore, by hypothesis, it contains also a Cauchy subsequence. Finally, by Theorem \ref{TB.sequences} the result follows.

\end{proof}
\end{theorem}

\medskip

\subsection{Further considerations}

\hspace{15pt} By Theorem \ref{totally.bounded.filters} and Theorem \ref{BB.sf}, Bourbaki-boundedness in a uniform space can be characterized by Cauchy filters of its star-finite modification. More precisely, a subset $B$ of a uniform space $(X,\mu)$ is Bourbaki-bounded in $X$ if and only if every filter $\mathcal{F}$ in $B$ is contained in some Cauchy filter $\mathcal{F}'$ of $(X, s_f \mu)$ in $B$. Moreover, Theorem \ref{star-finite2}  states that Bourbaki-completeness of a uniform space is equivalent to completeness of its star-finite modification.

However,  we next show  that, for a metric space $(X,d)$, Cauchy sequences of $(X, s_f\mu _d)$ are not strong enough to characterize Bourbaki-completeness, nor Bourbaki-boundedness of $(X,d)$ (see Example \ref{BB.non.seq.characterized}). Hence, we show in this way that the robustness of Bourbaki-Cauchy sequences is necessary.

\begin{lemma} {\rm(Efremovi\v{c}'s Lemma \cite{proximity})}\label{efremovic1} Let $\mathcal{V}\in \mu$ and $(x_n)_n$ and $(y_n)_n$ two sequences of a uniform space $(X,\mu)$  satisfying that for every $n\in \Nset$ there are $V(x_n),V(y_n)\in \mathcal{V}$ such that $x_n\in V(x_n)$, $y_n\in V(y_n)$ and $V(x_n)\cap V(y_n)=\emptyset$. Then, there exist subsequences $(x_{n_k})_k$ and $(y_{n_k})_k$ of $(x_n)_n$ and $(y_n)_n$, respectively, and some $\mathcal{U}\in \mu$, $\mathcal{U}^{*}<\mathcal{V}$, such that for every $k\in \Nset$ there are $U(x_{n_k}),U(y_{n_k})\in \mathcal{U}$ satisfying that $x_{n_k}\in U(x_{n_k})$, $y_{n_k}\in U(y_{n_k})$ and $U(x_{n_k})\cap U(y_{n_j})=\emptyset$  for every $k,j \in \Nset$.
\end{lemma}

\begin{theorem} \label{seq.complete}Let $(X,\mu)$ be a uniform space and $\nu$ a  uniformity on $X$ such that $f\mu \leq \nu \leq \mu$. If $(X,\mu)$ is sequentially complete then $(X,\nu)$ is also sequentially complete.

\begin{proof} Let $(x_n)_n$ be a Cauchy sequence in $(X,\nu)$. If  $(x_n)_n$ is Cauchy in $(X,\mu)$ there is nothing to prove. Otherwise, $(x_n)_n$ does not converge in $X$. More precisely, it does not cluster either, because it is a Cauchy sequence.

Therefore the sequence $(x_n)_n$ has a subsequence  $(y_n)_n$ satisfying that for some $\mathcal{U}\in \mu$, for every $n\in \Nset$ there is some $U(y_n)\in \mathcal{U}$ such that $y_n\in U(y_n)$ and $U(y_n)\cap U(y_k)=\emptyset$ for every $n,k\in \Nset$, $n\neq k$. Indeed, since $(x_n)_n$ is not Cauchy in $(X,\mu)$ there exists some $\mathcal{V}\in \mu$ such that for every $n\in \Nset$ we can take some $k_n>n$ for which $V(x_n)\cap V(x_{k_n})=\emptyset$, where $V(x_n), V(x_{k_n})\in \mathcal{V}$. Applying Lemma \ref{efremovic1} to the sequences $(x_n)_n$ and $(x_{k_n})_n$ we obtain the subsequence $(y_n)_n$. By all the foregoing, the subspace $Y=\{y_n:n\in \Nset\}$ of $(X,\mu)$ is uniformly homeomorphic to the uniformly discrete space $(\Nset,\chi)$. Then, $s_{\mu |_Y} Y$ is homeomorphic to $\beta \Nset$. In fact, since $(\Nset, \chi)$ is uniformly discrete then it is a fine space and $f\mu_{\chi} =f{\tt u}$, that is, $s_{\chi} \Nset$ is homeomorphic to $\beta \Nset$.

Now, as $(x_n)_n$ is Cauchy in $(X,\nu)$, then $(y_n)_n$ is also Cauchy in $(X, \nu)$ and in particular, it is Cauchy in $(X, f\mu)$. So $(y_n)_n$ converges to some $\xi$ in the Samuel compactification $s_\mu X$. In particular $\xi \notin X$ and, by Lemma \ref{Samuel3}, $\xi\in s_{\mu |_Y} Y={\rm cl}_{s_{\mu } X}Y$. But this is a contradiction, since $s_{\mu |_Y} Y$ is homeomorphic to $\beta \Nset$ and it is well-known that $\beta \Nset$ does not have non-trivial converging sequences \cite{engelkingbook}.

\end{proof}
\end{theorem}

\begin{corollary} If $(X,d)$ is a complete metric space then $(X,s_f \mu_d)$ is sequentially complete.
\end{corollary}


\begin{example}\label{BB.non.seq.characterized} {\it There exists a non-compact Bourbaki-bounded space $(X,d)$ such that $(X, s_f \mu _d)$ is sequentially complete and not every sequence of $X$ has a Cauchy subsequence in $(X, s_f \mu _d)$.} 

\begin{proof} It is enough to take any complete metric space which is in addition Bourbaki-bounded and not Bourbaki-complete. For instance, take the metric hedgehog $(H(\omega _0),\rho)$ (Example \ref{hedgehog}). Then, by Theorem \ref{seq.complete}, $(H(\omega _0),s_f \mu_\rho)$ is also sequentially complete. Next suppose that every sequence of $(H(\omega _0),\rho)$ has a Cauchy subsequence of $(H(\omega _0),s_f \mu_\rho)$. Then every sequence of $H(\omega _0)$ clusters by sequential completeness of $(H(\omega _0),s_f \mu_\rho)$. But this is not possible since $H(\omega _0)$ is not compact.

\end{proof}
\end{example}

Next, we study a problem that shows the strength of Bourbaki-bounded subsets against totally bounded subsets. Recall that by Theorem \ref{BB.sf}, Theorem \ref{star-finite2} and Theorem \ref{metric.BC}, {\it for a metric space $(X,d)$, the space $(X, s_f \mu_d)$ is complete if and only if ${\bf TB}_{s_f \mu _d}(X)={\bf CB}(X)$}.  Now we ask if the same is possible for the point-finite modification of a metric space, that is, if $(X, p_f \mu _d)$ is complete if and only if ${\bf TB}_{p_f \mu_d}(X)={\bf CB}(X)$. The answer is no as we will see with the help of the next result.

\begin{theorem}\label{tot.bounded.point.finite} {\rm(\cite[Theorem 3.5]{manisha} for metric spaces)} Let $B$ a subset of a uniform space $(X,\mu)$. The following statements are equivalent:
\begin{enumerate}
\item $B$ is totally bounded in $(X,\mu)$;

\item for every point-finite  cover $\mathcal{U}\in \mu$ there exists a finite subfamily $\{U_i:i=1,...,k\}\subset \mathcal{U}$ which covers $B$;

\item for every point-countable  cover $\mathcal{U}\in \mu$ there exists a finite subfamily $\{U_i:i=1,...,k\}\subset \mathcal{U}$ which covers $B$;

\item for every locally finite  cover $\mathcal{U}\in \mu$ there exists a finite subfamily $\{U_i:i=1,...,k\}\subset \mathcal{U}$ which covers $B$;

\item for every locally-countable  cover $\mathcal{U}\in \mu$ there exists a finite subfamily $\{U_i:i=1,...,k\}\subset \mathcal{U}$ which covers $B$;

\item for every star-countable cover $\mathcal{U}\in \mu$ there exists a finite subfamily $\{U_i:i=1,...,k\}\subset \mathcal{U}$ which covers $B$;

\item for every countable  cover $\mathcal{U}\in \mu$ there exists a finite subfamily $\{U_i:i=1,...,k\}\subset \mathcal{U}$ which covers $B$;
\end{enumerate}

\begin{proof} 

That $(1)$ implies all the other statements is trivial. We proof now  $(n)\Rightarrow (1)$, for every $n=2,...,7$, at once.

Assume that $B$ is not totally bounded.  Then it is not difficult to find some open $\mathcal{U}\in \mu$ for which there exist infinitely many $x_i \in B$, $i\in \Nset$, such that  $St^{4}(x_n,\mathcal{U})\cap St^{4}(x_m,\mathcal{U})=\emptyset$ whenever $n\neq m$. Write $A=\{x_n :n\in \Nset\}$. 
Consider the cover $\mathcal{G}=\{St^{3}(x_n,\mathcal{U}), X\backslash A: n\in \Nset\}$. Then $\mathcal{G}$ is uniform as $\mathcal{U}<\mathcal{G}$ and we have to prove that it is  locally finite (in particular locally countable, star-countable, countable, point-finite  and point-countable). Let $x\in X \backslash A$, if $x\in St^{2}(x_n, \mathcal{U})$ for some $n\in \Nset$, then $St(x, \mathcal{U})\subset St^{3}(x_n, \mathcal{U})$. Otherwise $St(x,\mathcal{U})\subset X\backslash A$. Moreover, $St(x,\mathcal{U})$ intersects only finitely many members of $\mathcal{G}$. Thus, $\mathcal{G}$ is a countable and locally finite. But $\mathcal{G}$ has no finite subfamily covering $B$ which contradicts the hypothesis.
\end{proof}

\end{theorem}



Now, recall Pelant's result in \cite{pelant.complete} stating that the point-finite modification of the Banach space $(\ell_{\infty}(\omega_1), ||\cdot ||_{\infty})$ is not complete. Therefore not every complete metric space $(X,d)$ satisfies that $(X, p_f\mu _d)$ is complete. This fact and the above Theorem \ref{tot.bounded.point.finite} answer negatively our question. Indeed, we wanted to know if for a metric space $(X,d)$, the point-finite modification $(X, p_f\mu _d)$ is complete if and only if ${\bf TB}_{p_f \mu_d}(X)={\bf CB}(X)$. However, by Theorem \ref{tot.bounded.point.finite} ${\bf TB}_{p_f \mu_d}(X)={\bf TB}_{d}(X)$, and since there exists complete metric spaces such that its point-finite modification $(X, p_f\mu _d)$ is not complete, it follows that the bornology ${\bf TB}_{p_f \mu_d}(X)$ is too weak to characterize the completeness of  $(X, p_f\mu _d)$ (recall Theorem \ref{complete.metric}).

Finally we give an example of a Bourbaki-complete metric space which does not have a point-finite base. This example is motivated by the fact that every cofinally Bourbaki-complete uniform space has a star-finite base for its uniformity (see Theorem \ref{uniform.strongly.paracompact}). Thus, we show in this way that, differently to cofinally Bourbaki-complete spaces, not every Bourbaki-complete uniform space has a star-finite base.

\begin{example} \label{point-finite} {\it There exists a Bourbaki-complete metric space which does not have a point-finite base for its uniformity. Therefore, the uniformity does not have a star-finite base either.}

For every $n\in \Nset$, let  $X_n=\ell_{\infty}(\omega_1)$ be the set of all bounded real-valued functions over a set of cardinality $\omega_1$, endowed with the metric  $$t_n (x,y)=\begin{cases}  \,\,\,  0  & \text{ if $x=y$}   \\ \,\,\,\frac{1}{2^{n}} +min \{ 1,||x-y||_{\infty}\} &     \text { if $x\neq y$.} \end{cases}$$ Then $(X_n, t_n)$ is Bourbaki-complete and has a point-finite base for its uniformity since it is a uniformly discrete metric space.

Now, let $X=\biguplus _{n\in \Nset} X_n$ the set given by the disjoint union of the above spaces, and let us endowed $X$ with the metric $$t (x,y)=\begin{cases}  \,\,\,  t_n (x,y)  & \text{ if for some $n\in \Nset$ $x,y\in X_n$}   \\ \,\,\,\ 2 &     \text { otherwise.} \end{cases}$$ Then $(X,t)$ is a Bourbaki-complete metric space because it is a disjoint union of uniformly separated Bourbaki-complete metric spaces.  However, $(X,t)$ fails to have a point-finite base for its uniformity as we are going to prove next. 

Let $d_{\infty}(x,y)=||x-y||_{\infty}$ the usual metric on $\ell _{\infty}(\omega _1)$. By Pelant's result \cite{pelant.complete} $(\ell_{\infty}(\omega _1), d_{\infty})$ does not have a point-finite base for its uniformity, so  we can choose  some $N\geq 2$ such that every uniform refinement, for the metric  uniformity induced by $d_{\infty}$, of the cover $\{B_{d_{\infty}}(x, 1/2^N):x\in \ell_{\infty}(\omega_1) \}$ fails to be point-finite. In particular, the same is true for the space $(\ell_{\infty}(\omega_1),\rho)$, where $\rho(x,y)=min\{1,d_{\infty}(x,y)\}$, by uniform equivalence of the metrics $\rho$ and $d_{\infty}$.

Now, observe that for every $n\in \Nset$ and every $x\in \ell_{\infty}(\omega_1)$, $$(\clubsuit) \text{ }B_{t_n}(x, 1/2^{n-1})= B_{\rho}(x, 1/2^{n}).$$ Moreover, for every $n,k\in \Nset$, $n>k$ and every $x\in \ell_{\infty}(\omega_1)$ $$(\heartsuit) \text{ }B_{t_n}(x, 1/2^{k})\subset B_{\rho}(x, 1/2^{k}).$$ Take the uniform cover $\mathcal{B}=\{B_{t}(x,1/2^N):x\in X\}$ of $X$. Then, it is clear that,  $$\mathcal{B}=\biguplus_{n\in \Nset} \{B_{t_n}(x,1/2^N):x\in X_n\}$$  as $N\geq 2$. We are going to prove that every uniform refinement $\mathcal{V}$ of $\mathcal{B}$  fails to be point-finite. Indeed, since $\mathcal{V}$ is uniform we can choose some $m\in \Nset$, $m>N$ such that  the cover $ \{B_{t}(x,1/2^{m}):x\in X\}=\biguplus_{n\in \Nset} \{B_{t_n}(x,1/2^m):x\in X_n\}$ refines $\mathcal{V}$. Thus, we can write also that $\mathcal{V}=\biguplus _{n\in \Nset} \mathcal{V}_n$ where each $\mathcal{V}_n$ is a uniform cover of $(X_n, t_n)$ that refines $\{B_{t_n}(x,1/2^N):x\in X_n\}$. Therefore, whenever $n>N$, $\mathcal{V}_n$ is a uniform cover of $(\ell_{\infty}(\omega_1),\rho)$ by $(\clubsuit)$, and it refines $\{B_{\rho}(x,1/2^N):x\in \ell_{\infty}(\omega _1)\}$ by $(\heartsuit)$. Then $\mathcal{V}_n$ fails to be point-finite for every $n>N$ which meas that $\mathcal{V}$ is not point-finite. Finally, we can conclude that $(X,t)$ does not have a point-finite base.

\qed
\end{example}

\begin{remark}Since every star-finite cover is point-finite, then  for every Bourbaki-complete uniform space $(X,\mu)$, $(X, p_f \mu)$ is complete. On the other hand, by \cite[Corollary 2.4]{pelant-c0}  a metric space has a point-finite base if and only if it can be uniformly embedded into $c_0 (\kappa)$ where  $\kappa$ is the density of $X$. Recall that $c_0 (\kappa)\subset \ell_{\infty}(\kappa)$ denotes the Banach space of function $f:\kappa \rightarrow \Rset$ such that the cardinality of the support $sup(f)$ of $f$  is at most countable and converges to $0$. Then, by Corollary \ref{banach.BC}, $(c_0 (\kappa), ||\cdot ||_{\infty})$, $\kappa \geq \omega_1$, is a example of non-Lindel\"of complete non-Bourbaki-complete space having a point-finite base.
\end{remark}

\medskip

\subsection{Cofinally Bourbaki-complete metric spaces}

\hspace{15pt} In spite of having proved that Cauchy sequences of $(X,s_f \mu _d)$ are not enough to characterize Bourbaki-boundedness of a subset of a metric space $(X,d)$, we have that cofinally Cauchy sequences of $(X,s_f \mu _d)$ are so. The reason is clear. Indeed, by Lemma \ref{ultrafilter} it is easy to see that cofinally Cauchy sequences of $(X, s_f \mu _d)$ and cofinally Bourbaki-Cauchy sequences of $(X,d)$ are the same thing. Thus, applying Theorem \ref{metric.BB} we get next result.

\begin{theorem} Let $(X,d)$ be a metric space. The following statements are equivalent:
\begin{enumerate}
\item $B$ is a Bourbaki-bounded subset in $(X,d)$;

\item every sequence $(x_n)_n$ of $B$ is a cofinal Cauchy sequence of $(X, s_f\mu _d)$.
\end{enumerate}

\end{theorem}

Next, we recall some facts about cofinal completeness of metric spaces. Let $nlc(X)$ be the subset of $X$ of all the points which does not have a compact neighborhood and $nlc(X)^{\varepsilon}=\bigcup \{B_{d}(x,\varepsilon):x\in nlc(X)\}$ for $\varepsilon >0$. Clearly a space is locally compact if and only if $nlc(X)=\emptyset$. In \cite{hohti-thesis} Hohti  studied uniformly paracompactness in the frame of metric spaces and gave the next precise metric characterization by means of the set $nlc(X)$.

\begin{theorem} \label{hohti.metric} {\rm (\cite[Theorem 2.1.1]{hohti-thesis})} A metric space is uniformly paracompact (cofinally complete) if and only if it is either uniformly locally compact, or either $nlc(X)$ is a non-empty compact set such that for every $\varepsilon >0$, the space $X\backslash nlc(X)^{\varepsilon}$ is  a uniformly locally compact.
\end{theorem}

In \cite{beer-between1} Beer studied also cofinal completeness in the frame of metric spaces, precisely, what we call sequentially cofinal completeness.

\begin{definition} A uniform space is {\it sequentially cofinally complete} if every cofinal Cauchy sequence clusters.
\end{definition}

\begin{theorem} {\rm (\cite[Theorem 3.2]{beer-between1})} \label{beer} A metric space $(X,d)$ is cofinally complete if and only if it is sequentially cofinally complete.
\end{theorem}
\bigskip

Now, we prove that sequential cofinal Bourbaki-completeness, defined below, and cofinal Bourbaki-completeness of a metric space are also equivalent. The proof is not trivial and it needs the $(\star)$-property of Theorem \ref{star} together with  the above Beer's result. 

\begin{definition} A uniform space is {\it sequentially cofinally Bourbaki-complete} if every cofinal Bourbaki-Cauchy sequence clusters.
\end{definition}

The following result is a metric generalization of the Efremovich's lemma \ref{efremovic1}  for the infinite countable case. We don't know if it is a known result, but it has been privately communicated by A. Hohti and H. Junnila.

\begin{lemma} \label{efremovich}  Let $\varepsilon >0$ and $\{E_n:n\in \Nset\}$ a countable family of infinite $\varepsilon$-discrete subsets of $(X,d)$. Then there are infinite subsets $E'_n\subset E_n$ such that $\bigcup_{n\in \Nset}E'_n$ is $\varepsilon /2$-discrete.
\begin{proof} We will construct the subsets $E'_n =\{x_{k}^{n}:k\in \Nset\}$ by induction. Choose any $x^{1}_{1}\in E_1$ and write $F_{1}=\{x^{1}_{1}\}$. If $d(x^{1}_1 ,y)<\varepsilon /2$ for every $y\in E_2$, then $d(y ,y' )<\varepsilon $ for any two $y ,y' \in E_2$, which would be impossible. Hence, we can find an $x_{2}^{2}\in E_2$ such that $d(x_{1}^{1}, x_{2}^{2})\geq \varepsilon /2$. In the same way, we can choose a point $x_{1}^{2}\in E_{1}-\{x_{1}^{1},x_{2}^{2}\}$ such that $F_2=\{x_{1}^{1},x_{2}^{2},x_{1}^{2} \}$ is $\varepsilon /2$-discrete.

In general, suppose that we have chosen a finite set $F_n$  of points of the sets $E_1,...,E_n$ satisfying: \begin{enumerate}

\item $F_n$ is $\varepsilon /2$-discrete;

\item $|F_n\cap E_i|=n-i+1$ for all $1\leq i\leq n$.
\end{enumerate}

We construct $F_{n+1}$ as follows. We claim that there is $x_{n+1}^{n+1}\in E_n+1$ such that $F_{n}\cup \{x^{n+1}_{n+1}\}$ is $\varepsilon /2$-discrete. If not, then for all $y\in E_{n+1}$ we have $d(y, F_n)<\varepsilon /2$. As $E_{n+1}- F_n$ is infinite and $F_n$ is finite, there are $y, y' \in E_{n+1}-F_n$, $y \neq y'$ such that $d(y, z)<\varepsilon /2$ and $d(y', z)<\varepsilon /2$ for some $z \in F_n$. Then $d(y, y')<\varepsilon$, which would be a contradiction. Repeating this process $n$ times for $i=n,(n-1),...,1$ and for the finite $\varepsilon /2$-discrete set $F_{n+1}=F_n\cup \{x^{n+1}_{n+1}, x^{n}_{n+1},...,x^{i+1}_{n+1}\}$ together with the infinite set $E_i$, we obtain the desired $\varepsilon /2$-discrete set $F_{n+1}=F_n\cup \{x^{n+1}_{n+1}, x^{n}_{n+1},...,x^{1}_{n+1}\}$.

Since clearly $\bigcup_{n\in \Nset}F_n$ is $\varepsilon /2$-discrete, we finish taking $E'_{n}=\{x_{n}^{n}, x_{n+1}^{n},...\}$.
\end{proof}
\end{lemma}

Next we give the correct proof of a result that can be found in \cite[Theorem~28]{merono.completeness}. This proof needs the previous lemma which wasn't set in  \cite{merono.completeness}.

\begin{theorem} \label{sequentially.cof.Bourbaki} Let $(X,d)$ be a metric space. The following statements are equivalent:
\begin{enumerate}
\item $(X,d)$ is cofinally Bourbaki-complete;

\item $(X,d)$ is sequentially cofinally Bourbaki-complete;

\item $(X,d)$ is (sequentially) cofinally complete and satisfies the $(\star)$-property. 

\item  either $(X,d)$ is uniformly locally compact or $nlc(X)$ is a non-empty compact set such that for every $\varepsilon >0$, $X\backslash nlc(X)^{\varepsilon}$
is a uniformly locally compact space and there exists $\delta >0$  satisfying that for every $x\in nlc(X)$, and every $n\in \Nset$ there exist finitely many $x_1,...,x_k\in X$ satisfying that $B^n _{d}(x,\delta)\subset\bigcup_{i=1}^{k}B_d(x_i,\varepsilon)$.  
\end{enumerate} 
\begin{proof} $(1)\Rightarrow (2)$ This implication is trivial.

$(2)\Rightarrow (3)$ If $X$ is sequentially cofinally Bourbaki-complete then it is sequentially cofinally complete. By Theorem \ref{beer}, $X$ is cofinally complete.

Next, assume that there exists $\varepsilon _0 >0$ satisfying that for every $\delta >0$ there exist some $x_\delta \in X$ and $m_\delta \in \Nset$ such that $B^{m_\delta}_{d}(x_\delta, \delta)$ cannot be covered by finitely many balls $B_{d}(x, \varepsilon _0)$. Take $\delta =\frac{1}{n}$, $n\in \Nset$. Then for every $n\in \Nset$ we can choose an infinite set  $E_{n}\subset B^{m_n}_{d}(x_n, 1/n)$ which is $\varepsilon _0$-discrete.

By Lemma \ref{efremovich}, there exists an infinite set $E'_n\subset E_n$ such that $\bigcup _{n\in \Nset}E'_n$ is $\varepsilon _{0}/2$-discrete. Now, consider a partition of $\Nset$, in a countable family of infinite subsets $\{M _n:n\in \Nset\}$ . Finally, if we enumerate every set $E'_n=\{x^{n}_j:j\in M_n\}$, and we define the sequence $y_j:=x^{n}_j$, if $j\in M_n$, then $(y_j)_j\in \Nset$ is a cofinally Bourbaki-Cauchy sequence which does not cluster. 

$(3)\Rightarrow (1)$ This follows from Theorem \ref{star}  and Theorem \ref{uniform.strongly.paracompact}.

$(3)\Rightarrow (4)$ This implication follows from Theorem \ref{beer} and the fact that $(X,d)$ satisfies the $(\star)$-property.

$(4)\Rightarrow (3)$ If $X$ is uniformly locally compact in particular it is cofinally complete. More precisely, it is cofinally Bourbaki-complete and by Theorem  \ref{star} and Theorem \ref{uniform.strongly.paracompact}, the result follows. Otherwise, $(X,d)$ must be cofinally complete by \cite[Theorem 3.2]{beer-between1} and we just need to prove the $(\star)$-property. So fix $\varepsilon >0$,  let $\delta >0$ from the hypothesis and take $x\in X$. If  $x\in B^{m}(y,\delta)$ for some $m\in \Nset$ and $y\in nlc(X)$ then, by hypothesis, $B^{n}(x,\delta)$ is covered by finitely many balls of the radius $\varepsilon$,  for every $n\in \Nset$. Otherwise, suppose that $x\in X\backslash \{\bigcup B^{\infty} _{d}(y,\delta):y \in nlc(X)\}$. Then $B_{d}^{n}(x,\delta)\subset X\backslash nlc(X)^{\delta}$ for every $n\in \Nset$. Since $X\backslash nlc(X)^{\delta}$ is uniformly locally compact in its relative metric, there exists $\gamma >0$ such that ${\rm cl _{X\backslash nlc(X)^{\delta}}} (B_{d}(z,\gamma)\cap (X\backslash nlc(X)^{\delta}))$ is compact for every $z\in X\backslash nlc(X)^{\delta}$. Let $\alpha <min\{\delta, \gamma\}$. Then using Lemma \ref{unif.loc.compact}, we have that $${\rm cl} _{X\backslash nlc(X)^{\delta}} (B_{d}^n (x,\alpha)\cap (X\backslash nlc(X)^{\delta}))={\rm cl} _{X\backslash nlc(X)^{\delta}} B_{d}^n (x,\alpha)$$ is also compact for every $n\in \Nset$. By compactness, for every $\varepsilon >0$ there exists finitely many balls $B_{d}(z,\varepsilon)$ covering $B_{d}^n(x,\varepsilon)$.

\end{proof}
\end{theorem}

Observe that statement $(4)$ in the previous result give us a good metric characterization of the cofinally Bourbaki-complete metric spaces. In particular, we can see that,

\begin{center} {\it cofinally Bourbaki-complete} $\nRightarrow$ {\it uniformly locally compact}
\end{center}
\noindent as it is shown in the next example.

\begin{example} {\it There exists a cofinally Bourbaki-complete metric space which is not uniformly locally compact.}
\begin{proof} Let $X$ be the following subspace of the metric hedgehog $(H(\omega _0),\rho)$ (Example \ref{hedgehog}): $$X=\{0\}\cup \{[(1/n, \alpha)]: n\in \Nset, \alpha <\omega _0\}.$$  Then, $(X,\rho|_X)$  is cofinally Bourbaki-complete by statement $4$ in Theorem \ref{sequentially.cof.Bourbaki}. However, $(X,\rho|_X)$ is not uniformly locally compact as the point $0$ has no compact neighborhood.
\end{proof}
\end{example}

\begin{corollary} Let $(X,d)$ be a metric space. Then $(X,s_f \mu _d)$ is cofinally complete if and only if $(X,s_f \mu _d)$ is sequentially cofinally complete.
\begin{proof} One implication is clear. So let us suppose that $(X,s_f \mu _d)$ is sequentially cofinally complete. Similarly to Lemma \ref{ultrafilter} (3), one can prove that a sequence $(x_n)_n$ is cofinally Cauchy in $(X,s_f \mu _d)$ if and only if it is cofinally Bourbaki-Cauchy in $(X,d)$. Therefore the result follows from Theorem \ref{uniform.strongly.paracompact} and Theorem \ref{sequentially.cof.Bourbaki}.
\end{proof}
\end{corollary}

Theorem \ref{sequentially.cof.Bourbaki} contrasts with the fact that sequential completeness of $(X,s_f \mu _d)$ is weaker than Bourbaki-completeness of $(X,d)$. The reason is that cofinal completeness of $(X,s_f \mu _d)$ implies that $s_f\mu_d=\mu_d$, that is, the uniformity $s_f\mu_d$ is metrizable.

\begin{example}{\it There is a uniform space which is sequentially cofinally Bourbaki-complete, totally bounded but not complete, nor cofinally complete.} 

\begin{proof}The space $[0,\omega _1)$ of all the countable ordinals is totally bounded and {\it sequentially compact}, that is, every infinite sequence clusters \cite{engelkingbook}.  Therefore it is sequentially cofinally Bourbaki-complete in its unique uniformity. However it is not compact and then not complete, nor Bourbaki-complete. 

\end{proof}
\end{example}

\bigskip

\bigskip

\begin{center} \Denarius \Denarius  \Denarius \end{center}
\bigskip

\bigskip

\newpage

\chapter{Embedding Bourbaki-complete spaces and Bourbaki-completely metrizable spaces} 
\thispagestyle{empty}
\newpage
\section{Embedding's results}

\subsection{``Universal space" for Bourbaki-complete uniform spaces}

\hspace{15pt} The main result  of this second part of the thesis is the  identification of a ``universal space'' for Bourbaki-complete uniform spaces. Now, we start solving the particular case of embedding complete metric spaces and complete uniform spaces having a base for  their uniformity by means of star-finite open covers, that is, satisfying that $\mu=s_f\mu$. Observe that, in this case, we can take a base of uniform star-finite open covers for $s_{f}\mu$ (see \cite[Proposition 28, Chapter IV]{isbellbook}).
\smallskip

We recall that along this thesis the real line $\Rset$ is  endowed with the euclidean metric $d_u$ and any discrete space $D$ is endowed with the uniformly discrete metric  $$\chi(d , e)=\begin{cases}  \,\,\,  0 &     \text { if $d= e$}  \\ \,\,\, 1  & \text{ if $d \neq e$.}\end{cases}$$ For a countable product of discrete spaces $\prod_{n\in \Nset}D _n$, $\rho$ denotes the ``first difference metric" , that is,
$$\rho(( d_n)_n , ( e_n)_n)=\begin{cases}  \,\,\,  0 &     \text { if $d_n= e_n$ for every $n\in \Nset$}  \\ \,\,\, 1/n  & \text{ if $d_j = e_j$ for every $j=1,...,n-1$ and $d_n\neq e_n$}\end{cases}$$ which is compatible with the product uniformity on it. Finally, by $t$ we will denote the product metric on $\Rset^{\omega _0}$ $$t((x_n)_n ,(y_n)_n)=\sum _{n=1}^{\infty}(|x_n -y_n|\wedge 1)/2^n$$ and  by $\pi$ we will denote the product uniformity over any product of uniform spaces. 

\begin{definition} Let $(X,\mu)$ be a uniform space. A {\it uniform partition} of $(X,\mu)$ is a partition $\mathcal{P}$ of the space $X$ such that $\mathcal{U}<\mathcal{P}$ for some $\mathcal{U} \in \mu$ (that is, $\mathcal{P}\in \mu$).
\end{definition}

Observe that, whenever $\mathcal{U}\in \mu$, the family of all the chainable components $\{St^{\infty}(x_i, \mathcal{U}): i\in I\}$ induced by $\mathcal{U}$ is a uniform partition of $(X,\mu)$. In particular, if $\mathcal{U}<\mathcal{P}$ for some uniform partition $\mathcal{P}$ of $(X,\mu)$ then $$\{St^{\infty}(x_i, \mathcal{U}): i\in I\}<\mathcal{P}.$$

Next, let $(X,\mu)$ be a uniform space and let $\mathfrak{P}$ the family of all the uniform partitions of $(X,\mu)$. We define $$\wp(X,\mu)=sup \{|\mathcal{P}|:\mathcal{P}\in \mathfrak{P}\}.$$ For a connected uniform space, or in general, for a  uniformly connected  uniform space $(X,\mu)$, we have that $\wp(X,\mu)=1$. Recall that a {\it uniformly connected space} (or {\it well-chained} space) is a uniform space $(X,\mu)$ such that for every $\mathcal{U}\in \mu $, $X=St^{\infty}(x_{\mathcal{U}}, \mathcal{U})$ for some $ x_{\mathcal{U}}\in X$.

In the next theorem we apply techniques that can be found in \cite{balogh} and in \cite{garrido2}.

\begin{theorem} \label{embedding.star-finite1}Let $(X,d)$ be a complete metric space such that $\mu _d=s_f \mu _d$. 
Then there exists an embedding $$\varphi: (X,d)\rightarrow \Big((\prod _{n\in \Nset}\kappa_n)\times\Rset ^{\omega _0} , \rho +t\Big)$$ where each $\kappa _n$ is a cardinal endowed with the uniformly discrete metric $\chi$, $\varphi$ is uniformly continuous and $\varphi (X)$ is a closed subspace of  $(\prod _{n\in \Nset}\kappa_n)\times\Rset ^{\omega _0}.$ Moreover, $\wp(X,d)\leq sup\{\prod_{j=1}^n\kappa _j :n\in \Nset\}$.
\begin{proof}
Let us take $\{\mathcal{U}_n:n\in \Nset\}$ a family of star-finite open covers being a base for the metric uniformity $\mu _d$ and such that $\mathcal{U}_{n+1}<\mathcal{U}_n$ for every $n\in \Nset$. Without loss of generality assume that for every $n\in \Nset$, $\mathcal{U}_{n}$ refines $\mathcal{B}_{1/n}=\{B_{d}(x, 1/n):x\in X\}$. Next, for every $n\in \Nset$, let $\mathcal{P}_n=\{St^{\infty}(x_{i_n},\mathcal{U}_n): i_n \in I_n\}$ be the family of all the chainable components of $X$ induced by $\mathcal{U}_n$. Notice that $\mathcal{P}_{n+1}<\mathcal{P}_n$ for every $n\in \Nset$. Take the cardinal $\kappa _1=|I_1|$  and order the elements of the partition $\mathcal{P}_1$ by writing $$\mathcal{P}_1:=\{P_{(\alpha _1)}: \alpha _1<\kappa _{1}\} \text{ (where }\alpha _1 <\kappa_1 \text{ means }0\leq \alpha _1<\kappa _1 \text{)}.$$ Then, for every $\alpha _1 <\kappa _1$ let $\mathcal{P}_{(\alpha _1)}=\{P\in \mathcal{P}_2: P\subset P_{(\alpha _1)}\} $ and $\kappa _{(\alpha _1)}=|\mathcal{P}_{(\alpha _1)}|$. Next, put $\kappa _2=sup\{\kappa_{(\alpha _1)}:\alpha _1 <\kappa _1\}$. In particular, it is clear that $$\kappa _1 \times \kappa _2\geq |\{(\alpha _1,\alpha _2):\alpha _1<\kappa _1, \alpha _2<\kappa _{(\alpha _1)}\}|\geq |\bigcup_{\alpha_1 <\kappa _1} \mathcal{P}_{(\alpha _1)}|=|I_2|$$ Moreover, let us order each $\mathcal{P}_{(\alpha _1)}$ as follows: $$\mathcal{P}_{(\alpha _1)}:=\{P_{(\alpha _1, \alpha _2 )}: \alpha _2 <\kappa _{(\alpha _1)}\}.$$   

Next, suppose that for $n\in \Nset$ we have that the families of sets $\mathcal{P}_{(\alpha _1, \alpha _2,...,\alpha _{n-1})}=\{P\in \mathcal{P}_n: P\subset P_{(\alpha _1, \alpha _2,...,\alpha _{n-1})}\}$, where $\alpha _1 <\kappa _1$ and $\alpha_ j<\kappa_{(\alpha _1,...,\alpha _{j-1})}$ for every $j=2,...,n-1$, are defined and let $\kappa _{(\alpha _1, \alpha _2,...,\alpha _{n-1})}=|\mathcal{P}_{(\alpha _1, \alpha _2,...,\alpha _{n-1})}|$. And also, suppose that we have ordered  each family $\mathcal{P}_{(\alpha _1, \alpha _2,...,\alpha _{n-1})}$, as follows: $$\mathcal{P}_{(\alpha _1, \alpha _2,...,\alpha _{n-1})}:=\{P_{(\alpha _1,\alpha _2,...,\alpha _{n-1}, \alpha _n)}:\alpha _n <\kappa _{(\alpha _1, \alpha _2,...,\alpha _{n-1})}\}.$$ 

Then, by induction, we put $\mathcal{P}_{(\alpha _1, \alpha _2,...,\alpha _{n})}=\{P\in \mathcal{P}_{n+1}:P\subset P_{(\alpha _1,\alpha _2,...,\alpha _n})\} $ and $\kappa _{(\alpha _1, \alpha _2,...,\alpha _{n})}=|\mathcal{P}_{(\alpha _1,\alpha _2,...,\alpha _n})|$. Finally we order each family of sets $\mathcal{P}_{(\alpha _1, \alpha _2,...,\alpha _{n})}$  as before: $$\mathcal{P}_{(\alpha _1, \alpha _2,...,\alpha _{n})}:=\{P_{(\alpha _1,\alpha _2,...,\alpha _{n}, \alpha _{n+1})}:\alpha _{n+1} <\kappa _{(\alpha _1, \alpha _2,...,\alpha _{n})}\}.$$ 

In addition, for every $n\in \Nset$ we put $$\kappa _{n+1}=sup\{\kappa _{(\alpha _1,\alpha _2,...,\alpha _n)}: \alpha _1 <\kappa_1,\alpha_ j<k_{(\alpha _1,...,\alpha _{j-1})}, j=2,...,n\}.$$ Observe that $$\kappa_1 \times \kappa _2 \times ...\times \kappa_n \times \kappa _{n+1}\geq $$ $$|\{(\alpha _1,\alpha _2,...,\alpha _n,\alpha _{n+1}): \alpha _1 <\kappa_1,\alpha_ j<k_{(\alpha _1,...,\alpha _{j-1})}, j=2,...,n+1\}\geq$$

$$|\bigcup \{ \mathcal{P}_{(\alpha _1,...,\alpha _{n}}): \alpha _1 <\kappa_1,\alpha_ j<k_{(\alpha _1,...,\alpha _{j-1})}, j=2,...,n\}|=|I_{n+1}|.$$ Moreover, as $\{\mathcal{U}_n:n\in \Nset\}$ is a base for the uniformity $\mu_d$, $\wp (X,d)\leq sup\{|I_n|:n\in \Nset\}\leq sup\{\prod_{j=1}^n\kappa _j :n\in \Nset \}$. 

Next, notice that for every $n\in \Nset$ there exists a unique $(\alpha _1,...,\alpha _n)\in \prod_{j=1}^n \kappa _j$ such that $x\in \mathcal{P}_{(\alpha _1,...,\alpha _n)}$. Besides, $(\alpha _1,...,\alpha _n,\alpha _{n+1})$ extends $(\alpha _1,...,\alpha _n)$, so there exists a unique $\sigma(x)\in \prod_{n\in \Nset}\kappa _n$ such that the restriction $\sigma(x)|n$ of $\sigma(x)$ over the first $n$'s coordinates is exactly $(\alpha _1,...,\alpha _n)$.  Let us denote by $\sigma: (X,d)\rightarrow (\prod_{n\in \Nset}\kappa _n, \rho)$ the map that sends $x$ to $\sigma(x)$.  Recall, that for every $x\in X$ and every $n\in \Nset$ there exists a unique $i_n \in I_n$ such that $P_{\sigma (x)|n}=St^{\infty}(x_{i_n},\mathcal{U}_n)= P_{\sigma (x_{i_n})|n}$.

Now, for every $n\in \Nset$ let $\mathcal{A}(\mathcal{U}_n)=\{ A_{m,i _n}:m\in \Nset,  i _n\in I _n\}$ the open cover from Lemma \ref{corona} induced by $\mathcal{U}_n$, and define the sets $A_m ^n =\bigcup \{A_{m, i_n}:i_n \in I_n\}$, $m\in \Nset$. Then the cover $\mathcal{A}_n=\{A_m ^n:m\in \Nset\}$ is uniform and linear. Take  $\varepsilon _n >0$ such that $\{B_{d}(x,\varepsilon _n):x\in X\}<\mathcal{U}_n<\mathcal{A}_n $. Applying the same techniques than in \cite[Lemma 1.2]{garrido-montalvo} there  exists a uniformly continuous function $h_n: (X,d)\rightarrow (\Rset, d_u)$ such that $h_{n}^{-1}((m-1, m+1))=A_{m}^n$ for every $m\in \Nset$. Moreover, the following is always satisfied (see \cite[Lemma 1.2]{garrido-montalvo}):

\begin{center} if $d(x,y)\leq\varepsilon _n$ we have that $|h_{n}(x)- h_{n}(y)|\leq (10/\varepsilon _n ^2)\cdot d(x,y)$ 
\end{center}

Next, recall that we can write $\mathcal{U}_n=\{U_{j,i_n}, j\in \Nset, i_n \in I _n\}$ where $U_{j,i _n}\cap U_{j',i' _n}=\emptyset $ if $i _n\neq i' _n$ as the covers $\mathcal{U} _n$ are star-finite and every chainable component $\mathcal{P}_{\sigma(x_{i_n})|_n}$, $i_n\in I_n$, contains at most countable many $U\in \mathcal{U}_n$.
\smallskip

Now, since  for every $x\in X$ and for every $n\in \Nset$, there exists a unique $i_n \in I_n$ such that $\sigma (x)|_n=\sigma (x_{i_n})|n$, we define the map \begin{align*}\varphi: (X,\mu)&\rightarrow \Big((\prod_{n\in \Nset} \kappa _n )\times (\Rset \times \Rset^{\omega _0})^{\omega_0}, \pi\Big) \\
x&\mapsto\varphi(x)=\big(\sigma (x), \langle h_{n}(x), \langle d(x,X\backslash U_{j,i_n})\rangle_{j\in \Nset}\rangle_{n\in \Nset}\big) .\end{align*}  

$\bullet$ \underline{The map $\varphi$ is injective.} Let $x,y \in X$, $x\neq y$. Then we can take some $\varepsilon <d(x,y)$ such that $y\notin B_{d} (x,\varepsilon)$. Since $\{\mathcal{U}_n:n\in \Nset\}$ is a base for the uniformity inducing the topology on $X$, for some $n\in \Nset$ we can choose $U_{j,i _n}\in \mathcal{U}_n$ such that  $x\in U_{j,i _n}\subset B_{d}(x,\varepsilon )$.  Then $d(x, X\backslash U_{j,i _n})>0$ and $d(y, X\backslash U_{j,i _n})=0$. Therefore $\varphi(x)\neq \varphi(y)$ and $\varphi$ is an injective map.

$\bullet$ \underline{The map $\varphi$ is uniformly continuous.} We check that the map $\varphi$ is uniformly continuous by showing that it is uniformly continuous when we compose it with the projections. First, it is easy to see that $\sigma$ is a uniformly continuous map since whenever $d(x,y)<\varepsilon _n$  then $P_{\sigma(x)|n}=P_{\sigma(x_{i_n})|n}=P_{\sigma(y)|n}$ for a unique $i_n\in I_n$. Therefore, $\sigma(x)|n=\sigma(y)|n$ and then $\rho(x,y)<\frac{1}{n+1}$. Next, let $d(x,y)<\varepsilon _n$ again, then $$d_u(h_{n}(x),h_{n}(y)) +t( \langle d(x,X\backslash U_{j,i_n})\rangle_{j\in \Nset},\langle d(y,X\backslash U_{j,i_n})\rangle_{j\in \Nset})=$$  $$|h_{n}(x)-h_{n}(y)|\ +\sum_{j=1}^{\infty}|d(x,X\backslash U_{j,i _n })-d(y,X\backslash U_{j,i _n })|/2^j\leq$$ $$|h_{n}(x)-h_{n}(y)|\ +\sum _{j=1}^{\infty}d(x,y)/2^j \leq  (10/\varepsilon _n ^2)\cdot d(x,y)+d(x,y)=(10/\varepsilon _n ^2 +1)\cdot d(x,y).$$

$\bullet$ \underline{The map $\varphi$ is closed.} Before proving that $\varphi$ is closed, we need to prove the following claim.

{\bf CLAIM.} Let $Y\subset \varphi (X)$ and $\mathcal{F}$ a Cauchy filter of the subspace $(Y,\pi |_Y)$. Then $\varphi ^{-1}(\mathcal{F})$ is a Bourbaki-Cauchy filter of $(X,d)$.

\begin{proof}[Proof of the claim] Let $\mathcal{F}$ be Cauchy filter of $(Y, \pi |_{Y})$. Then, fixed $k\in \Nset$, since $\mathcal{F}$ is Cauchy, there is some $W\in \mathcal{F}$, $W=\big(V\times(\prod _{n\in \Nset} U_n)\big)\cap Y$, where, for some $x_0\in X$ and some $i_k \in I_k$, $V=B_{\rho}(\sigma (x_0),1/k)$, $\sigma(x_0)|k=\sigma(x_{i_k})|k$ and $U_{k}=B_{d_u}(( h_{k}(x_0),1/k)\times \Rset ^{\omega _0}$,  and $U_n= \Rset \times \Rset^{\omega _0} $ for every $n\neq k$.

As the fixed $i_k\in I_k$ such that $\sigma(x_0)|k=\sigma (x_{i_k})|k$ is unique, then  
$$\varphi ^{-1}(W)=\{x\in P_{\sigma(x_{i_k})|k}:|h_{k}(x)-h_{k}(x_0)|<1/k\}\subset$$ $$ h_{i _k}^{-1}((h_{k}(x_0)-1/k, h_{k}(x_0)+1/k))\cap P_{\sigma(x_{i_k})|k}.$$ By the construction of $h_{k}$ there is some $m\in \Nset$ such that $$ h_{k}^{-1}((h_{k}(x_0)-1/k, h_{k}(x_0)+1/k))\cap P_{\sigma(x_{i_k})|k} \subset A_{m,i _k}\subset St^{m+1}(x_{i_k},\mathcal{U}_k).$$ Therefore $$\varphi^{-1}(W)\subset St^{m+1}(x_{i _k}, \mathcal{U}_k) $$ and we have proved that $\varphi ^{-1}(\mathcal{F})$ is a Bourbaki-Cauchy filter in $(X,d)$.
\end{proof}

Now we follow with the proof that $\varphi$ is a closed map. Let $C\subset X$ a closed subset and let $\mathcal{F}$ be an ultrafilter in $\varphi (C)$ which converges to some $z \in (\prod_{n\in \Nset}\kappa _n) \times (\Rset \times \Rset ^{\omega _0} )^{\omega _0}$.  Then $\mathcal{F}$ is a Cauchy ultrafilter of the subspace $(\varphi(C), \pi|_{\varphi(C)})$. By maximality of $\mathcal{F}$ and the above claim, $\varphi ^{-1}(\mathcal{F})$ is a Bourbaki-Cauchy ultrafilter of $(X,d)$. Then, $\varphi ^{-1}(\mathcal{F})$ converges in $X$ because, by Theorem \ref{star-finite2}, $(X,d)$ is in particular Bourbaki-complete. By continuity of $\varphi$, $\varphi (\varphi ^{-1}(\mathcal{F}))$ converges in $\varphi (X)$. Since $\mathcal{F}= \varphi (\varphi ^{-1}(\mathcal{F}))$,  by maximality, $\mathcal{F}$ converges in $\varphi(X)$, that is, $z\in \varphi(C)$. Thus,  $\varphi(C)$ is a closed subspace of $(\prod_{n\in \Nset}\kappa _n )\times (\Rset \times \Rset ^{\omega _0} )^{\omega _0}$.

$\bullet$ \underline{The image $\varphi(X)$ is a closed subspace of $(\prod_{n\in \Nset}\kappa _n )\times (\Rset \times \Rset ^{\omega _0} )^{\omega _0}.$} Since $\varphi$ is a closed map, this is clear.
\smallskip

Finally,  observe that, by the results in Bourbaki \cite[II.2.3 Prop 5 p. 177, p.180]{bourbaki}, the spaces $\Big((\prod_{n\in \Nset}\kappa _n )\times (\Rset \times \Rset ^{\omega _0} )^{\omega _0},\pi \Big)$  and $\Big((\prod_{n\in \Nset}\kappa _n) \times \Rset^{\omega _0},\rho+t\Big)$ are uniformly equivalent, and this complete the proof.

\end{proof}
\end{theorem}

Now we solve the analogous problem for uniform space having a base of star-finite covers for the uniformity.

\begin{lemma}\label{pseudometric} {\rm (\cite{willard})} Let $\langle\mathcal{G}_n \rangle _{n\in \Nset}$ be a normal sequence of open (uniform) covers of a (uniform) space $X$. Then there exists a (uniformly) continuous pseudometric $d:X\times X\rightarrow [0,\infty)$ such that $$\mathcal{B}_{1/2^{n+1}}<\mathcal{G}_n<\mathcal{B}_{1/2^{n-1}} \text{ for every }n\in \Nset.$$ 
\end{lemma}

\begin{theorem} \label{embedding.star-finite2}Let $(X,\mu)$ be a complete uniform space such that $\mu =s_f\mu$.   Then there exists an embedding $$\varphi: (X,\mu)\rightarrow \Big((\prod_{ i\in I, n\in \Nset} \kappa _n ^i)\times \Rset ^{\alpha},\pi\Big) $$ where each $\kappa _n ^i$ is a cardinal endowed with the uniformly discrete metric $\chi$, $\alpha \geq \omega _0$, 
$\varphi$ is uniformly continuous and $\varphi (X)$ is a closed subspace of $(\prod_{ i\in I, n\in \Nset} \kappa _n ^i)\times \Rset ^{\alpha}.$ Moreover,  $\wp(X,\mu)\leq sup \{(\prod_{j=1} ^n\kappa _{j} ^i):n\in \Nset, i\in I\}$.
\begin{proof} Let $(X,\mu)$ be a uniform space having a base of star-finite open covers for the uniformity, and let $\{\mathcal{U}^i:i\in I\}$ be a base for $\mu$. Then, for every $\mathcal{U}^{i}\in \mu$, $i\in I$ there is a normal sequence $\langle\mathcal{U}_n ^i\rangle_{n\in \Nset}$ of star-finite uniform open covers such that $\mathcal{U}_1 ^i<\mathcal{U}^i$. This can be obtained applying the axioms of uniformity and \cite[Proposition 8, Chapter IV]{isbellbook}.

For every $i\in I$, let $\rho_i$ be the pseudometric on $X$ from Lemma \ref{pseudometric} generated by the normal sequence $\langle \mathcal{U}_n ^i\rangle _{n\in \Nset}$. Then,  the family of covers $\{\mathcal{U}_n ^i:n\in \Nset \}$ is a base of star-finite open covers of the space $(X,\rho _i)$.  

Let $(Y_i,\widehat{\rho _i})$ be the metric space obtained by doing the usual metric identification $\sim$ on $(X,\rho_{i})$:

\begin{center} $x_1\sim x_2$ if and only if $\rho _{i}(x_1,x_2)=0.$\end{center} If $\psi _i:(X,\rho_{i})\rightarrow (Y_i,\widehat{\rho_i})$ denotes the quotient map induced by $\sim$, then $\psi_i ^{-1}(\widehat{x})=\{z\in X: \rho_{i}(x,z)=0\}$, $\widehat{A}:=\psi_i (A)=\{\widehat{x}:x\in A\}$ and $\psi_i^{-1}(B_{\widehat{\rho _i}}(\widehat{x}, \varepsilon))=B_{\rho_i}(x,\varepsilon)$. Hence,  the family of covers $\{\widehat{\mathcal{U}}_n ^i:n\in \Nset \}$ is  a base of star-finite  covers for the metric uniformity on $Y_i$ induced by $\widehat{\rho _i}$. In addition, the map $\psi_i$ preserves the uniform partitions induced by the covers $\mathcal{U}_n ^i$, $n\in \Nset$.

Let $(Z_i, d_i)$ denote the completion of $(Y_i, \widehat{\rho _i})$ and $\mathcal{V} _n  ^i$ denote the  extension to $(Z_i, d_i)$ of the covers $\widehat{\mathcal{U}} _n ^i$.  Then $\{\mathcal{V}_{n} ^i:n\in \Nset\}$ is a base of star-finite open covers for the metric uniformity of $(Z_i,d _i)$ (\cite[Lemma p. 370]{reynolds}). 

By \cite[Theorems 39.11 and  39.12]{willard} $(X,\mu)$ is uniformly homeomorphic to subspace of the product $\prod _{i\in I}(Z_{i}, d _{i})$. In particular it is closed by completeness. Denote by $\varphi _{i}$ the embedding of $(Z_i, d_i)$ into $\Big((\prod_{n\in \Nset}\kappa _n ^i )\times \Rset^{\omega _0}, \rho+t\Big)$ from Theorem \ref{embedding.star-finite1}, and let $$\varphi: \prod_{i\in I}(Z_i,d_{i})\rightarrow \prod_{i\in I}\Big((\prod_{n\in \Nset}\kappa _n ^i )\times \Rset^{\omega _0}, t+\rho\Big)$$ be the product map $\varphi=\prod_{i\in I} \varphi _{i}$. Then, the restriction of $\varphi$ over the uniform homeomorphic image of $(X,\mu)$ in $\prod_{i\in I}(Z_i,d_{i})$ is the desired map. Indeed, notice that $\varphi (X)$ is closed in $\varphi (\prod_{i\in I} Z_i)=\prod_{i\in I} \varphi _i (Z_i)$. Moreover,  by \cite[II.2.3 Prop 5 p.177, p.180]{bourbaki}, the spaces $\prod_{i\in I}\Big((\prod_{n\in \Nset}\kappa _n ^i )\times \Rset^{\omega _0}\Big)$ and $(\prod_{ i\in I, n\in \Nset}\kappa _n ^i )\times \Rset^{\alpha}$ for $\alpha =sup\{ |I|,\omega _0\}$ are uniformly equivalent when they are  endowed with their respective product uniformities. Finally, by Theorem \ref{embedding.star-finite1}, $$\wp(X,\mu)\leq sup \{\wp(X,\rho_i):i\in I\}$$ $$= sup \{\wp(Z_i,d_i);i\in I\}\leq sup\{\prod_{j=1}^n\kappa_j ^i: n\in \Nset, i\in I\}$$ as the quotient map $\psi _i$ and the operation of completion preserve uniform partitions.
\end{proof}
\end{theorem}

Next, recall that by Theorem \ref{star-finite2}, $(X,s_f\mu)$ is complete if and only if $(X,\mu)$ is Bourbaki-complete. Therefore, if we compose the embedding $$\varphi: (X,s_f \mu) \rightarrow \Big((\prod_{i\in I, n\in \Nset}\kappa _n ^i)\times \Rset ^{\alpha} ,\pi\Big)$$ from Theorem \ref{embedding.star-finite2}, with the identity map $id:(X,\mu)\rightarrow (X, s_f\mu)$ we have that the following result, characterizing a universal space for Bourbaki-complete uniform spaces, is immediate. 

\begin{theorem} \label{universal.uniform} Let $(X,\mu)$ be a Bourbaki-complete uniform space. Then there exists  an embedding $$\varphi: (X,\mu)\rightarrow \Big((\prod_{i\in I, n\in \Nset}\kappa _n ^i)\times \Rset ^{\alpha} ,\pi\Big),$$ where each $\kappa _n ^i$ is a cardinal endowed with the uniformly discrete metric $\chi$, $\alpha \geq \omega _0$,  $\varphi$ is uniformly continuous and $\varphi (X)$ is a closed subspace of $(\prod_{i\in I, n\in \Nset}\kappa _n ^i)\times \Rset ^{\alpha}.$ Moreover,  $\wp(X,\mu)\leq sup \{(\prod_{j=1} ^n\kappa _{j} ^i):n\in \Nset, i\in I\}$.
\end{theorem}

\begin{remark}\label{0-dim} As we have said at the beginning of this section, whenever $(X,\mu)$ is a connected or uniformly connected space then every uniform partition has cardinal 1  and hence, in the above embeddings, the discourse on the chainable components is clearly not needed, and we can straightly embed our Bourbaki-complete uniform space in a product of real-lines. 

On the other hand, we have a class of uniform spaces which represents the opposite situation. Recall that a uniform space  is {\it uniformly 0-dimensional} if the uniformity has a base composed of partitions. Observe that from Theorem \ref{embedding.star-finite2} any complete uniformly 0-dimensional space can be uniformly embedded, as closed subspace, in a product of uniformly discrete spaces. Indeed, the embedding is given by the map $\sigma:(X,\mu)\rightarrow (\prod _{i\in I, n\in \Nset} \kappa_n ^i, \pi)$ in  Theorem \ref{embedding.star-finite2}. We just need to see that the inverse map $\sigma ^{-1}$ from Theorem \ref{embedding.star-finite1} is uniformly continuous. Observe that, since $(X,d)$ is uniformly 0-dimensional the family of all the chainable components $\{\mathcal{P}_n:n \in \Nset\}$ is a base for the uniformity of $(X,d)$. Fix $n\in \Nset$, suppose that for $x,y\in X$, $\rho(\sigma(x),\sigma(y))<1/(n+1)$. Then $\sigma(x)|n=\sigma(y)|n$ and this implies that $x,y$ belong to the same chainable component of $\mathcal{P}_k$ for every $k=1,...,n$.  Therefore we can conclude that $\sigma ^{-1}$ is uniformly continuous. Finally, the general case for uniformly 0-dimensional uniform spaces proceeds like in Theorem \ref{embedding.star-finite2} and taking into the account that the product map of uniformly continuous functions is uniformly continuous.

\end{remark}


\subsection{Embedding Bourbaki-complete metric spaces}

\hspace{15pt} By Theorem \ref{universal.uniform}, it is clear that a universal space for Bourbaki-complete metric spaces is also determined. However, in the metrizable case, by future technical reasons, we need that this universal space is also metrizable. Indeed, observe that for a Bourbaki-complete metric space $(X,d)$, the star-finite modification $s\mu _d$ is not metrizable in general and therefore the universal space obtained from Theorem \ref{universal.uniform} is not metric. In particular we have the following result.

\begin{theorem} Let $(X,d)$ be a Bourbaki-complete metric space. Then the uniform space $(X,s_f\mu _d)$ is metrizable if and only if $s_f\mu _d=\mu _d$.
\begin{proof} One implication is clear, so let $\rho$ be a metric on $X$ such that $\mu _{\rho}=s_f\mu _d$. In particular $(X,\rho)$ is complete. By Lemma \ref{Samuel2}, the Samuel compactification $s_d X$ coincides with the Samuel compactification $s_\rho X$. Therefore, by \cite[Corollary 3]{garrido.Banach-Stone}, it is known that  $(X,d)$ and $(X,\rho)$ are  uniformly homeomorphic. Hence, $s_f\mu _d=\mu _\rho=\mu _d$. 
\end{proof}
\end{theorem}

In order to find a universal metric space for Bourbaki-complete metric spaces we will need first some technical results.

\begin{definition} (\cite{engelkingbook}) A sequence $\langle\mathcal{A}_n\rangle _{n\in \Nset}$ of open covers of a topological space $X$ is a {\it complete sequence of covers} if, for every filter $\mathcal{F}$ of $X$ satisfying that $\mathcal{F}\cap \mathcal{A}_n \neq \emptyset$ for every $n\in \Nset$, then $\mathcal{F}$ has a cluster point.
\end{definition}

\begin{definition} A sequence of covers $\langle\mathcal{C}_n\rangle _{n\in \Nset}$ of a set $X$ is a {\it decreasing sequence of covers} if for every $n\in \Nset$, $\mathcal{C}_{n+1}<\mathcal{C}_n$  and for each $C\in \mathcal{C}_n$, we have that $C=\bigcup\{C'\in \mathcal{C}_{n+1}:C'\subset C\}$.
\end{definition}

For (open) covers $\mathcal{G}_1, \mathcal{G}_2,...,\mathcal{G}_n$ of a (space) set $X$, we denote by $\mathcal{G}_1\wedge\mathcal{G}_2\wedge...\wedge \mathcal{G}_n$ the (open) cover $\{G_1\cap G_2\cap ...\cap G_n:G_{i}\in \mathcal{G}_i, i=1,2,...,n\}$. In particular $\mathcal{G}_1 \wedge \mathcal{G}_2\wedge ... \wedge \mathcal{G}_n $ refines $\mathcal{G}_i$ for each $i=1,...,n$. More precisely, if $\langle \mathcal{A}_n\rangle _{n\in \Nset}$ is a sequence of covers of $X$ and for every $n\in \Nset$, we define $\mathcal{C}_n=\mathcal{A}
_1\wedge ...\wedge \mathcal{A}_n$ then, $\langle\mathcal{C}_n\rangle _{n\in \Nset}$ is a decreasing sequence of covers.

\begin{lemma}\label{complete.sequence} {\rm (\cite[Lemma 2.8]{merono.hh})} Let $ \langle\mathcal{U}_n\rangle _{n\in \Nset}$ be a decreasing complete sequence of a topological space $X$. Then the family $\bigcup _{n\in \Nset}\mathcal{U}_n$ contains a refinement of every directed open cover of $X$.

\begin{proof} Let $\mathcal{G}$ be a directed open cover of $X$. We show that the family $\mathcal{W}=\{U\in \bigcup _{n\in \Nset}\mathcal{U}_n: U\subset G \text{ for some } G\in \mathcal{G}\}$ covers $X$. Assume on the contrary that there exists a point $x$ in the set $X\backslash \bigcup \{W: W\in \mathcal{W}\}$. Since $\langle\mathcal{U}_n\rangle_{n\in \Nset}$ is a decreasing sequence of covers of $X$, there exists $U_1, U_2,...,U_n,...$ such that $x\in U_{n+1}\subset U_n \in \mathcal{U}_n$ for every $n$. Since $x\notin \bigcup \{W:W\in \mathcal{W}\}$ none of the sets $U_1,U_2,...,U_n,...$ is contained in any member of $\mathcal{G}$. It follows that the family $\mathcal{L}=\{U_n\backslash G:n\in \Nset \text{ and }G\in \mathcal{G}\}$ is a filter base. Let $\mathcal{F}$ be the filter of $X$ generated by $\mathcal{L}$. We have $U_n\in \mathcal{F}$ for every $n$ and it follows, since $\langle\mathcal{U}_n\rangle _{n\in \Nset}$ is a complete sequence, that $\mathcal{F}$ has a cluster point $z$. This, however, is impossible: there exists $G\in \mathcal{G}$ with $z\in G$ and now we have $U_1\backslash G \in \mathcal{F}$ but $z\notin {\rm cl}_X (U_1\backslash G)$. This contradiction shows that $\mathcal{W}$ covers $X$. As a consequence, $\mathcal{W}$ is a refinement of $\mathcal{G}$ contained in $\bigcup _{n\in \Nset}\mathcal{U}_n$.
\end{proof}

\end{lemma}

Observe that the next result is a metric and uniform extension of \cite[Theorem 2.16]{merono.hh}.

\begin{theorem} \label{Bourbaki-complete.strongly-metrizable} Let $(X,d)$ be a Bourbaki-complete metric space. Then there exists a complete sequence $\langle \mathcal{V}_n\rangle_{n\in \Nset}$ of uniform star-finite open covers of $(X,d)$ such that $\bigcup \mathcal{V}_n$ is a base of the topology of $X$.  
\begin{proof}  For every $n\in \Nset$, let $\mathcal{A}_n:=\mathcal{A}(\mathcal{B}_{1/n})=\{A_{m,i _n}:m\in \Nset, i_n\in \Nset\}$ the cover from Lemma \ref{corona} induced by the cover of open balls $\mathcal{B}_{1/n}=\{B_{d}(x, 1/n):x\in X\}$.

Next, define $\mathcal{U}_n=\mathcal{A}_1 \wedge \mathcal{A}_2 \wedge ...\wedge \mathcal{A}_n$, $n\in \Nset$. It is clear that $\langle \mathcal{U}_n\rangle_{n\in \Nset}$ is a decreasing sequence of star-finite uniform open covers of $X$ since finite intersection of star-finite uniform open covers is again star-finite, open and uniform. We prove now that $\langle \mathcal{U}_n\rangle_{n\in \Nset}$ is a complete sequence. Let $\mathcal{F}$ be a filter in $X$ such that for every $n\in \Nset$ there exists some $U\in \mathcal{U}_n$ such that $F\subset U$ for some $F\in \mathcal{F}$.  In particular, $\mathcal{F}$ is a Bourbaki-Cauchy filter. Indeed, if $U\in \mathcal{F}$ for some $U\in\mathcal{U}_n$, then $$U\subset A_{m ,i _n}\subset B^{m+1}_{d}(x_{i _n}, 1/n)$$ for some $m\in \Nset$ and $i _n\in I _n$. Therefore, $\mathcal{F}$  clusters in $X$ and $\langle \mathcal{U}_n\rangle_{n\in \Nset}$ is a complete sequence.

Next, let $\mathcal{G}$ be an open cover of $X$ and $\mathcal{G}^f$ the directed open cover given by finite unions of elements of $\mathcal{G}$. By all the foregoing and by Lemma \ref{complete.sequence}, there exists a cover $\mathcal{U}\subset \bigcup_{n\in \Nset}\mathcal{U}_n$ such that $\mathcal{U}<\mathcal{G}^{f}$. Now, for every $U\in \mathcal{U}$ fix $\mathcal{G}_{U}$ a finite subfamily of $\mathcal{G}$ such that $U\subset \bigcup \mathcal{G}_U$. Note that for each $n\in \Nset$, the family $\mathcal{U}_n(\mathcal{G})=\mathcal{U}_n \cup \{U\cap G: U\in \mathcal{U}\cap\mathcal{U} _n \text{ and } G\in \mathcal{G}_U\}$ is a star-finite open cover of $X$. Moreover, it is also uniform since $\mathcal{U}_n\subset \mathcal{U}_n (\mathcal{G})$.  Therefore, the cover $\mathcal{G}$ has a refinement which is contained in $\bigcup _{n\in \Nset} \mathcal{U}_n(\mathcal{G})$.

Finally, let $f:\Nset \times \Nset \rightarrow \Nset$ be any bijection, and, for every $n,j\in \Nset$, define the covers $\mathcal{V}_{f((n,j))}=\mathcal{U}_n(\mathcal{B}_{1/j})$. We check now that $\bigcup _{(n,j)\in \Nset\times \Nset} \mathcal{V}_{f((n,j))}$ is a base for the topology of $X$. Indeed, let $G$ be an open set of $X$ and $x\in G$. Then we can choose some $k\in \Nset$ such that $x\in B_{d}(x,1/k)\subset G$. Consider the open cover of balls $\mathcal{B}_{1/2k}$. By all the foregoing, $\mathcal{B}_{1/2k}$ has a refinement contained in $\bigcup _{n\in \Nset} \mathcal{U}_n(\mathcal{B}_{1/2k})=\bigcup _{n\in \Nset} \mathcal{V}_{f((n,2k))}$. Choose $V\in \bigcup _{n\in \Nset} \mathcal{V}_{f((n,2k))}$ such that $x\in V$. Then $x\in V\subset B_{d}(y, 1/2k)$ for some $y\in X$. Since $y\in B_{d}(x, 1/2k)$ then $$x\in V\subset B_{d}(y, 1/2k)\subset B_{d}(x, 1/k)\subset G.$$ Thus, we conclude that $\bigcup _{(n,j)\in \Nset\times \Nset} \mathcal{V}_{f((n,j))}$ is a base for the topology of $X$.

\end{proof}
\end{theorem}

\begin{theorem} \label{metrizable.star-finite} Let $(X,d)$ be a Bourbaki-complete metric space. Then there exists a  complete metric $d'$ on $X$ which is compatible with the topology of $X$ such that the metric uniformity $\mu _{d'}$ has a base of star-finite covers and $\mu _d \geq \mu _{d'} $. Moreover, $\wp (X,d)=\wp(X,d')$.
\begin{proof}  By Theorem \ref{Bourbaki-complete.strongly-metrizable}, there exists a complete sequence $\langle \mathcal{V}_n\rangle _{n\in \Nset}$ of uniform star-finite open covers of $(X,d)$ such that $\bigcup \mathcal{V}_n$ is a base of the topology of $X$. Observe that we can take a complete normal sequence $\{\mathcal{W} _n:n\in \Nset\}$  of open covers from $s_f \mu _d $ such that  $\mathcal{W}_{n} <\mathcal{V}_n$ for every $n\in \Nset$ and $\bigcup_{n\in \Nset} \mathcal{W}_n$ is a base for the topology of $X$. Indeed, take  $\mathcal{W}_1=\mathcal{V}_1$. By \cite[Proposition 8, Chapter IV]{isbellbook}, we can take an open cover $\mathcal{A}_1\in s_f \mu_d$ such that $\mathcal{A}_1 ^{*}<\mathcal{W}_1$. Put $\mathcal{W}_2=\mathcal{V}_1\wedge \mathcal{V}_2 \wedge \mathcal{A}_1$. Then $\mathcal{W}_2 ^{*} <\mathcal{W}_1$, $\mathcal{W}_2<\mathcal{V}_2$ and clearly, $\mathcal{W}_2$ is a uniform star-finite open cover. Again, by \cite[Proposition 8, Chapter IV]{isbellbook}, we can take an open cover $\mathcal{A}_2\in s_f \mu _d$ such that $\mathcal{A}_2 ^* <\mathcal{W}_2$. Put $\mathcal{W}_3=\mathcal{V}_1\wedge \mathcal{V}_2 \wedge \mathcal{V}_3 \wedge \mathcal{A}_2$. Thus, proceeding by induction we obtain such a normal sequence.

By Lemma \ref{pseudometric}, there exists a uniformly continuous pseudometric $d'$ on $X$ such that $(X, d')$ has a base of star-finite open covers for the pseudometric uniformity. Moreover, $(X, d')$ is complete because the sequence $\langle \mathcal{W}_n\rangle_{n\in \Nset}$ is complete.  Since $\bigcup_{n\in \Nset} \mathcal{W}_n$ is a base for the topology of $X$, then $d'$ is compatible with the topology of $X$. More precisely, since $X$ is Hausdorff, then $d'$ is a metric.  Finally, $\mu _d \geq \mu _{d'} $ since $\mathcal{W}_n\in s_f \mu _d$ for every $n\in \Nset$.

Next, we check that $\wp(X,d)=\wp(X,d')$. Since $\mu_d\geq \mu_{d'}$ then $\wp(X,d)\geq\wp(X,d')$. Therefore, we just need to prove that $\wp(X,d)\leq\wp(X,d')$. To that purpose, let us denote by $\mathcal{P}_n$, $n\in \Nset$, the families of all the chainable components induced by the covers $\{B_d(x, 1/n): x\in X\}$,  and by $\mathcal{Q}_m$, $m\in \Nset$, the family of all the chainable components induced by the above covers $\mathcal{W}_m$. Observe that, if we want to prove that $\wp(X,d)\leq\wp(X,d')$, it is enough to see that for every $n\in \Nset$, there exists $m\in \Nset$ such that $\mathcal{Q}_m<\mathcal{P}_n$, that is, $|\mathcal{P}_n|\leq |\mathcal{Q}_m|$. 

Fix $n\in \Nset$ and consider $\mathcal{P}_n$. Looking into the end of the proof of Theorem \ref{Bourbaki-complete.strongly-metrizable}, there exists some $m\in \Nset$ such that the cover $\mathcal{V}_m:=\mathcal{V}_{f((n,j))}$, from the beginning of this proof, induces the same chainable components than the cover $\{B_d(x, 1/n): x\in X\}$, $n\in \Nset$, that is, precisely the family $\mathcal{P}_n$. Since $\mathcal{W}_m<\mathcal{V}_m$. Then, it is clear that $\mathcal{Q}_m <\mathcal{P}_n$ as we claimed.
\end{proof}
\end{theorem}

\begin{theorem}\label{embedding.metric} Let $(X,d)$ be a Bourbaki-complete metric space. Then, there exists an embedding $$\varphi: (X,d)\rightarrow \Big((\prod_{n\in \Nset}\kappa _n )\times \Rset ^{\omega _0} , \rho +t\Big)$$ where each $\kappa _n$ is a cardinal endowed with the uniformly discrete metric $\chi$, $\varphi$ is uniformly continuous and $\varphi (X)$ is a closed subspace of  $(\prod_{n\in \Nset}\kappa _n )\times \Rset ^{\omega _0}.$ Moreover, $\wp(X,d)\leq sup\{\prod_{j=1}^n\kappa _j :n\in \Nset\}$.
\begin{proof}  By Theorem \ref{metrizable.star-finite} there exists a compatible metric $d'$ on $X$ such that $(X,d')$ is complete, the metric uniformity $\mu _{d'}$ has a star-finite base, the identity map $id:(X,d)\rightarrow (X,d')$ is uniformly continuous and $\wp(X,d)=\wp(X,d')$. On the other hand, consider the embedding $$\varphi: (X,d')\rightarrow \Big((\prod_{n\in \Nset}\kappa _n )\times \Rset ^{\omega _0} , \rho +t \Big)$$ form Theorem \ref{embedding.star-finite1}. Then the composition $\varphi \circ id=\varphi$ is the desired embedding. Finally, as $\wp(X,d)=\wp(X,d')$, by Theorem \ref{embedding.star-finite1} it follows that $\wp(X,d)\leq sup\{\prod_{j=1}^n\kappa _j :n\in \Nset\}$.

\end{proof}
\end{theorem}

\medskip

\begin{remark} Next, we would like to know if the above embedding $\varphi$ of $X$ into $ (\prod _{n\in \Nset}\kappa _n )\times \Rset ^{\omega _0}$ could be stronger. More precisely, we want to know if $\varphi$ preserves ``proximities", that is, if the Samuel compactification $s_d X$ of $(X,d)$  is homeomorphic to the Samuel compactification of its image in $\Big((\prod _{n\in \Nset}\kappa _n )\times \Rset ^{\omega _0}, \rho+t \Big)$.  However, this requirement is too strong. Indeed, recall that the metric space  $(\prod _{n\in \Nset}\kappa _n )\times \Rset ^{\omega _0}$ has a star-finite base for its uniformity as $(\prod _{n\in \Nset}\kappa _n ,\rho)$ is uniformly $0$-dimensional, that is, it has a base of partitions for its uniformity, and the uniformity on $\Rset ^{\omega _0}$ induced by $t$ is exactly the weak uniformity $wU_{t}(\Rset ^{\omega _0})$, which is of course star-finite (see Theorem \ref{real.weak}). Therefore, $(\varphi(X), \rho+t|_{\varphi(X)})$ has a star-finite base for its uniformity. By \cite[Corollary 3]{garrido.Banach-Stone}, if the Samuel compactification of $(X,d)$ is homeomorphic to the Samuel compactification of $(\varphi(X), \rho+t|_{\varphi(X)})$, then $(X,d)$ is uniformly homeomorphic to $(\varphi(X), \rho+t|_{\varphi(X)})$, and in particular it has a star-finite base for its uniformity too.  But this is not true for every Bourbaki-complete metric space as we have shown in Example \ref{point-finite}.
\end{remark}

\bigskip

\bigskip

\begin{center} \Denarius \Denarius  \Denarius \end{center}
\bigskip

\bigskip

\section{Metrization results}

\subsection{Bourbaki-completely metrizable spaces and related properties}

\hspace{15pt} Next, we are going to characterize those metric spaces that are topologically metrizable by a Bourbaki-complete metric.  This problem was solved and deeply studied in the work by Junnila, Hohti and Mero\~{n}o \cite{merono.hh} trough a property called {\it strong \v{C}ech-completeness.} Here we have decided to present a shorter proof than the one in \cite{merono.hh} that strongly depends of the embeddings from the previous section. 

Our characterization of the Bourbaki-completely metrizable spaces is related to the following properties.

\begin{definition} A space $X$ is {\it \v{C}ech-complete} if there exist open sets $G_n$, $n\in \Nset$ of $\beta X$ containing $X$ such that $X=\bigcap _{n\in \Nset}G_n$.
\end{definition}

It is well-known that {\it a space $X$ is \v{C}ech-complete if and only if it has a complete sequence of covers} (see \cite[Theorem 3.9.2 and notes on page 199]{engelkingbook}). Moreover, {\it in the frame of metrizable spaces,  \v{C}ech-complete spaces are the completely metrizable spaces, that is, those  spaces that are metrizable by a complete metric}.

Now, consider the following topological property.

\begin{definition} A space $X$ is {\it completely paracompact} if every open cover $\mathcal{G}$  has a refinement $\mathcal{V}$ which is a subcollection of a family of sets $\bigcup _{n\in \Nset} \mathcal{V}_n$ where each $\mathcal{V}_n$ is an open star-finite cover of $X$.
\end{definition}

It is clear that {\it every strongly paracompact space is completely paracompact}. Moreover, by Remark \ref{star-countable-cozero} we can  write star-countable instead of star-finite in the above definition.

Next, {\it every completely paracompact space is paracompact}. This follows from the fact that every star-countable cover is $\sigma$-{\it discrete}, that is, a countable union of discrete families of sets, and from the result of Michael \cite[Theorem 1]{michael} that states that if a space satisfies that every open cover has a $\sigma$-star-finite refinement then is paracompact. In particular, it is useful to recall the result of Morita \cite{morita.star-finite} that every {\it locally compact paracompact space is strongly paracompact and therefore also completely paracompact}.

\smallskip

Now, we need to consider complete paracompactness in the frame of metrizable spaces.

\begin{definition} A  space is  {\it strongly metrizable} if it has a base for the topology which consists of the union of countably many star-finite open covers.
\end{definition}

\noindent By the following result of Zarelua, completely paracompact metrizable spaces are exactly the {\it strongly metrizable spaces}.

\begin{theorem}\label{zarelua} {\rm (\cite[Lemma 5]{zarelua})} Let $X$ be a metrizable space. Then $X$ is strongly metrizable if and only if it is completely paracompact.
\end{theorem}

Again, by Remark \ref{star-countable-cozero} we could also write, in the above definition of strongly metrizable space, star-countable instead of star-finite.  Moreover, since star-finite covers are $\sigma$-discrete it follows from the Nagata-Smirnov-Bing metrization Theorem \cite{engelkingbook} that every strongly metrizable space is metrizable. 

An example of metrizable space (paracompact space) which is not strongly metrizable (completely paracompact) is provided by the following result of Wiscamb.

\begin{theorem} \label{wiscamb} {\rm (\cite[Theorem 5.2]{wiscamb})} A connected space is completely paracompact if and only if it is Lindel\"{o}f.
\end{theorem}

Observe that the family of, non necessarily metrizable, \v{C}ech-complete and completely paracompact spaces are well studied in \cite{merono.hh}. On the other hand, in the last subsection, we will analyze a uniform version of complete paracompactness in the same line than uniform strong paracompactness and uniform paracompactness.
\smallskip

From Theorem \ref{Bourbaki-complete.strongly-metrizable} it is immediate that {\it every Bourbaki-complete metric space is completely metrizable and strongly metrizable}. Now we prove the converse. 

\begin{theorem} \label{Bourbaki.completely-metrizable} Let $X$ be a completely metrizable and strongly metrizable space. Then $X$ is metrizable by a complete metric $D$ such that $\mu_D=s_f\mu _D$. In particular $(X,D)$ is Bourbaki-complete.
\begin{proof}  Since $X$ is strongly metrizable there exists a countable family $\{\mathcal{V}_n:n\in \Nset \}$ of star-finite open covers of $X$ such that $\bigcup _{n\in \Nset}\mathcal{V}_n$ is a base for the topology of $X$. In particular, by paracompactness of $X$, the family of all the star-finite open covers of $X$ form a base for the uniformity $s_f {\tt u}$.
Therefore we can apply the axioms of uniformity to the countable family of star-finite open covers $\{\mathcal{V}_n:n\in \Nset \}$. Thus, there exists a normal sequence of star-finite open covers $\langle \mathcal{U}_n \rangle _{n\in \Nset}$, such that $\mathcal{U}_n <\mathcal{V}_n$ for every $n\in \Nset$ (as in the proof of Theorem \ref{Bourbaki-complete.strongly-metrizable}).  Since  $\bigcup _{n\in \Nset}\mathcal{V}_n$ is a base for the topology of $X$, then $\bigcup _{n\in \Nset}\mathcal{U}_n$ is also a base for the topology of $X$.  

Indeed, let $G$ an open set of $X$ and $x\in X$. Since $\bigcup _{n\in \Nset}\mathcal{V}_n$ is a base, there exists some $V\in \mathcal{V}_n$ for some $n\in \Nset$, such that $x\in V\subset G$. Next, let $\rho$ be any metric on $X$, and choose $k\in \Nset$ such that $B_\rho (x, 1/k)\subset V$. Again by the base condition of $\bigcup _{n\in \Nset}\mathcal{V}_n$, there exists some $m\in \Nset$, such that for some $V'\in \mathcal{V}_m$
 $$x\in V'\subset B_\rho (x, 1/2k)\subset B_{\rho}(x,1/k)\subset V\subset G.$$ Next, consider the cover $\mathcal{U}_m$ and choose some $U\in \mathcal{U}_m$, such that $x\in U$. Since $\mathcal{U}_m<\mathcal{V}_m$ then $$x\in U\subset St^2(V',\mathcal{V}_m)\subset B_\rho(x,1/k)\subset V\subset G.$$ We have proved in this way that $\bigcup _{n\in \Nset}\mathcal{U}_n$ is also a base for the topology of $X$.

Next, apply Lemma \ref{pseudometric} to $\langle \mathcal{U}_n\rangle _{n\in \Nset}$ and let  $d$ be the pseudometric obtained. Then $d$ is compatible with the topology $X$, and in particular, $d$ is a metric. Moreover, the uniformity induced by $d$ has a star-finite base, that is, $\mu _d=s_f \mu_d$.

Consider now the completion $(\widetilde{X}, \widetilde{d})$ of $(X,d)$. Then $(\widetilde{X}, \widetilde{d})$  is complete and has a star-finite base by \cite[Lemma p. 370]{reynolds}, that is, $\mu _{\widetilde{d}}=s_f \mu_{\widetilde{d}}$. Now, since $X$ is completely metrizable, by \cite[Theorem 4.3.24]{engelkingbook}, $X$ is a $G_{\delta}$-set of $\widetilde{X}$. Thus, by \cite[Theorem 4.3.22]{engelkingbook}, $X$ is homeomorphic to a closed subspace of $(\widetilde{X}\times \Rset ^{\omega _0}, \widetilde{d}+ t)$.  Hence $X$ is metrizable by a complete metric $D$ satisfying that $\mu _D =s_f\mu _D$, precisely, the restriction of $\widetilde{d}+t$ over $X$. Indeed, it is not difficult to see that the product of two uniform spaces satisfying that both uniformities have a star-finite base, has also a star-finite base for the product uniformity. It is also clear the preservation of this property to subspaces. Finally, observe that, in particular, $(X,D)$ must be Bourbaki-complete.

\end{proof}

\end{theorem}

The next result collect several characterizations of the spaces metrizables by a Bourbaki-complete metric.

\begin{definition} The {\it cellularity} of a space $X$ is the supremum of the cardinal of all the partitions by open sets of the space.
\end{definition}

\begin{theorem} \label{summary} Let $X$ be a space. The following statements are equivalent:
\begin{enumerate}
\item $X$ is metrizable by a Bourbaki-complete metric;

\item $X$ is metrizable by a complete metric $d$, such that $\mu _d=s_f\mu _d$;

\item $X$ is homeomorphic to a closed subspace of $ D^{\omega_0} \times\Rset^{\omega _0}$ where $D$ is a discrete space of cardinal the cellularity of $X$;

\item $X$ is homeomorphic to a closed subspace of a countable product of locally compact metric spaces;

\item $X$ is \v{C}ech-complete and strongly metrizable;

\item $X$ is metrizable by a complete metric $d$ such that $\mu _d=s_c \mu _d$
\end{enumerate}

\begin{proof} $1)\Rightarrow 2)$ This is Theorem \ref{metrizable.star-finite}.

$2)\Rightarrow 3)$ Let $D$ a set a set of cardinal the cellularity of $X$. Then, it is clear that $|D|\geq sup \{\kappa _n:n\in \Nset\}$ where the $\kappa _n$'s are the cardinals from \ref{embedding.star-finite1} such that $X$ can be embedded as a closed subspace of $\prod_{n\in \Nset}\kappa  _n \times \Rset ^{\omega _0}$. Then for every $n\in \Nset$, $\kappa _n$ can be identified with a subspace of $D$ and thus $\prod_{n\in \Nset}\kappa _n$ is a closed subspace of the product space $D^{\omega _0}$. Then, the result follows from Theorem \ref{embedding.star-finite1}.

$3)\Rightarrow 4)$ This is trivial since $D$ and $\Rset$ are locally compact.

$4)\Rightarrow 1)$ First we show that every locally compact metric space $(Y,d)$ is metrizable by a Bourbaki-complete metric. For every $y\in Y$ let $V^y$ an open neighborhood of $y$ such that ${\rm cl} _Y V^y$ is compact. Put $\mathcal{V}=\{V^y:y\in Y\}$. By \cite[Lemma 38.1]{willard}, there exists a continuous pseudometric $\rho$ on $Y$ such that $\{B_{\rho}(y, 1):y\in Y\}<\mathcal{V}$. In particular ${\rm cl}_{Y} B_{\rho}(y, 1)$ is compact for every $y\in Y$. Now, consider the metric $d+\rho$ on $Y$. This metric is compatible with the topology of $Y$ and, it is clear that every set ${\rm cl}_{Y} B_{d+\rho}(y, 1)$, $y\in Y$ is compact, that is, $(Y, d+\rho)$is uniformly locally compact.

Next, since every uniformly locally compact space is Bourbaki-complete (Theorem \ref{unif.local.compact.th}), a countable product of Bourbaki-complete  metric spaces is a Bourbaki-complete metric space (Theorem \ref{product}), and Bourbaki-comple- teness is inherited by closed subspaces, the result follows.

$1)\Rightarrow 5)$ This is Theorem \ref{Bourbaki-complete.strongly-metrizable}.

$5)\Rightarrow 1)$ This is Theorem \ref{Bourbaki.completely-metrizable}.

$2)\Rightarrow 6)$ This is immediate.

$6\Rightarrow 5)$ We just need to prove complete paracompactness (see Theorem \ref{zarelua}). If $\mu_d=s_c\mu_d$ then for every $n\in \Nset$ we can take a uniform star-countable cover $\mathcal{V}_n$ such that $\mathcal{V}_n$ refines the uniform cover $\{B_d(x, \frac{1}{n}):x\in X\}$. Then $\bigcup _{n\in \Nset} \mathcal{V}_n$ is a base of for the topology of $X$ and by Remark \ref{star-countable-cozero}, it follows that $X$ is completely paracompact.

\end{proof}

\end{theorem}

From the result of Wiscamb, Theorem \ref{wiscamb}, the following corollary follows.

\begin{corollary} A connected  space $X$ is metrizable by a Bourbaki-complete metric if and only it is Lindel\"of and completely metrizable.
\end{corollary}

\noindent Therefore, any non separable, connected and complete metric space is an example of completely metrizable space which is not metrizable by a Bourbaki-complete metric.

\smallskip

In the next result we give several characterizations of the property of strong metrizability that can be deduced from the above result on Bourbaki-complete metrization. 

\begin{theorem} \label{summary2} Let $X$ be a space. The following statements are equivalent:
\begin{enumerate}
\item $X$ is metrizable by a metric $d$ such that ${\bf TB}_{d}(X)={\bf BB}_{d} (X)$;

\item $X$ is metrizable by a metric $d$ such that $\mu _d=s_f \mu _d$ ;

\item $X$ is homeomorphic to a subspace of $ D^{\omega_0} \times\Rset^{\omega _0}$ where $D$ is a discrete space of cardinality the cellularity of $X$;

\item $X$ is homeomorphic to a subspace of a countable product locally compact metric spaces;

\item $X$ is strongly metrizable;

\item $X$ is metrizable by a metric $d$ such that $\mu _d=s_c \mu _d$.
\end{enumerate}
\begin{proof} $1)\Rightarrow 2)$ If $(X,d)$ satisfies that ${\bf TB}_{d}(X)={\bf BB}_{d} (X)$, then its completion $(\widetilde{X}, \widetilde{d})$ is Bourbaki-complete by Corollary \ref{completion.BC}. By Theorem \ref{metrizable.star-finite}, $\widetilde{X}$ is metrizable by a complete metric $\rho$ such that $(\widetilde{X},\rho)$ satisfies that $\mu _\rho=s_f \mu _\rho$. Since this last property is clearly hereditary, the restriction of $\rho$ over $X$ is the desired metric.

$2)\Rightarrow 3)$ If $(X,d)$ satisfies that $\mu _d=s_f \mu _d$ then its completion $(\widetilde{X}, \widetilde{d})$ too by \cite[Lemma p. 370]{reynolds}. Then $\widetilde{X}$ is homeomorphic to a closed subspace of $ D^{\omega_0} \times\Rset^{\omega _0}$ by Theorem \ref{embedding.star-finite1}. Hence, $X$ is homeomorphic to a subspace of $ D^{\omega_0} \times\Rset^{\omega _0}$.

$3)\Rightarrow 4)$ This is trivial.

$4)\Rightarrow 1)$ Consider the closure of $X$ in the countable product of locally compact spaces and apply Theorem \ref{summary}. The result follows by Corollary \ref{completion.BC}.

$1) \Rightarrow 5)$ This implication follows from Theorem \ref{Bourbaki-complete.strongly-metrizable}, Corollary \ref{completion.BC} and from the fact that strong metrizability is an hereditary property.

$5) \Rightarrow 1)$ Let $X$ be a strongly metrizable space. Then for any compatible metric $d$ on $X$, if $(\widetilde{X}, \widetilde{d})$ denoted its completion, then $\widetilde{X}$ is strongly metrizable. In fact the open covers of $X$ are extended to $(\widetilde{X}, \widetilde{d})$. Thus $(\widetilde{X}, \widetilde{d})$ is strongly metrizable and complete. The result is then immediate from Theorem \ref{Bourbaki.completely-metrizable} and Corollary \ref{completion.BC}.

$2)\Rightarrow 6)$ This is immediate.

$6)\Rightarrow 5)$ This is like the proof of the implication $6)\Rightarrow 5)$ in Theorem \ref{summary}.
\end{proof}
\end{theorem}

\begin{remark} That a space $X$ is strongly metrizable if and only if it is homeomorphic to a subspace of $D^{\omega_0}\times \Rset^{\omega _0}$ is a well-known result (see \cite{nagata}, \cite{pears} an \cite{balogh}).
\end{remark}

\begin{example} {\it There is a complete metric space $(X,d)$ such that $\mu _d=s_f\mu _d$ but which is not strongly paracompact. In particular it is an example of strongly metrizable space which is not strongly paracompact}  

\begin{proof} Consider the Example \ref{nagata}, a take the product $(D^{\omega _0} \times\Rset,\rho +d _u)$, where $|D|=\omega _1$. It is a product of complete metric spaces so it is complete. In addition, it is easy to see that the uniformity induced by $\rho+d_u$ has a star-finite base. In particular it is Bourbaki-complete and then strongly metrizable (Theorem \ref{summary}). However it is not strongly paracompact since it is homeomorphic to $ D^{\omega _0} \times (0,1)$ which was shown not being strongly paracompact (\cite[Remark p. 169]{nagata}). \end{proof}
\end{example}

Finally, we characterize those spaces which are uniformizable by a Bourbaki-complete uniformity. These are exactly the $\delta$-complete spaces of Garc\'ia-M\'aynez (\cite{garcia.delta-complete}, see Definition 1.2.33).

\begin{theorem}\label{summary.uniform} For a space $X$ the following statements are equivalent:
\begin{enumerate}
\item $X$ is uniformizable by a Bourbaki-complete uniformity;

\item $X$ is $\delta$-complete;

\item $X$ is homeomorphic to a closed subspace of $ D^{\alpha} \times \Rset ^{\alpha}$ where $D$ is a discrete space of cardinality the cellularity of $X$;

\item $X$ is homeomorphic to a closed subspace of a product of locally compact metric spaces.
\end{enumerate} 
\begin{proof} $(1)\Rightarrow (2)$. This implications follows from the fact that every Bourbaki-uniform space is $\delta$-complete.

$(2)\Rightarrow (3)$ Since $(X, s_f \tt{u})$ is complete then, by Theorem \ref{embedding.star-finite2}, $X$ can be embedded as a closed subspace of $ \prod_{i\in I, n\in \Nset} \kappa _n ^i\times \Rset^{\alpha}$ where each $\kappa _n ^i$  are the cardinals endowed with the discrete topology. Moreover, the cardinal of $D$ is the cellularity of $X$ and then $|D|\geq sup\{\kappa _n ^i: i\in I, n\in \Nset\}$.  Thus, each $\kappa _n $ can be identified with a subset of $D$ and  $\prod_{i\in I, n\in \Nset} \kappa _n ^i$ is a closed subspace of the product space $D^{\alpha}$.

$(3)\Rightarrow (4)$ This is trivial.

$(4)\Rightarrow (1)$ Every locally compact metrizable space is metrizable by a uniformly locally compact metric by as in the proof of Theorem \ref{summary}. In particular, every uniformly locally compact metric space is Bourbaki-complete by Theorem \ref{unif.local.compact.th}.  Since Bourbaki-completeness is a productive property (Theorem \ref{product})  and hereditary by closed subspaces the result follows.
\end{proof}
\end{theorem}

\begin{remark} The equivalence of $2)$, $3)$ and $4)$ in the previous theorem was well-known by Garc\'ia-M\'aynez (see \cite{garcia.delta-complete} and \cite{garcia.equality}).
\end{remark}

\begin{example} {\it There exists a $\delta$-complete space which is not completely paracompact.} 
\begin{proof} Since $e{\tt u}\leq s_f{\tt u}$ (Lemma \ref{morita}), any realcompact non paracompact space is such an example. For instance, the product of two Sorgenfrey lines.
\end{proof}
\end{example}

\medskip

\subsection{Metric spaces that are Bourbaki-complete and cofinally complete at the same time}

{\hspace{15pt}} In this section we study those metric spaces that are metrizable by a metric which is Bourbaki-complete and cofinally complete at the same time (Theorem \ref{sigma.unif3}). The interest in these spaces lies in the fact that these are close to the metric spaces metrizable by a cofinally Bourbaki-complete metric (Theorem \ref{cof.B.complete.metrizable}) but they are topologically weaker (Example \ref{ex.strongly.metrizable}). 

On the other hand, recall that complete paracompactness is a property lying between strong paracompactness and paracompactness. Therefore it is interesting to ask if there exists a uniform extension of completeness lying between the properties of uniform strong paracompactness and uniform paracompactness, in such a way that it coincides with complete paracompactness when the space is endowed with the fine uniformity ${\tt u}$. To that purpose we we are going to introduce the following family of covers.

\begin{definition} An open cover $\mathcal {U}$ of a space $X$ is {\it $\sigma$-star-finite} ({\it $\sigma$-star-countable}) if there exists a countable family $\mathcal{U}_n$, $n\in \Nset$, of star-finite (star-countable) open covers of $X$ such that  $\mathcal{U}\subset\bigcup_{n\in \Nset} \mathcal{U}_n$. 

\end{definition}

We will denote by $\sigma$-$s_f\mu$ the family of all the  uniform covers  from a uniform space $(X,\mu)$ having a $\sigma$-star-finite uniform open refinement. We cannot tell if this family of uniform covers is in general a base for some compatible uniformity on $X$. However, we have some particular examples of uniform spaces such that the original uniformity has a base of $\sigma$-star-finite open covers.

\begin{theorem} \label{sigma.unif2} The following statements are equivalent for a space $X$
\begin{enumerate}

\item $X$ is completely paracompact; 

\item $X$ is paracompact and ${\tt u}=\sigma$-$s_f{\tt u}$;

\end{enumerate}

\begin{proof} The proof is clear from the definition of completely paracompact space. Moreover, recall that if a space $X$ is paracompact, a base of the uniformity ${\tt u}$ is given by all the open covers of $X$. 
\end{proof}

\end{theorem}

Another example of uniform space having a $\sigma$-star-finite open base for its uniformity is any uniform space $(X,\mu)$ such that $\mu=s_f\mu$. Indeed, any open star-finite cover is clearly $\sigma$-star-finite. Moreover, we next prove that, if $(X,\mu)$ has a countable base for its uniformity, then it has also a $\sigma$-star-finite open base for its uniformity.

\begin{theorem} \label{sigma.unif1} Let $(X,\mu)$ be a uniform space such that $\mu=e\mu$. Then $\mu$ has a a base of $\sigma$-star-finite open covers.
\begin{proof} If  $\mathcal{U}=\{U_n:n\in \Nset\}$ is a countable uniform cover then each cover $\mathcal{V}_n=\{U_n, \bigcup_{j\neq n} U_j\}$ is uniform and star-finite and $\mathcal{U}\subset \bigcup_{n\in \Nset} \mathcal{V}_n$, that is, $\mathcal{U}$ is $\sigma$-star-finite.
\end{proof}
\end{theorem}

\begin{remark} It is not difficult to see that for a connected uniform space $(X,\mu)$, the uniformity $\mu$ has a countable base, that is $\mu =e\mu$, if and only if $\mu$ has a base of $\sigma$-star-finite open covers. The proof is like Lemma \ref{cozero.connected}. Moreover, by this fact, we can easily deduce the result by Wiscamb Theorem \ref{wiscamb}.
\end{remark}

\smallskip

Recall that in Theorem \ref{uniform.strongly.paracompact}  it is proved that a uniform space $(X,\mu)$ is uniformly strongly paracompact if and only if $(X,\mu)$ is cofinally complete and $\mu$ has a star-finite base. Moreover, by Theorem \ref{point.finite}, every cofinally complete uniform space has a point-finite base for its uniformity.  In parallel to these results we propose the following definition of uniformly completely paracompact space.


\begin{definition} A uniform space $(X,\mu)$ is {\it uniformly completely paracompact} if it is cofinally complete (equivalently, uniformly paracompact) and the uniformity $\mu$ has a base of $\sigma$-star-finite open covers.
\end{definition}

From Theorem \ref{sigma.unif2} it is clear that a topological space $X$ is completely paracompact if and only if, when it is endowed with the fine uniformity ${\tt u}$,  the space $(X, {\tt u})$ is uniformly completely paracompact. 


\smallskip

Next, we are going to characterize those metric spaces that are metrizable by a uniformly completely paracompact metric. To that purpose recall the following result by Romaguera (a proof can be found also in \cite[Theorem 4.1]{beer-between1}).

\begin{theorem} \label{romaguera} {\rm (\cite{romaguera})} A metrizable space is metrizable by a cofinally complete metric if and only if the family $nlc(X)$ of points of $X$ without a locally compact neighborhood,  is compact.

\end{theorem}



\begin{theorem} \label{sigma.unif3}Let $X$ be a metrizable space. The following statements are equivalent:

\begin{enumerate}

\item $X$ is metrizable by a uniformly completely paracompact metric;

\item $X$ is strongly metrizable and $nlc(X)$ is compact;

\item $X$ is metrizable by a metric which is Bourbaki-complete and cofinally complete at the same time.
\end{enumerate}

\begin{proof}
$(1)\Rightarrow (2)$ This implication follows from Theorem \ref{romaguera}, Corollary \ref{sigma.unif2} and Theorem \ref{zarelua}.

$(2)\Rightarrow (1)$ Suppose first that $nlc(X)=\emptyset$. Then $X$ is locally compact and, as in the proof $(4)\Rightarrow (1)$ of Theorem \ref{summary}, $X$ is metrizable by a uniformly locally compact metric $d$. Therefore, by Theorem \ref{unif.local.compact.th}   $(X,d)$ is cofinally Bourbaki-complete. in particular it is uniformly completely paracompact, by Theorem \ref{uniform.strongly.paracompact}.

Otherwise, assume that $nlc(X)\neq \emptyset$ and let $\rho$ a metric on $X$. Since $nlc(X)$ is compact, there exists a countable family of open sets $\{W_1,...,W_k,...\}$ in $X$ such that for every open subset $A$ of $X$ containing $nlc(X)$ there exists $k\in \Nset$ satisfying that $nlc(X)\subset W_k\subset A$. For instance, consider $W_k=nlc(X)^{1/k}$, $k\in \Nset$. Now for every $x\notin nlc(X)$ take $V^x$ an open neighborhood of $x\in X$ with compact closure. For every $k\in \Nset$, let $\mathcal{G}_k=\{V^x :x\notin W_k\}\cup \{W_k\}$. 

Now, we start by $\kappa=1$. By strong metrizability (equivalently complete paracompactness), the cover $\mathcal{A}_1=\mathcal{G}_1$ has a $\sigma$-star-finite open refinement $\mathcal{U}_1$. Next, consider the open cover $\mathcal{U}_1\wedge \mathcal{G}_2$, and take an open cover $\mathcal{A}_2$ such that $\mathcal{A}_2^* <\mathcal{U}_1 \wedge \mathcal{G}_2$.  Again, by complete paracompactness we can take an open refinement $\mathcal{U}_2$ of $\mathcal{A}_2$ being $\sigma$-star-finite. 

Thus, proceeding in this way we obtain a normal sequence $\langle{\mathcal{U}_n} \rangle$ of $\sigma$-star-finite open covers such that $\mathcal{U}_k<\mathcal{G}_k$ for every $k\in \Nset$. Now, applying Lemma \ref{pseudometric} there exists compatible metric $d$ on $X$ such that $$\mathcal{B}_{1/2^{k+1}}<\mathcal{U}_k<\mathcal{B}_{1/2^{k-1}} \text{ for every }k\in \Nset$$ (where $\mathcal{B}_{\varepsilon}=\{B_{d}(x,\varepsilon):x\in X\}$). Therefore the metric uniformity $\mu_d$ has a base of $\sigma$-star-finite open covers. 

Next, we prove that $(X,d)$ is cofinally complete. Let $(x_n)$ be a cofinally Cauchy sequence in $(X,d)$. If for some $\kappa \in \Nset$, $(x_n)$ is cofinally in some $V^x\in \mathcal{G}_k$, where $V^x$ is one of the above sets with compact closure, then $(x_n)$ clusters in ${\rm cl}_{X} V^x$ by compactness. Otherwise, it follows that for every $k\in \Nset$, $(x_n)$ is eventually in $W_k$. Indeed, by cofinal-Cauchyness, as $\mathcal{B}_{1/2^{k+1}}<\mathcal{U}_k<\mathcal{G}_k$, then for every $k\in \Nset$ there must be some $G_k\in \mathcal{G}_k$ such that $(x_n)$ is cofinally in some $G\in \mathcal{G}_k$. 

Suppose that $(x_n)$ does not cluster.  Thus, $\mathcal{A}=\{X\backslash {\rm cl}_{X} \{x_n:n\geq k\}:k\in \Nset\}$ is an open cover of $X$ and in particular it is and open cover of $nlc(X)$. Since $nlc(X)$ is compact, there exists a finite subfamily $\mathcal{A}'\subset \mathcal{A}$ such that $nlc(X)\subset \bigcup \mathcal{A}'$. Hence, for some $k\in \Nset$ we have that $nlc(X)\subset W_k\subset \bigcup \mathcal{A}'$, and, since $\mathcal{A}'$ is finite, this implies that for some $A\in \mathcal{A }'$, $A\cap \{x_n:n\geq k\}\neq \emptyset$ for every $k\in \Nset$. But this is a contradiction.

$(2)\Rightarrow (3)$. Let $X$ be a strongly metrizable such that $nlc(X)$ is compact. In particular $X$ is completely metrizable. Therefore, by Theorem \ref{romaguera} and Theorem \ref{summary}, $X$ is metrizable by a metric $\rho$ and a metric $t$, which are cofinally complete and Bourbaki-complete respectively. Let $d: X\times X \rightarrow [0,\infty)$,  $d(x,y)=\max\{\rho (x,y), t(x,y)\}$. Then it is easy to check that $d$ is a metric compatible with the topology of $X$ which is cofinally complete and Bourbaki-complete, as for every $x\in X$ and every $\varepsilon >0$, $B_{d}(x,\varepsilon)=B_{\rho}(x,\varepsilon)\cap B_{t}(x,\varepsilon)$.

$(3)\Rightarrow (2)$ This follows from Theorem \ref{Bourbaki.completely-metrizable} and Theorem \ref{romaguera}.
\end{proof} 
\end{theorem}

Similarly to the above theorem of metrization by a uniformly completely paracompact metric we have the next theorem for spaces which are metrizable by a cofinally Bourbaki-complete metric.

\begin{theorem}\label{cof.B.complete.metrizable} {\rm (\cite[Theorem 33]{merono.completeness})} Let $X$ be a metrizable space. The following statements are equivalent:
\begin{enumerate}
\item $X$ is metrizable by a cofinally Bourbaki-complete (uniformly strongly paracompact) metric;

\item $X$ is strongly paracompact and $nlc (X)$ is compact.
\end{enumerate}
\begin{proof} 
$1)\Rightarrow 2)$ This implication follows from Theorem \ref{uniformly.strong.fine}.

$2)\Rightarrow 1)$ This proof is similar to the proof of the implication $(2)\Rightarrow (1)$ in the above Theorem \ref{sigma.unif3}. Just change complete paracompactness by strong paracompactness and $\sigma$-star-finite covers by star-finite covers.

\end{proof}
\end{theorem}

\begin{example}  \label{ex.strongly.metrizable} {\it There exists a uniformly completely paracompact metric space $(X,d)$ which is not Bourbaki-complete and not strongly paracompact. In particular, $(X,d)$ is not metrizable by a cofinally Bourbaki-complete metric, even if, by Theorem \ref{sigma.unif3} it is metrizable by a metric which is Bourbaki-complete and cofinally complete at the same time.}

\begin{proof} [Construction]
This is a subspace of the metric hedgehog (Example \ref{hedgehog}) $(H(\omega _1), \rho)$. Let $\{A_n:n\in \Nset\}$ be a partition of $\omega _1$ such that $|A_n|=\omega _1$ for every $n\in \Nset$. Let $L=\{0\}\cup \bigcup_{n\in \Nset} E_n$ where $E_n=\{[(x, \alpha)]: \frac{1}{n}\leq x\leq 1, \alpha \in A_n\}$. Then $(L,\rho)$ is a cofinally complete metric space since it is a closed subspace of $(H(\omega _1),\rho)$. However, it is not Bourbaki-complete because the Bourbaki-bounded subset given by taking just one point of the form $[(1,\alpha)]$ in each $E_n$ is a closed but not compact.

On the other hand $L$ is strongly metrizable. Indeed, fixed $k\in \Nset$, the open cover $\mathcal{B}_{1/k}$ has the following $\sigma$-star-finite refinement. For every $n\in \Nset$ and every $\alpha _n\in A_n$, in the spine $S_{n,\alpha}=\{[(x,\alpha)]:\frac{1}{n}\leq x\leq \alpha\}$ we can choose, by compactness, a finite cover $\mathcal{G}_{n,\alpha}$,  of open balls of radius $\varepsilon _n =min\{1/2k, 1/2n\}$ and center in $ S_{n,\alpha}$. For every $n\in \Nset$ define the open covers $$\mathcal{A}_n=\{B_{\rho}(0,\varepsilon _n)\}\cup \big(\bigcup_{j\leq n}\bigcup_{ \alpha \in A_j}\mathcal{G}_{j,\alpha}\big)\cup\Big{\{}\bigcup_{j>n} E_j\Big{\}}.$$ Then, it is clear that 
each $\mathcal{A}_n$ is open and star-finite. Moreover, $\bigcup_{n\in \Nset}\mathcal{A}_n$ contains a refinement $\mathcal{A}$ of $\mathcal{B}_{1/k}$, $$\mathcal{A}=\{B_{\rho}(0,1/2k)\}\cup \bigcup_{j\in \Nset}\bigcup_{ \alpha \in A_j}\mathcal{G}_{j,\alpha}.$$

However, $X$ is not strongly paracompact. For instance, take again the open cover $\mathcal{B}_{1/k}$. Then for every $0<\varepsilon \leq 1/k$, the open ball $B_{\rho}(0,\varepsilon)$ meets always uncountably many parwise disjoint open balls $B_{\rho}(x, \delta)$,  of center $x\notin B_{\rho}(0,\varepsilon)$, for any $0<\delta\leq \varepsilon$. Therefore it is easy to deduce that  $\mathcal{B}_{1/k}$ cannot have a star-finite (star-countable) open refinement.
\smallskip

By Theorem  \ref{sigma.unif3} and Theorem \ref{cof.B.complete.metrizable}, $L$ is metrizable by a metric which is Bourbaki-complete and cofinally complete at the same time but not cofinally Bourbaki-complete

\end{proof}
\end{example}

\smallskip

{\bf OPEN PROBLEMS:} In the above example, is it possible to construct explicitly a metric which is Bourbaki-complete and cofinally complete at the same time but which is not cofinally Bourbaki-completely metrizable? Is every metric space $(X,d)$, cofinally complete and Bourbaki-complete at the same time, also uniformly completely paracompact? 

\smallskip

\begin{example}{\it There exists a complete metric space $(X,d)$ having a base $\sigma$-star-finite open covers for the uniformity, but which is not cofinally complete, nor Bourbaki-complete. }

\begin{proof} By Theorem \ref{sigma.unif2}, any separable Banach space, as the Hilbert space $\ell _2$ is such example. Indeed by theorem \ref{banach.BC}, $\ell_2$ is not Bourbaki-complete because it is infinite dimensional. Moreover, by Theorem \ref{hohti.metric} it is not cofinally complete.
\end{proof}
\end{example}

\noindent Related to the above example, we have the following topological question. 
\smallskip

{\bf OPEN PROBLEM:} Is there an example of  completely paracompact space not being $\delta$-complete? 
\smallskip

We close this section with the following observations. Recall that in \cite{musaev}, a uniform extension of complete paracompactness and strong paracompactness, called R-completely paracompactness and R-strong paracompactness, respectively, were given. More precisely, a uniform space $(X,\mu)$ is {\it R-completely paracompact} ({\it R-strongly paracompact}) if for every open cover there exists a $\sigma$-star-finite open refinement (star-finite open refinement) which is in addition uniformly locally finite. 

Observe that the property of $R$-strong paracompactness is clearly uniformly weaker than Hohti's uniform strong paracompactness definition. Indeed, applying Theorem \ref{cof.B.complete.metrizable}, every separable  cofinally complete space metric space which is not Bourbaki-complete is a counterexample. 

We strongly believe that the same happens for the definition of $R$-completely paracompact and our definition of uniform complete paracompactness. However, we don't have a counterexample. By this reason we propose the following open problems.

\smallskip

{\bf OPEN PROBLEMS:} Give a definition of uniform complete paracompactness by means of covers. Give an example of $R$-completely paracompact space which is not uniformly completely parcompact.


\medskip

\bigskip

\bigskip

\begin{center} \Denarius \Denarius  \Denarius \end{center}
\bigskip

\bigskip

\newpage

\chapter{The Samuel realcompactification of a uniform space}
\thispagestyle{empty}

\newpage

\section{Primary results}

\subsection{Basic facts about realcompactifications}

{\hspace{15pt}} The definition of realcompact space that we use is the given in Definition 1.2.35, that is, a space $X$ is {\it realcompact} if it is homeomorphic to a closed subspace of a product of real lines. Moreover, similarly to compactification we get the notion of realcompactification.


\begin{definition} A {\it realcompactification} of a space $X$ is a realcompact space $Y$ in which $X$ is densely embedded.
\end{definition}


Next, we give some well-known results on realcompactifications. Many of them can be found in \cite{garrido1}.

A classical  way of generating realcompactifications of a space $X$ is the following. First, we take a family $\mathcal L\subset C(X)$ of real-valued continuous functions, \textit{separating points and closed sets of $X$}. Then, we embed homeomorphically $X$ into the product space of real lines $\Rset ^{\mathcal{L}}$, through the evaluation map $$e: X\rightarrow \Rset ^{\mathcal{L}}$$ $$\hspace{2.7cm} x\mapsto e(x)=(f(x))_{f\in \mathcal{L}}.$$ Then, the closure of $e(X)$ in $\Rset ^{\mathcal{L}}$ is known as the \textit{realcompactification of $X$ generated by $\mathcal{L}$}. Usually, $\mathcal{L}$ has some algebraic structure. Here, we will suppose that $\mathcal{L}$ is at least a \textit{unital vector lattice}, also because we will work with families of real-valued functions that are lattices but which are not necessarily an algebra (ring), as they are in the earliest paper by Isbell \cite{isbell.algebras}. We will denote by $H(\mathcal{L})$ the realcompactification generated by $\mathcal{L}$  because it is exactly the set of all the real unital vector lattice homomorphisms on $\mathcal{L}$. 
\smallskip

On the other hand, for a unital vector lattice  $\mathcal{L}\subset C(X)$ we can consider $w {\mathcal{L}}$ the weak uniformity in $X$ which is the weakest uniformity making each function in $\mathcal{L}$ uniformly continuous \cite{willard}. When $\mathcal{L}$ separates points and closed sets in $X$, then $w\mathcal{L}$ is a Hausdorff uniformity compatible with the topology of $X$. If we endow $X$ with the weak uniformity $w{\mathcal{L}}$ and $\Rset ^{\mathcal{L}}$ with the usual product uniformity $\pi$, then the evaluation map $e:X\rightarrow \Rset ^{\mathcal{L}}$ is now uniformly  continuous and the inverse map $$e^{-1}: e(X)\rightarrow X$$ is also uniformly continuous. Thus, $(X,w\mathcal{L})$ is uniformly embedded in $(\Rset ^{\mathcal{L}},\pi)$. Since $(\Rset ^{\mathcal{L}},\pi)$ is a complete uniform space, the closure  $H(\mathcal{L})$ of $X$ in $\Rset ^{\mathcal{L}}$ is the completion of $(X, w{\mathcal{L}})$ (by uniqueness of the completion).  We can summarize all of this  as follows.

\begin{theorem} \label{completion.realcompactification} Let  $\mathcal{L}\subset C(X)$ be  a unital vector lattice  separating points and closed sets in the space $X$. The realcompactification $H(\mathcal{L})$ of $X$  is (topologically) homeomorphic to the completion of the uniform space $(X, w{\mathcal{L}})$ where $w{\mathcal{L}}$ is the weak uniformity generated  by $\mathcal{L}$.
\end{theorem}

If we just take the bounded functions $\mathcal{L}^{*}=\mathcal{L}\cap C^{*}(X)$  in $\mathcal{L}$ (where $C^{*}(X)$ denotes the family of bounded real-valued continuous functions) we get that $H(\mathcal{L}^{*})$ is now a compactification of $X$.
\smallskip

Likewise compactifications, we can consider a partial order $\leq$ on the set $\mathfrak{R}(X)$ of all the realcompactifications of $X$ \cite{engelking.realcompact}. Namely, for two realcompactifications $\alpha _1 X$ and $\alpha _2 X$, we write $\alpha_1 X\leq \alpha_2 X$ whenever there is a continuous mapping $h:\alpha_2X\to \alpha_1X$ leaving $X$ pointwise fixed.  We say that $\alpha_1 X$ and $\alpha_2X$  are equivalent whenever $\alpha_1 X\leq \alpha_2 X$ and $\alpha_2 X\leq \alpha_1 X$, and this implies the existence of a homeomorphism between $\alpha_1 X$ and $\alpha_2 X$ leaving $X$ pointwise fixed. 

\begin{theorem} {\rm (\cite{garrido1})} The realcompactification $H(\mathcal L)$ (resp. $H(\mathcal L^*)$) is characterized  (up to equivalence) as the smallest realcompactification (resp. compactification) of $X$ such that every function $f\in \mathcal{L}$ can be continuously extended to it.\end{theorem}

For the next result, recall that, if $f\in \mathcal{L}$, we can regard $f$ as a continuous function from X into the one-point compactification $\Rset \cup \{\infty\}$ of $\Rset$. 

\begin{theorem}\label{realcompact.subspace} The realcompactification  $H(\mathcal L)$ can be considered as a topological subspace of $H(\mathcal L^*)$. Thus, we can write $$X\subset H(\mathcal{L})\subset H(\mathcal{L}^{*}).$$ In particular, every $f\in \mathcal{L}$ can be extended to a unique continuous function $f^{*}: H(\mathcal{L}^{*})\rightarrow \Rset \cup \{\infty\}$. Moreover, $$H(\mathcal{L})=\{\varphi \in H(\mathcal{L}^{*}): f^{*}(\varphi)\neq\infty \text{ for all }f\in \mathcal{L}\}.$$ 

\end{theorem}

The partial order set $(\mathfrak{R}(X), \leq)$ is a complete upper semi-lattice where the largest element is exactly the Hewitt-Nachbin realcompactification $H(C(X))=\upsilon X$. Hewitt  himself observed \cite{hewitt} that the realcompactification that bears his name can be obtained as the completion of the space $X$ endowed with the weak $wC(X)$ uniformity generated by the family of real-valued continuous functions. Moreover, Nachbin \cite{nachbin} developed his theory of realcompactness, 
considering those spaces which are complete in the above uniformity. Later, Shirota \cite{shirota1}, \cite{shirota2} proved that the Hewitt realcompactification is homeomorphic to the completion of the space $X$ together with the uniformity generated by all the countable normal covers $e{\tt u}$, and hence that realcompactness is equivalent to the completeness of $(X, e{\tt u})$ (Theorem \ref{shirota}).

\smallskip

Recall that, in the case of the compactifications $\mathfrak{K}(X)$ of the space $X$, the partial order set $(\mathfrak{K}(X), \leq)$ is also a complete upper semi-lattice, where the largest element is now the Stone-\v{C}ech compactification $H(C^*(X))=\beta X$. However, we have that, in the partial order set of the realcompactifications, $$\upsilon X \geq \beta X$$ and by Theorem \ref{realcompact.subspace} $$X\subset \upsilon X\subset \beta X.$$ Precisely, $\upsilon X$ is the smallest realcompact subspace of $\beta X$ containing $X$ \cite[Theorem 8.5.b]{gillman}.

The partial order set of realcompactifications, as well as the partial order set of compactifications, of a space $X$, is a complete lattice if and only if $X$ is locally compact. In this case, the smallest element in both lattices is the Alexandroff compactification, also called the one-point compactification of $X$ which is generated by the lattice of all the real-valued  functions which are constant at infinity \cite{mack}. 
\smallskip

Next, let $\alpha X$ be a realcompactification of $X$, which is not explicitly generated by a lattice of real-valued continuous functions. Then if $C(\alpha X)$ is the algebra of all the real-valued continuous functions on $\alpha X$ then the uniform space $(\alpha X, wC(\alpha X))$ is complete \cite{gillman}. Precisely, it is the completion of $(X, wC(\alpha X)|_{X})$ where $C(\alpha X)|_{X}$ consists of the restrictions to $X$ of the functions in $C(\alpha X)$.  


The following theorems on extensions of maps are well-known. 


\begin{theorem} \label{cc}{\rm (\cite{borsik})} Let $(X,\mu)$ and $(Y,\nu)$ be uniform space spaces and let $(\widetilde{X},\widetilde{\mu})$ and $(\widetilde{Y},\widetilde{\nu})$ denote their respective completions. Let $f: (X,\mu)\rightarrow (Y,\nu)$. Then $f$ can be extended to a (unique) continuous map $F:(\widetilde{X},\widetilde{\mu})\rightarrow (\widetilde{Y},\widetilde{\nu})$ if and only if $f$ is Cauchy-continuous.
\end{theorem}

It is clear that every uniformly continuous map  is Cauchy-continuous and that every Cauchy continuous map is continuous.

\begin{theorem} \label{dimaio} {\rm (\cite{dimaio})} If a uniform space is complete then every real-valued  continuous function is Cauchy-continuous. The converse is in general false. However it is true in the frame of metric spaces.
\end{theorem}

An example of non-complete uniform space such that every continuous function is Cauchy-continuous is given by any fine uniform space which is not complete as the space of all the countable ordinals $[0,\omega _1)$.

\smallskip

By all the foregoing, we can describe  $C(\alpha X)_{|X}$ as the family of those continuous real-valued functions $X$ which are Cauchy-continuous when  $X$ is endowed with the weak uniformity $wC(\alpha X)_{|X}$  (see  \cite{colebunders}). Moreover, $C(\alpha X)$ is isomorphic to the subalgebra  $C(\alpha X)_{|X}$ of $C(X)$ \cite{hager-johnson}. From now on, we will not distinguish  between $C(\alpha X)$ and the corresponding subalgebra $C(\alpha X)_{|X}$ of $C(X)$.

\begin{theorem}\label{determined}{\rm (\cite[Theorem 3]{engelking.realcompact})} Given two realcompactifications $\alpha _1 X$ and $\alpha _2 X$ of a space $X$, $\alpha _1 X\leq \alpha _2 X$ if and only if $C(\alpha _1 X)\subset C(\alpha _2 X)$ and $\alpha _1 X$ and $\alpha _2 X$ are equivalent realcompactifications if and only if $C(\alpha _1 X)=C(\alpha _2 X)$.
\end{theorem}

By the previous theorem the algebra $C(\alpha X)$ is uniquely determined by $\alpha X$ and vice versa. However, different algebras of functions can generate the same realcompactification (\cite{blasco},\cite{redlin}). 





Moreover, given $\mathcal{L}\subset C(X)$ and the realcompactification $H(\mathcal{L})$ generated by $\mathcal{L}$, it is clear that it is homeomorphic  to the realcompactification $H(C(H(\mathcal{L})))$ generated by $C(H(\mathcal{L}))$. However, they are not necessarily uniformly homeomorphic when they are endowed respectively with the weak uniformities, $w\mathcal{L}$ and $wC(H(\mathcal{L}))$ as we will point later.

\medskip

\subsection{The Samuel realcompactification}

{\hspace{15pt}} In the frame of uniform spaces $(X,\mu)$ we have uniformly continuous functions. The family of all the real-valued uniformly continuous functions $U_{\mu}(X)$ is of course a unital vector lattice which separates points from closed sets of $X$. Thus we get the realcompactification $H(U_{\mu}(X))$ which is exactly the smallest realcompactification of $X$ such that every uniformly continuous function can be continuously extended. 

\begin{definition} \textit{The Samuel realcompactification} of a uniform space $(X,\mu)$ is the space $H(U_{\mu}(X))$. Whenever $X=H(U_{\mu}(X))$ we will say that the space $X$ is \textit{Samuel realcompact.} 
\end{definition}


Recall that we have defined the Samuel compactification $s_\mu X$ of a uniform space $(X,\mu)$ as the completion of $(X,f\mu)$. Since $f\mu=wU^{*}_{\mu}(X)$, $s_\mu X$ is exactly $H(U^{*}_{\mu}(X))$, the smallest compactification (and realcompactification) of $(X,\mu)$ such that every bounded  uniformly continuous function can be continuously extended. By Theorem \ref{realcompact.subspace} we know that  $$X\subset H(U_{\mu}(X))\subset s_{\mu} X$$ and that $$H(U_{\mu}(X))=\{\varphi \in s_{\mu} X: f^{*}(\varphi)\neq\infty\text{ for all }f\in U_{\mu}(X)\}$$ where $f^{*}: s_{\mu} X \rightarrow \Rset \cup \{\infty\}$ is the unique continuous extension of $f$ to the above domain and range.

Next, as well as every compactification is the Samuel compactification for some suitable uniformity on the space \cite{gal}, every realcompactification is equivalent to the Samuel realcompactification for some uniformity on $X$. In fact, let $\alpha X$ be a realcompactification of X and observe that $U_{wC(\alpha X)}(X)=C(\alpha X)$, because every uniformly continuous functions is Cauchy-continuous \cite{borsik}. Thus, we have the following result.

\begin{theorem}\label{all.samuel} Let $\alpha X$ be a realcompactification of the space $X$. Then $\alpha X$ is (topologically) homeomorphic to the Samuel realcompactification of the uniform space $(X, wC(\alpha X))$.
\end{theorem}

For instance, the Hewitt-Nachbin realcompactification $\upsilon X$ can be considered the Samuel realcompactification of the uniform space $(X, w_{C(X)})$ \cite{gillman}. Equivalently, $\upsilon X$ is the Samuel realcompactification of the uniform space $(X,{\tt u})$ where ${\tt u}$ is the universal uniformity on $X$, because the family of real-valued uniformly continuous functions on $(X,{\tt u})$ is exactly $C(X)$. In the same way, the Stone-\v{C}ech compactification $\beta X$ is the Samuel compactification of the uniform space $(X, w_{C^{*}(X)})$.

By Theorem \ref{real.weak} the above result is immediate

\begin{theorem} \label{R-Samuel.realcompact} Let $\alpha$ be a cardinal and $(\Rset ^{\alpha},\pi)$ be the product of $\alpha$ real lines endowed with the product uniformity $\pi$ of the euclidean uniformities on each factor. Then, $(\Rset ^{\alpha},\pi)$ is Samuel realcompact.
\end{theorem}

Observe that in the above result the Hewitt realcompactification $\upsilon \Rset^{\alpha}$ and the Samuel realcompactification $H(U_{\pi}(\Rset^{\alpha}))$ are equivalent since both realcompactifications are homeomorphic to $\Rset^{\alpha}$. However, $wU_{\pi}(\Rset^{\alpha})\neq wC(\Rset^{\alpha})$  $(\Rset^{\alpha},\pi)$ is not a UC space.

In general, it is clear that every Samuel realcompact space $(X,\mu)$ satisfies that $\upsilon X$ and $H(U_{\mu}(X))$ are equivalent. However, there are uniform space spaces $(X,\mu)$ satisfying that both realcompactifications are equivalent  which are not necessarily Samuel realcompact.

\begin{definition} A uniform space $(X,\mu)$ is a {\it UC space} if every real-valued continuous function is uniformly continuous.
\end{definition}

Every UC uniform space satisfies trivially that $\upsilon X$ and $H(U_{\mu}(X))$ are equivalently realcompactifications. For instance, every uniformly discrete metric space is a UC space and there are uniformly discrete space which are not realcompact and hence not Samuel realcompact. Whether or not a uniformly discrete metric space is realcompact depends on some particular assumption about the cardinality of the space, as we will explain later. However, there are also examples of non-realcompact UC spaces which do not depend of any cardinality property. Indeed, observe that every fine uniform space is UC by Theorem \ref{map.fine}. So any pseudocompact non-compact fine uniform space, as $([0,\omega _1), {\tt u})$, does the work because a pseudocompact space is realcompact if and only if it is compact \cite{engelkingbook}. 

In addition, observe that every metric UC space is a fine space (\cite{dimaio}). However, this is not true in the frame of uniform spaces. For instance, take a uniformly discrete space $(D,d)$ of uncountable cardinality and endowed it with the weak uniformity $wC(D)$. This space is still UC but it is not fine since ${\tt u}$ does not have a countable base.

\begin{remark} It is well-known \cite{woods} that a metric space $(X,d)$ is $UC$ if and only if $s_d X$ and $\beta X$ are equivalent compactifications.  Then, every metric space with $s_d X=\beta X$ satisfies also that $H(U_{d}(X))$ and $\upsilon X$ are equivalent. However, it is clear that there are metric spaces $(X,d)$ satisfying that $H(U_{d}(X))$ and $\upsilon X$ are equivalent  but such that $s_\mu X$ and  $\beta X$ are not equivalent. Just consider the euclidean space $(\Rset, d)$, which is not UC . 
On the other hand, there are uniform spaces $(X,\mu)$ satisfying that $s_\mu X=\beta X$ but $H(U_{\mu}(X))\neq\upsilon X$. This is due to the fact that in the general frame of uniform spaces even if $C^{*}(X)=U^{*}_{\mu}(X)$ then $(X,\mu)$ is not necessarily $UC$. Indeed, just endowed the real-line with the weak uniformity $\mu=wC^{*}(\Rset)$. Then, $C^{*}(X)=U^{*}_{\mu}(X)=U_{\mu}(X)\subsetneq C(X)$.
\end{remark}

In spite of Proposition \ref{R-Samuel.realcompact} and of the UC spaces, in general $\upsilon X$ and $H(U_{\mu}(X))$ are not equivalent realcompactifications. Similarly, not every realcompact space is Samuel realcompact. In order to give such example, we characterize first the uniform spaces such that the Samuel realcompactification and the Samuel compactification coincide.

\begin{theorem}  Let $(X,\mu)$ be a uniform space. Then $H(U_{\mu}(X))=s_{\mu}X$ if and only if $U_{\mu}(X)=U^{*}_{\mu}(X)$, that is, if and only if $(X,\mu)$ is Bourbaki-bounded.
\begin{proof} If $H(U_{\mu}(X))= s_{\mu} X$ then $H(U_{\mu}(X))$ is compact and then every $f\in C(H(U_{\mu}(X)))$ is bounded. In particular since $U_{\mu}(X)\subset C(H(U_{\mu}(X)))$, then, $U_{\mu}(X)=U_{\mu}^{*}(X)$, that is, $X$ is Bourbaki-bounded \cite{hejcman}. The converse is trivial.
\end{proof}

\end{theorem}

Thus, every realcompact Bourbaki-bounded space, which fails compactness is not Samuel realcompact. For instance, take the metric hedgehog (Example \ref{hedgehog}) $(H(\omega _0),\rho)$ which is realcompact because in particular it is Lindel\"{o}f.




\smallskip

Now, we finish  with some easy results in this topic, that we will use later.

\begin{theorem} \label{eq.Samuel} Let $(X,\mu)$ be a uniform space. For every uniformity $\nu$ on $X$ satisfying that $w_{U_{\mu}(X)} \preceq \nu \preceq \mu$ it is satisfied that $H(U_{\mu}(X))=H(U_{\nu}(X))$. 
\begin{proof} It follows easily from the fact that $U_{\mu}(X) =U_{w_{U_{\mu}(X)}} (X)\subset U_{\nu}(X)\subset U_{\mu}(X)$. 
\end{proof}
\end{theorem}

\begin{theorem} \label{eq.Samuel.complete} Let $(\widetilde{X},\widetilde{\mu})$ be the completion of a uniform space $(X,\mu)$. Then $H(U_{\mu}(X))=H(U_{\widetilde{\mu}}(\widetilde{X}))$. 
\begin{proof} It is known that the functions in $U_{\mu}(X)$ are exactly the restrictions of the real-valued uniformly continuous functions of the completion $(\widetilde{X}, \widetilde{\mu})$ (see \cite{willard}). Thus, by density of $X$ is $\widetilde{X}$, the result follows. 
\end{proof}
\end{theorem}

\medskip

\subsection{Samuel realcompact spaces}

\hspace{15pt}In this section we solve the problem of characterizing those uniform spaces $(X,\mu)$ which are Samuel realcompact, that is, those spaces satisfying that $X=H(U_{\mu}(X))$, or equivalently, that $(X,wU_{\mu}(X))$ is complete. This result is strictly related to the classical Kat\v{e}tov-Shirota Theorem.

\begin{theorem}{\rm (\cite{shirota1} and \cite{katetov} (only for paracompact spaces))} A space $X$ is realcompact if and only it is topologically complete and there is no closed discrete subspace of Ulam-measurable cardinal.
\end{theorem}

Next, recall the definition of Ulam-measurable cardinal.

\begin{definition} A filter $\mathcal{F}$ on a set $S$ satisfies the {\it countable intersection property}  if for every countable family $\{F_n:n\in \Nset\}\subset \mathcal{F}$, $$\bigcap_{n\in \Nset} F_n\neq \emptyset.$$
\end{definition}

\noindent In particular, if an ultrafilter $\mathcal{F}$ satisfies the countable intersection porperty, then, by maximality, $\bigcap_{n\in \Nset}F_n \in \mathcal{F}$, that is, $\mathcal{F}$ is {\it closed under countable intersection}. 

\begin{definition} An infinite cardinal $\kappa$ is Ulam-measurable if there is a free (non-principal) ultrafilter $\mathcal{F}$ (that is, $\bigcap \mathcal{F}=\emptyset$) satisfying the countable intersection property on any set of cardinal $\kappa$. 
\end{definition}

\noindent As we have explained in the introduction, the above definition is equivalent to say that any discrete space of cardinal $\kappa$ is not realcompact. Moreover given two cardinals $\kappa \leq \lambda $, if $\kappa$ is Ulam-measurable then $\lambda$ is also Ulam-measurable.

\smallskip
 
We are going to see now that  a Kat\v{e}tov-Shirota type theorem characterizing Samuel realcompactness can be obtained where Bourbaki-completeness will play the role of topological completeness in the classical one. The proof that we give here is different than the proofs in \cite{merono.Samuel.metric}, \cite{merono.Samuel.uniform} and even in \cite{husek} where a categorical method is used. It uses the embedding results from Part 2.

\begin{lemma} \label{function}If $f:(X,\mu)\rightarrow (Y,\nu)$ is uniformly continuous then $$f:(X, wU_{\mu} (X))\rightarrow (Y,wU_\nu (Y))$$ is uniformly continuous.
\begin{proof} That $f:(X,\mu)\rightarrow  (Y,wU_\nu (Y))$ is uniformly continuous is clear.  Now, recall that linear covers are a subbase for the weak uniformity $wU_\nu (Y)$. Thus, take a linear cover $\mathcal{V}=\{V_n:n\in \Nset\}\in wU_\nu (Y)$ then $f^{-1}(\mathcal{V})=\{f^{-1}(V_n): n\in \Nset\}\in \mu$ by uniform continuity of $f$. Next, it is clear that $f^{-1}(\mathcal{V})$ is countable. Moreover, $f^{-1}(V_n)\cap f^{-1}(V_m)=f^{-1}(V_n\cap V_m)$ is always satisfied. Therefore, $f^{-1}(\mathcal{V})$ is linear because $f^{-1}(V_n)\cap f^{-1}(V_m)=\emptyset$ whenever $|n-m|>1$.
\end{proof}
\end{lemma}

\begin{theorem} \label{subspace}Let $(X,\mu)$ be a uniform space and $(Y,\nu)$ be a Samuel realcompact space such that there exists  an embedding $$\varphi :(X,\mu)\rightarrow (Y,\nu)$$ such that $\varphi$ is uniformly continuous and $\varphi(X)$ is closed in $Y$. Then $(X,\mu)$ is Samuel realcompact.
\begin{proof} We are going to prove that $(X, wU_{\mu}(X))$ is complete. Let $\mathcal{F}$ be a Cauchy filter of $(X, wU_{\mu}(X))$. By Lemma \ref{function} then $\varphi (\mathcal{F})$ is a Cauchy filter of $(Y, wU_\nu(Y))$. By Samuel realcompactness, $\varphi (\mathcal{F})$ clusters in $Y$.  Since $\varphi (X)\in \varphi(\mathcal{F})$ and $\varphi(X)$ is a closed subspace then $\varphi (\mathcal{F})$ cluster in $\varphi(X)$. Finally $\mathcal{F}$ clusters in $X$ because $\varphi$ is an embedding.
\end{proof}
\end{theorem}

\begin{theorem} \label{product.Samuel} Let $(X_i,\mu_i)$, $i\in I$, be uniform spaces. Then $$(\prod_{i\in I} X_i,\prod_{i\in I}\mu_i)$$ is Samuel realcompact if and only if each factor $(X_i,\mu_i)$ is Samuel realcompact.
\begin{proof} $\Rightarrow)$ This implication is clear from Theorem \ref{subspace} since each factor $(X_i,\,\mu_i)$ is uniformly homeomorphic to a closed subspace of $(\prod_{i\in I} X_i,\prod_{i\in I}\mu_i)$.

$\Leftarrow)$  We prove that $(\prod_{i\in I} X_i, wU_{\prod_{i\in I}\mu _i}(\prod_{i\in I}X_i))$ is complete. By Lemma \ref{function} the projection maps $$p_i : (\prod_{i\in I} X_i, wU_{\prod_{i\in I}\mu _i}(\prod_{i\in I}X_i))\rightarrow (X_i,wU_{\mu_i}(X_i))$$ are uniformly continuous. Let $\mathcal{F}$ be a Cauchy ultrafilter of $$(\prod_{i\in I} X_i, wU_{\prod_{i\in I}\mu _i}(\prod_{i\in I}X_i))$$ then $p_i (\mathcal{F})$ is a Cauchy ultrafilter of $(X_i,wU_{\mu_i}(X_i)).$ By Samuel realcompactness $p_i (\mathcal{F})$ converges to some $x_i\in X_i$. Then, $\mathcal{F}$ converges to the the point $x=(x_i)\in \prod_{i\in I} X_i$ by maximality.
\end{proof}
\end{theorem}

The next theorem characterizes those uniform spaces having no uniform partition of Ulam-measurable cardinal. This is the clue to our Kat\v{e}tov-Shirota Theorem \ref{measurable2} characterizing Samuel realcompact uniform spaces.

\begin{theorem} \label{measurable1}  Let $(X,\mu)$ be a uniform space and let $(\widetilde{X}, \widetilde{s_{f}\mu})$ be the completion of  $(X, s_f\mu)$. The following statements are equivalent:
\begin{enumerate}

\item There is no uniform partition of $(X,\mu)$ having Ulam-measurable cardinal;

\item $\widetilde{X}=H(U_{\mu}(X))$  (that is, $(\widetilde{X}, \widetilde{s_{f}\mu})$ is Samuel realcompact);

\item $\widetilde{X}$ is realcompact;

\item there is no uniformly discrete subspace of $(X, s_f\mu)$ having Ulam-measurable cardinal;

\item there is no uniform partition of $(X,s_f\mu)$ having Ulam-measurable cardinal.

\end{enumerate}
\begin{proof} $(1)\Rightarrow (2)$ By Lemma \ref{eq.Samuel} and Lemma \ref{eq.Samuel.complete} $$H(U_\mu (X))=H(U_{s_{f}\mu}(X))=H(U_{\widetilde{s_f\mu}}(\widetilde{X})).$$ Therefore it is enough to prove that $\widetilde{X}= H(U_{\widetilde{s_f\mu}}(\widetilde{X}))$, that is, that $(\widetilde{X}, \widetilde{s_{f}\mu})$ is Samuel realcompact. 

By \cite[Lemma p. 370] {reynolds} the uniformity $\widetilde{s_{f}\mu}$ has a star-finite base. Then applying Theorem \ref{embedding.star-finite2}, there exists an embedding $$\varphi:(\widetilde{X}, \widetilde{s_f\mu})\rightarrow  \Big((\prod_{i\in I, n\in \Nset} \kappa_n ^i) \times \Rset^{\alpha},\pi\Big)$$ where each cardinal $\kappa _n  ^i$ satisfies that $\kappa _n ^i\leq \wp((\widetilde{X}, \widetilde{s_f\mu}))$. It is clear that $\wp((\widetilde{X}, \widetilde{s_f\mu}))=\wp((X, s_f\mu))=\wp ((X,\mu))$. Therefore each cardinal $\kappa _n ^i$ has no Ulam-measurable cardinal by hypothesis and then by Theorem \ref{R-Samuel.realcompact} and Theorem \ref{product.Samuel}, the space $\Big((\prod_{i\in I, n\in \Nset} \kappa_n ^i )\times \Rset^{\alpha},\pi\Big)$ is Samuel realcompact. Thus by Theorem \ref{subspace} it follows that $(\widetilde{X}, \widetilde{s_{f}\mu})$ is Samuel realcompact.

$(2)\Rightarrow (3)$ This implication is trivial.

$(3)\Rightarrow (4)$ Since every uniformly discrete subspace of $(X, s_f\mu)$ is a closed discrete subspace of the completion $(\widetilde{X}, \widetilde{s_f\mu})$, then it follows that there is no uniformly discrete subspace of $(X, s_f\mu)$ having Ulam-measurable cardinal.

$(4)\Rightarrow (5)$. This follows taking into account that every uniform partition determines a uniformly discrete subspace with the same cardinal.

$(5)\Rightarrow (1)$. It is clear since  $(X,\mu)$ and $(X, s_f\mu)$ share the same uniform partitions.
\end{proof}
\end{theorem}

From the last result we deduce easily our Kat\v{e}tov-Shirota type Theorem.

\begin{theorem} \label{measurable2} {\rm ({\sc Kat\v{e}tov-Shirota type theorem)}} A uniform space $(X,\mu)$ is Samuel realcompact if and only if it is Bourbaki-complete and there is no uniform partition having Ulam-measurable cardinal.
\begin{proof} If $(X,\mu)$ is Samuel realcompact, that is, $(X, wU_{\mu}(X))$ is complete, then it is clear that $(X, s_f\mu)$ is complete as $wU_{\mu}(X)\leq s_f \mu$. Therefore, by Theorem \ref{star-finite2}, $(X,\mu)$ is Bourbaki-complete. On the other hand,  Samuel realcompactness implies realcompactness. Hence, given a uniform partition $\mathcal{P}=\{P_i :i\in I\}$ of $(X,\mu)$, the cardinal of $I$ is no Ulam-measurable. Indeed, if we take representative points $x_i \in P_i$ for every $i\in I$, then the subspace $\{x_i:i \in I\}$ of $X$ is closed and discrete. Hence, by realcompactness, $\{x_i:i \in I\}$ has no Ulam-measurable cardinal, and so does $I$.

Conversely, let $(\widetilde{X}, \widetilde{s_f\mu})$ denote the completion of $(X,s_f\mu)$. Then, by Theorem \ref{measurable1},  $\widetilde{X}=H(U_{\mu}(X))$. Next, by Bourbaki-completeness (Theorem \ref{star-finite2}), the space $(X,\mu)$ is homeomorphic with the completion  $(\widetilde{X}, \widetilde{s_f\mu})$. Therefore $X=H(U_\mu (X))$, that is, $(X,\mu)$ is Samuel realcompact.
\end{proof}
\end{theorem}

As a corollary of the above theorem we have the following result for uniformly 0-dimensional spaces from \cite{husekpulga}. Recall that uniformly 0-dimensional spaces are those spaces having a base of partitions for their uniformity, or equivalently, that are uniform subspaces of a product of uniformly discrete spaces (see Remark \ref{0-dim}). 

\begin{corollary} \label{pulga} A uniformly 0-dimensional space is Samuel realcompact if and only if it is complete and it does not have a  uniformly discrete subset of Ulam-measurable cardinal.
\begin{proof} Every complete uniformly 0-dimensional space is Bourbaki-complete because it has a base of partitions for its uniformity (hence, $\mu=s_f\mu$). By the same reason, uniformly discrete subsets determine the uniform partitions. Therefore, the result follows from Theorem \ref{measurable2}.
\end{proof} 
\end{corollary}

We must notice now that the above Kat\v{e}tov-Shirota type theorem can also be  found in the paper by Rice and Reynolds \cite{reynolds}. However their result is stated as follows.

\begin{theorem} \label{reynolds} {\rm (\cite{reynolds})} Let $(X,\mu)$ be a uniform space such that no uniform cover has Ulam-measurable cardinal. Then $(X, w_{U_\mu(X)})$ is complete if and only if $(X, s_f\mu)$ is complete.
\end{theorem}

Observe that Theorem \ref{star-finite2}, that is, the equivalence between Bourbaki-completeness of a uniform space and completeness of the star-finite modification, is the link between our Kat\v{e}tov-Shirota result (Theorem \ref{measurable2}) and the result by Rice and Reynolds. 
However, we point out now a slight difference between them. It is clear that only the non Ulam-measurability of the uniform partitions is needed  in last Theorem \ref{reynolds}, and not the (stronger) condition of the non Ulam-measurability  of the cardinality of the uniform covers. In fact Theorem \ref{measurable1} characterizes precisely uniform spaces having no uniform partition of Ulam-measurable cardinal, which is not in general equivalent to satisfy that no uniformly discrete subspace has Ulam-measurable cardinal, even if this equivalence is true whenever we consider the uniformity $s_f\mu$. Examples of uniform spaces such that no uniform partition has Ulam-measurable cardinal but having uniform covers of Ulam-measurable cardinal are easy to find. Just think on connected spaces containing a uniformly discrete subspace of Ulam-measurable cardinal as, for instance, the metric hedgehog $H(\kappa)$ or any $\ell_\infty(\kappa)$ where $\kappa$ is an Ulam-measurable cardinal.

\smallskip

Now, if we apply Theorem \ref{measurable2} to the fine uniformity ${\tt u}$ on a space $X$ we obtain the following result which is a refinement of the classical Kat\v{e}tov-Shirota Theorem. 

\begin{corollary} \label{garcia} A space $X$ is realcompact if and only if $X$ is $\delta$-complete and there is no open partitions of $X$ having Ulam-measurable cardinal.
\end{corollary}

\noindent Observe that from the above result we can also deduce Corollary \ref{connected.realcompact}.


\smallskip

\medskip

\subsection{Uniformities having a base of countable uniform covers}

\hspace{15pt} In his famous paper  \cite{shirota1}, Shirota proved that the Hewitt-Nachbin realcompactification $\upsilon X$ of a space $X$ is homeomorphic to the completion of $(X, e{\tt u})$ where $e{\tt u}$ is the uniformity having as a base all the cozero countable covers of $X$. Since for the fine uniformity ${\tt u}$, $C(X)=U_{{\tt u}}(X)$, we have that the completion of $(X,wU_{{\tt u}}(X))$ and the completion of $(X, e{\tt u})$ are homeomorphic.

There is another class of uniform spaces $(X,\mu)$ satisfying that the Samuel realcompactification $H(U_{\mu}(X))$ is homeomorphic to the completion of $(X, e\mu)$, precisely, the uniformly 0-dimensional  spaces. Indeed, for these spaces $e\mu=wU_{\mu}(X)$ as 0-dimensional spaces have a base of partitions for their uniformity.

In general, for a uniform space $(X,\mu)$, it is not true that the completion of $(X, wU_{\mu}(X))$ and the completion $(X,e\mu)$ are homeomorphic. Indeed, there are complete metric spaces which are separable, and therefore, satisfying that their metric uniformity has a base of countable covers, but which are not Bourbaki-complete, and hence not Samuel realcompact by Theorem \ref{measurable2}. For instance the metric hedgehog $(H(\omega _0),\rho)$.

That the completion of $(X,e\mu)$ is a realcompactification of $X$ follows from the fact that it is also complete with the uniformity having as a base all the countable uniform covers of its fine uniformity, and therefore, by the Shirota result,  it is  realcompact.  Moreover, it represents a uniform generalization of realcompactness to the frame of uniform spaces as the Samuel realcompactification does (\cite{husek}). Precisely it is the realcompactification of $(X,\mu)$ generated by the subalgebra of $C(X)$ given by all the Cauchy-continuous functions of $(X, e\mu)$ (see Theorem \ref{cc}).

In the next result we show how the Samuel realcompactification $H(U_{\mu}(X))$ of a uniform space and the completion of $(X,e\mu)$ are related even if they are different.

\begin{theorem} \label{countable1} Let $(X,\mu)$ be a uniform space. The following spaces are (topologically) homeomorphic:
\begin{enumerate}
\item The Samuel realcompactification $H(U_\mu (X))$;

\item  the subspace of $s_{\mu} X$ of all the cluster points of the Bourbaki-Cauchy filters of $(X, e\mu)$;

\item  the completion of $(X, s_f(e\mu))$;

\item  the subspace of $s_{\mu} X$ of all the cluster points of the Bourbaki-Cauchy filters of $(X,e(s_f\mu))$.

\end{enumerate}
\begin{proof} Applying Theorem \ref{eq.Samuel} several times, we have that $$H(U_\mu (X))=H(U_{e\mu} (X))=H(U_{s_f\mu} (X))=H(U_{s_f(e\mu)} (X))=H(U_{e(s_f\mu)} (X)).$$ The same is true for the Samuel compactification.  Since every uniform partition of $(X,e\mu)$ is countable then, by Theorem \ref{measurable1}, the completion of $(X, s_f(e\mu))$ coincides with its Samuel realcompactification which is exactly $H(U_{\mu}(X))$ as we have previously shown. Therefore $(1)$ and $(3)$ are homeomorphic. Next, by Theorem \ref{star-finite1}, the subspace of $s_{\mu} X$ of all the cluster points of the Bourbaki-Cauchy filters of $(X, e\mu)$ is homeomorphic to the completion of  $(X, s_f(e\mu))$. Hence, $(2)$ and $(3)$ are homeomorphic. Finally applying similarly Theorem \ref{measurable1} and Theorem \ref{star-finite1} to the uniform space $(X,e(s_f\mu))$, the homeomorphisms of $(1)$ and $(4)$ follows also at once. 
\end{proof}
\end{theorem}

\begin{theorem} \label{countable2}Let $(X,\mu)$ be a uniform space. The following statements are equivalent:
\begin{enumerate}
\item $X$ is Samuel realcompact.

\item $(X, e\mu)$ is Bourbaki-complete.

\item $(X, s_f (e\mu))$ is complete.

\item $(X,e(s_f\mu))$ is complete.
\end{enumerate}
\begin{proof} We only need to prove $(4) \Rightarrow (1)$. If $(X,e(s_f\mu))$ is complete then $(X, s_f\mu)$ is complete since the identity map $id:(X, s_f\mu)\rightarrow (X,e(s_f\mu))$ is uniformly continuous. In addition, $X$ is a realcompact space since, as we have said before, the completion of every uniformity of type $e\mu$ is a realcompactification of $X$. Thus, every uniform cover of $(X, s_f\mu)$ has no Ulam-measurable cardinal and then, by Theorem \ref{reynolds}, $X$ is Samuel realcompact.
\end{proof}
\end{theorem}

In \cite{reynolds}, Rice and Reynolds characterize those uniform spaces $(X,\mu)$ satisfying that $(X,e\mu)$ is complete.

\begin{theorem} \label{reynolds2} {\rm (\cite{reynolds})} Let $(X,\mu)$ be a uniform space. If there is no uniform cover of $(X,\mu)$ having Ulam-measurable cardinal, then $(X,e\mu)$ is complete if and only if $(X, p_f\mu)$ is complete.
\end{theorem} 

\noindent By the results  of Pelant in \cite{pelant.complete} (see also \cite{husek} for more information), it is known that the Banach space $(\ell _\infty (\omega _1), ||\cdot||_1)$  satisfies that its point-finite modification  $(\ell _\infty  (\omega _1), p_f \mu _{||\cdot||_1})$, where $\mu _{||\cdot||_1}$ denotes the uniformity induced by the norm,  is not complete, even if $(\ell _\infty  (\omega _1), ||\cdot||_1)$ is complete. Therefore, we can deduce from Theorem \ref{reynolds2} that not every complete uniform space, satisfying that every uniform cover has no Ulam-measurable cardinal, is  complete whenever it is endowed with the countable modification $e\mu$.

Finally, observe that Reynolds and Rice asked if it is possible to prove their result Theorem \ref{reynolds} from Theorem \ref{reynolds2}. The problem was that they did not know if  the modification $e(s_f\mu)$ has a base of star-finite covers, or more precisely, if $e(s_f\mu)\leq s_f (e(s_f\mu))\leq s_f (e\mu)$, as we don't know either. Related to the above question of Rice and Reynolds, observe that in condition $(4)$ of  Theorem \ref{countable1}, it would be interesting to be able to consider the completion of $(X, e(s_f \mu))$  instead of the set of cluster points of its Bourbaki-Cauchy filters, as it happens in condition $(4)$ of  Theorem \ref{countable2}. However, since we cannot assure that the uniformity $e(s_f\mu)$  has a star-finite base we cannot refine the result. Thus, according to Theorem \ref{countable2}, we wonder if the difference between the above results lies in the non Ulam-measurability of the cardinality of the uniform covers. Nevertheless, for uniformly 0-dimensional we have that $e(s_f\mu)=e\mu=wU_{\mu}(X)=s_f(e\mu).$

\bigskip

\bigskip

\begin{center} \Denarius \Denarius  \Denarius \end{center}
\bigskip

\bigskip

\section{A step forward: relating the Samuel realcompactification and the Hewitt realcompactification. }

\subsection{Certain family of real-valued continuous functions}

\hspace{15pt} The first problem that we want to solve in this section is to know if every metric space $(X,d)$ is Bourbaki-complete, whenever the Hewitt realcompactification and the Samuel realcompactification of $(X,d)$ are equivalent. Recall that in \cite{merono.Samuel.metric} it was already proved the  necessity of completeness. However, this fact is clearly false in the general frame of uniform spaces since there are fine uniform spaces which are not even complete. For instance, the space $([0,\omega _1),{\tt u})$ is not complete but $\upsilon ([0,\omega _1))=[0, \omega _1]=H(U_{\tt u} ([0,\omega _1))).$
\smallskip

To that purpose, we are going to characterize Bourbaki-completeness by means of certain family of real-valued functions. Recall that by Theorem \ref{dimaio} a metric space $(X,d)$ is complete if and only if every real-valued continuous functions preserves Cauchy filters of $(X,d)$, that is, $f$ is Cauchy-continuous. Next, let us denote by $CC_{\mu}(X)$ the subalgebra of all the Cauchy-continuous functions on a uniform space $(X,\mu)$. By Theorem \ref{star-finite2} and Theorem \ref{dimaio}, every Bourbaki-complete uniform space satisfies that $$CC_{s_f \mu} (X)=C(X).$$ 

Now we ask it the converse is true for a metric space $(X,d)$, that is, if $CC_{s_f \mu_d} (X)=C(X)$ is satisfied, then $(X,d)$ is Bourbaki-complete.

\begin{theorem} \label{Bourbaki-complete.functions} For a metric space $(X,d)$ the following statements are equivalent:

\begin{enumerate}
\item $(X,d)$ is Bourbaki-complete;

\item $C(X)=CC_{s_f \mu_d}(X)$


\end{enumerate}
\begin{proof}
$(1) \Rightarrow (2)$ This follows at one by Theorem \ref{star-finite2} and Theorem \ref{dimaio}. 
\smallskip

$(2)\Rightarrow (1)$ If $CC_{s_f\mu_d}(X)=C(X)$ then $CC_{\mu_d}(X)=C(X)$ because, since $s_f \mu _d \leq \mu _d$, then $CC_{s_f\mu_d}(X)\subset CC_{\mu_d}(X)$. Hence, by Theorem \ref{dimaio}, the metric space $(X,d)$ must be complete. 

Assume by contradiction that $(X,d)$ is not Bourbaki-complete and take $B$ a Bourbaki-bounded subset in $X$ which fails to be totally bounded (see Theorem \ref{metric.BC}). Then, for some $\delta >0$, $B$ contains a subset $\{x_n:n\in \Nset\}$ such that $B_{d}(x_{n}, \delta )\cap B_{d}(x_{j}, \delta )=\emptyset$ for every $n,j\in \Nset$, $n\neq j$. Let us define  $$f(x)=\begin{cases}  \,\,   n-\frac{n}{\delta}\cdot d(x, x_n) &     \text {if $d(x,x_n)<\delta$ for some  $n\in \Nset$;}  \\ \,\,\, 0  & \text{ otherwise.}\end{cases}$$ Then $f$ is  continuous because it is uniformly continuous when restricted over every ball $B_{d}(x, \delta)$, $x\in X$. 

Moreover, the set $S=\{x_n:n\in \Nset\}$ is Bourbaki-bounded in $(X,d)$ and, in particular, it is a totally bounded subset of $(X, s_f\mu _d)$ (see Theorem \ref{BB.sf}). Next, consider the filter base $\{\{x_{j}:j\geq n\}: n\in \Nset \}$ and let $\mathcal{F}$ the filter of $S$ induced by it. By Theorem \ref{totally.bounded.filters}, $\mathcal{F}$ is contained in some filter $\mathcal{F}'\subset S$ which is Cauchy in $(X,s_f\mu _d)$. Observe that $\mathcal{F}'$ does not cluster because $\mathcal{F}$ does not cluster either. Thus, since $f$ is unbounded on $S$,  it cannot map the filter $\mathcal{F}'$ to a Cauchy filter of $(\Rset, d_u)$. Hence, $f$ does not belong to $CC_{s_f\mu_d}(X)$ contradicting the hypothesis (see \cite{garrido-beer1} and \cite{garrido-beer2} for similar techniques.)



\end{proof}
\end{theorem}


\smallskip

\smallskip

Now consider the Samuel realcompactification $H(U_{\mu}(X))$ of a uniform space $(X,\mu)$. By Theorem \ref{cc}, the subalgebra  $C(H(U_{\mu}(X)))$ is isomorphic to the subalgebra of Cauchy continuous functions $CC_{wU_{\mu}(X)}(X)$. Moreover, observe that in general, for a uniform space $(X,\mu)$, $$ C(H(U_{\mu}(X)))=CC_{wU_{\mu}(X)}(X)\subset CC_{s_f\mu}(X)$$ because $wU_\mu (X) \leq s_f \mu$. 

By all the foregoing and Theorem \ref{Bourbaki-complete.functions}, the following theorem is immediate.

\begin{theorem}\label{realcompactifcations->Bourbaki-complete} Let $(X,\mu)$ be a uniform space satisfying that the realcompactifications $\upsilon X$ and $H(U_{\mu}(X))$ are equivalent. Then $C(X)=CC_{wU_{\mu}(X)}(X)$. In particular, for a metric space, this implies that the space is Bourbaki-complete.
\end{theorem}

\medskip

\subsection{A positive result}

\hspace{15pt} Now, we want solve the converse problem form the previous section, that is, if for every Bourbaki-complete uniform space the Hewitt realcompactification and the Samuel realcompactification are equivalent. This is clearly true in the frame of Bourbaki-complete uniform spaces  having no uniform partition of Ulam-measurable cardinal, since by Samuel realcompactness (Theorem \ref{measurable2})  both realcompactifications are homeomorphic to the space. However, we are interested here in examples and results not implying Samuel realcompactness. For instance, any uniformly discrete metric space, independently of its cardinality satisfies that both realcompactifications are equivalent.  

We start next, by analyzing the case of  the product of two uniform spaces, one of them being Samuel realcompact. From it we will derive a first answer in the frame of metric spaces.

\begin{lemma}\label{lemma1f} Let $(X,\mu)$ be a Samuel realcompact space and $(Y,\nu)$ be a uniform space. Then every uniformly continuous function $f\in U_{\mu \times \nu}(X\times Y)$ can be continuously extended to $X\times H(U_{\nu}(Y))=H(U_{\mu}(X))\times H(U_{\nu}(Y))$.
\begin{proof} (See \cite{woods} for similar techniques.) Observe that for every $x\in X$ the function $f_x :(Y,\nu)\rightarrow ( \Rset, d_u)$ defined by $f_x(y)=f(x,y)$, for every $y\in Y$, is uniformly continuous and hence, it can be extended to a unique continuous function $f^{*}_x: H(U_\nu (Y))\rightarrow \Rset$. Define $$f^{*}:X\times H(U_\nu (Y))\rightarrow \Rset$$ by $f^{*}(x,\xi)=f^{*}_x (\xi)$. We are going to prove that $f^*$ is continuous. To that purpose it is enough to prove that $f^*$ is continuous on $X\times Y \cup \{(x_0, \xi)\}$ for every $x_0\in X$ and $\xi \in H(U_{\nu}(Y))\backslash Y$ (see \cite[Exercise 6H]{gillman}).

Fix $x_0 \in X$, $\xi \in H(U_{\nu}(Y))\backslash Y$, and $\varepsilon >0$. Since $f^{*}_{x_0}$ is continuous  there exists $V^\xi\subset H(U_{\nu}(Y))$ open set such that $\xi \in V^\xi$ and $$f^{*}_{x_0}(V^\xi)\subset B_{d_u}(f^{*}_{x_0}(\xi),\varepsilon/2).$$ Moreover, if $y\in V^{\xi}\cap Y$ then $f(x_0,y)\in  B_{d_u}(f^{*}(x_0,\xi),\varepsilon/2)$.

On the other hand, as $f$ is uniformly continuous, there exists a uniform cover $\mathcal{U}\times \mathcal{W}\in \mu \times \nu$ such that $|f(x,y)-f(u,w)|<\varepsilon/2$ whenever $(x,y), (u,w)\in U\times W$ for every $U\times W\in \mu \times \nu$. 

Now, take $U\in \mathcal{U}$ such that $x_0\in U$ and $W\in \mathcal{W}$ such that $W\cap V^{\xi}\neq \emptyset$. Then, if $(x,y)\in (U\times ( V^{\xi} \cap Y))$ we have that by all the foregoing $$|f(x,y)-f^* (x_0,\xi)|\leq|f(x,y)-f(x_0,y)|+|f (x_0,y)-f^* (x_0,\xi)|<\varepsilon/2+\varepsilon/2=\varepsilon$$ and $f^{*}$ is continuous in $(x_0,\xi)$ as we claimed. 
\end{proof}
\end{lemma}

As a consequence of the above lemma we have that for a Samuel realcompact space $(X,\mu)$ and for any uniform space $(Y,\nu)$ the following is satisfied:

$$H(U_{\mu\times \nu}(X\times Y))\leq H(U_{\mu}(X))\times H(U_{\nu}(X)).$$ Next we see that the reverse inequality is, in general, always satisfied.

\begin{lemma} \label{lemma3f} Let $(X,\mu)$ and $(Y,\nu)$ be uniform spaces. Then $$H(U_{\mu\times \nu}(X\times Y))\geq H(U_{\mu}(X))\times H(U_{\nu}(Y)).$$
\begin{proof}
Let $i:(X\times Y, wU_{\mu \times \nu}(X\times Y))\rightarrow (X\times Y, wU_{\mu}(X)\times wU_{\nu} (Y))$ be the identity map then, it is uniformly continuous. Indeed, the projections maps $$p_1:((X\times Y, wU_{\mu \times \nu}(X\times Y))\rightarrow (X, wU_{\mu}(X))$$ and $$p_2:((X\times Y, wU_{\mu \times \nu}(X\times Y))\rightarrow (Y, wU_{\nu}(Y))$$ are uniformly continuous by Lemma \ref{function} and hence, $i$ is uniformly continuous as the product uniformity $wU_{d}(X)\times wU_{\nu}(Y)$ is the weakest uniformity on $X\times Y$ making the above projections onto $(X, wU_{\mu}(X))$ and $(Y, wU_{\nu}(Y))$ uniformly continuous (see \cite{willard}). Therefore, $wU_{\mu}(X)\times wU_{\nu}(Y) \leq wU_{\mu \times \nu}(X\times Y)$, that is, the identity map $$id:(X\times Y, wU_{\mu \times \nu}(X\times Y)) \rightarrow (X \times Y, wU_{\mu}(X)\times wU_{\nu}(Y))$$ is uniformly continuous. Hence, the the map $id$ can be extended to a continuous map between the  completions of both spaces, that is, $H(U_{\mu \times \nu}(X\times Y))\geq H(U_{ \mu}(X))\times H(U_{\nu}(Y))$.
\end{proof}
\end{lemma}

\begin{theorem} \label{Samuel.product} Let $(X,\mu)$ be a Samuel realcompact space and $(Y,\nu)$ any uniform space. Then $$H(U_{\mu \times \nu}(X\times Y))=H(U_{\mu}(X))\times H(U_{\nu}(Y))=X\times H(U_{\nu}(Y)).$$
\begin{proof} The proof follows at once by Lemma \ref{lemma1f} and Lemma \ref{lemma3f}.
\end{proof}
\end{theorem}

Next, we study the problem for a product of a Samuel realcompact space and a uniformly discrete space. To that purpose we recall the following result by Hu\v{s}ek.

\begin{theorem}\label{product.husek} {\rm(\cite[Theorem 3]{husek.product})}
 Let $D$ be a discrete space. Then $\upsilon (D\times X)=\upsilon D \times \upsilon X$ if and only if either $D$ or $X$ do not have Ulam-measurable cardinal.

\end{theorem}

\begin{theorem} \label{product.Samuel.discrete} Let $(X,\mu)$ be a Bourbaki-complete uniform space having no Ulam-measurable cardinal, and $(D,\chi)$ be any uniformly discrete space. Then $$H(U_{\mu\times \mu_\chi} (X\times D))=H(U_{\mu}(X))\times H(U_{\chi}(D))=X\times \upsilon D=\upsilon X\times \upsilon D= \upsilon (X\times D).$$ 
\begin{proof}
Since $(X,\mu)$ is Bourbaki-coomplete and has not Ulam-measurable cardinal, $(X,\mu)$ is Samuel realcompact. Moreover, $(D,\xi)$ is a uniformly discrete space, in particular, a UC space. Then  by Theorem \ref{product.Samuel} and Theorem \ref{product.husek} $$H(U_{\mu\times \mu_\chi} (X\times D))=H(U_{\mu}(X))\times H(U_{\chi}(D))=X\times \upsilon D=\upsilon X\times \upsilon D= \upsilon (X\times D).$$

\end{proof}
\end{theorem}

\begin{corollary}\label{product.Samuel.discrete2} Let $(X,d)$ be a Samuel realcompact metric space and $(D,\chi)$ be any uniformly discrete space. Then $$H(U_{d+\chi} (X\times D))=H(U_{d}(X))\times H(U_{\chi}(D))=X\times \upsilon D=\upsilon X\times \upsilon D= \upsilon (X\times D).$$ 

\begin{proof} If $(X,d)$ is any Samuel realcompact metric space, then by Theorem \ref{measurable2}, $(X,d)$ is Bourbaki-complete  and there is no uniform partition having Ulam-measurable cardinal. Moreover, by Theorem \ref{embedding.metric} there exists an embedding $$\varphi :(X,d)\rightarrow (\prod_{ n\in \Nset} \kappa _n\times\Rset^{\omega _0} , \rho+t)$$ where no cardinal $\kappa _n$ is Ulam-measurable, $\varphi$ is uniformly continuous and $\varphi(X)$ is a closed subspace of $\prod_{ n\in \Nset} \kappa _n\times\Rset^{\omega _0}$. Then, $X$ does not have Ulam-measurable cardinal and we can apply Theorem \ref{product.Samuel.discrete}.
\end{proof}
\end{corollary}

From the above result and the embedding Theorem \ref{embedding.metric} we are going to characterize the class of Bourbaki-complete metric spaces satisfying that the Samuel realcompactification and the Hewitt realcompactification are equivalent. But first, we introduce a couple of useful lemmas.

\begin{lemma} \label{subspace.lemma} Let $(Y,\mu)$ be a uniform space and $X\subset Y$. Then $$H(U_{\mu|_X}(X))\subset H(U_\mu (Y))\cap s_{\mu|_X} X.$$
\begin{proof} By Lemma \ref{Samuel3} and Theorem \ref{realcompact.subspace} $$X\subset H(U_{\mu|_X}(X))\subset s_{\mu|_X} X={\rm cl}_{s_{\mu}Y} X\subset s_{\mu}Y.$$

Suppose that there is some $\xi \in H(U_{\mu|_X}(X))$ such that $\xi \notin H(U_\mu (Y))$. Then for some $F \in U_\mu (Y)$, $F^{*}(\xi)=\infty$. But $F|_X \in U_{\mu|_X}(X)$ and by uniqueness of the extension $(F|_X)^{*}(\xi)=F^{*}(\xi)=\infty$, which is a contradiction as $\xi \in H(U_{\mu|_X}(X))$.
\end{proof}
\end{lemma}

\begin{lemma}\label{mapping.Samuel} Let $(Y,\rho)$ be a metric space satisfying that $\upsilon Y$ and $H(U_\rho(Y))$ are equivalent realcompactifications. Let $(X,d)$ be a metric space such that there exists an embedding $\varphi:(X,d)\rightarrow (Y,\rho)$ satisfying that $\varphi$ is uniformly continuous and $\varphi(X)$ is a closed subspace of $Y$. Then $\upsilon X$ and $H(U_d(X)))$ are also equivalent realcompactifications.
\begin{proof} Let $f\in C(X)$ and identify $X$ with its image $\varphi(X)$. We are going to prove that $f$ can be continuously extended to $H(U_{d}(X))$. Indeed, since $X$ is a closed subspace of $Y$  then, by normality, $f$ can be extended to a continuous function $F\in C(Y)$. Then, by hypothesis,  $F$ can be continuously extended to $H(U_{\rho}(Y))$. Therefore, by Lemma \ref{subspace.lemma}, $f=F|_X$ can be continuously extended to $H(U_{\rho|_X}(X))$. Hence, $\upsilon X=H(U_{\rho|_X}(X))$. Moreover, since $\varphi$ is uniformly continuous then $H(U_{\rho|_X}(X))\leq H(U_{d}(X))$. Thus, $\upsilon X=H(U_{d}(X)).$
\end{proof}
\end{lemma}

\begin{lemma}\label{partitions} Let $(X,d)$ be a Bourbaki-complete metric space, $\mathcal{P}_{n}$ be the family of all the chainable components induced by the cover of open balls $\mathcal{B} _{1/n}=\{B_{d}(x,1/n):x\in X\}$, $n\in \Nset$,  and consider the embedding $$\varphi:(X,d)\rightarrow (\prod_{n\in \Nset}\kappa _n \times \Rset ^{\omega _0}, \rho+t)$$ from Theorem \ref{embedding.metric}. Then, there is some $n_0\in \Nset$ satisfying that for every $P\in \mathcal{P}_{n_0}$ and for every $n> n_0$, the subfamily of chainable components $\{Q\in \mathcal{P}_n:Q\subset P\}$ does not have Ulam-measurable cardinal if and only if for some $j_0\in \Nset$, $\kappa _n$ is not Ulam-measurable  for every $n>j_0$. 
\begin{proof} Recall that for a uniform space $(X,\mu)$, the family of uniform partitions of $(X,s_f \mu)$ coincides with the family of uniform partitions of $(X,\mu)$. In addition, recall that from Theorem \ref{Bourbaki-complete.strongly-metrizable} there exists a complete sequence $\langle{ \mathcal{W}_n}\rangle _{n\in \Nset}$ of star-finite covers $\mathcal{W}_{n}\in \mu _d$ such that $\bigcup_{n\in \Nset} \mathcal{W}_n$ is a base for the topology of $X$. From the proof of Theorem \ref{metrizable.star-finite} and Theorem \ref{Bourbaki-complete.strongly-metrizable}, it is easy to realize that for every $n\in \Nset$, the family $\mathcal{Q}_n$ of chainable  components induced by $\mathcal{W}_n$ refines $\mathcal{P}_n$. Moreover, since every $\mathcal{W}_n$ is a uniform cover for the metric uniformity $\mu _d$, for every $n\in \Nset$ there exists some $j\geq n$ such that $\mathcal{P}_j$ refines $\mathcal{Q}_n$. Thus, it is clear that there is some $n_0\in \Nset$ satisfying that for every $P\in \mathcal{P}_{n_0}$ and for every $n> n_0$, the subfamily of chainable components $\{Q\in \mathcal{P}_n:Q\subset P\}$ does not have measurable cardinal if and only if there is some $j_0 \in \Nset$ such that for every $n>j_0$ and every $Q\in \mathcal{Q}_{j_0}$ the subfamily of chainable components $\{R\in \mathcal{Q}_n:R\subset Q\}$ does not have Ulam-measurable cardinal.

Next, if we apply the construction in the proof of  Theorem \ref{embedding.star-finite1} to the sequence $\langle \mathcal{W}_n\rangle$ we obtain the desired embedding $\varphi$ (Theorem \ref{embedding.metric}) and we can deduce the result.
\end{proof}
\end{lemma}

\begin{theorem} \label{positive} Let $(X,d)$ be a Bourbaki metric space and $\mathcal{P}_{n}$ be the family of all the chainable components induced by the cover of the open balls $\mathcal{B} _{1/n}=\{B_{d}(x,1/n):x\in X\}$. Suppose that for some $n_0\in \Nset$, for every $P\in \mathcal{P}_{n_0}$ and for every  $n> n_0$, the subfamily of chainable components $\{Q\in \mathcal{P}_n:Q\subset P\}$ does not have Ulam-measurable cardinal. Then $\upsilon X$ and $H(U_{d}(X))$ are equivalent realcompactications.
\begin{proof} By the above Lemma \ref{partitions} there exists and embedding  $$\varphi:(X,d)\rightarrow (\prod_{n\in \Nset}\kappa _n \times \Rset ^{\omega _0}, \rho+t)$$ where each cardinal $\kappa _n$ is endowed with the discrete uniformity, $\varphi$ is uniformly continuous, $\varphi (X)$ is a closed subspace of  $\prod_{n\in \Nset}\kappa _n \times \Rset ^{\omega _0}$, satisfying the additional property that for some $j\in \Nset$, $\kappa _n$ is not Ulam-measurable for every $n>j$. 

Now, we apply  Theorem \ref{Samuel.product} twice to the metric space $\Big((\prod_{n\in \Nset}\kappa _n )\times \Rset ^{\omega _0}, \rho +t\Big)$. Thus, $$H\Big(U_{\rho+t}(\prod_{n\in \Nset}\kappa _n \times \Rset ^{\omega _0})\Big)= H\Big(U_{\rho}(\prod_{n\in \Nset}\kappa _n)\Big) \times \Rset ^{\omega _0}=$$ $$ H\Big(U_{\chi}(\prod_{n=1}^j \kappa _n) \Big)\times \prod_{n>j} \kappa _n \times \Rset ^{\omega _0}.$$ The last equality is possible since it is clear that the space $(\prod_{n\in \Nset}\kappa _n,\rho)$ is uniformly homeomorphic to the product space $(\prod_{n=1} ^j \kappa _n \times \prod_{n>j}\kappa _n, \chi +\rho)$ where $(\prod_{n=1} ^j \kappa _n , \chi)$ is of course a uniformly discrete space. Therefore, by Corollary \ref{product.Samuel.discrete2}, $$H\Big(U_{\rho+t}(\prod_{n\in \Nset}\kappa _n \times \Rset ^{\omega _0})\Big)= \upsilon(\prod_{n=1}^j \kappa _n) \times \Big(\prod_{n>j} \kappa _n \times \Rset ^{\omega _0}\Big)=$$ $$ \upsilon\Big(\prod_{n=1}^j \kappa _n \times \prod_{n>j} \kappa _n \times \Rset ^{\omega _0}\Big)=\upsilon\Big(\prod_{n\in \Nset}\kappa _n \times  \Rset ^{\omega _0}\Big).$$ 

Finally, by Lemma \ref{mapping.Samuel} the result follows.
\end{proof}
\end{theorem}

\medskip

\subsection{The $G_{\delta}$-closure of a uniform space in its Samuel compactification}

\hspace{15pt} By Theorem \ref{measurable2} it is clear that if $(X,\mu)$ is Bourbaki-complete and no uniform partition has Ulam-measurable cardinal, then the Samuel realcompactification and the Hewitt realcompactification of $(X,\mu)$ are equivalent by Samuel realcompactness. Therefore, in order to find a counterexample we need to work with Ulam-measurable cardinals.

\begin{example} \label{example.betaD} {\it There exists a Bourbaki-complete uniform space such that its Samuel realcompactification is not equivalent to its Hewitt realcompactification}.
\begin{proof}[Construction] Let $D$ be a discrete space of Ulam-measurable cardinal and let $X=D\times \beta D$. Endow $X$ with product uniformity $\pi$ of the discrete metric uniformity on $D$ and the unique uniformity on $\beta D$. Then $\upsilon (D\times \beta D)$ is not equivalent to $\upsilon D \times \beta D$ by  the Theorem \ref{product.husek}. However, $H(U_{\pi}(D\times \beta D))=\upsilon D \times \beta D$ by Theorem \ref{Samuel.product}. Thus, both realcompactifications are clearly not equivalent.
\end{proof}
\end{example}

The above example does not close the problem. It is just a restriction. Indeed, even if $\upsilon (D\times\beta D)$ and $H(U_\pi (D\times \beta D))$ are not equivalent compactifications, we are going to see that $H(U_\pi (D\times \beta D))$ is equivalent to another realcompactification of $(D\times \beta D, \pi)$, namely the $G_{\delta}$-closure of $D\times \beta D$ in its Samuel compactification $s_{\pi}(D\times \beta D)$.
\smallskip

Consider the below generalization of closure.

\begin{definition}\label{def.Gdelta}Let $X\subset Y$ two spaces, then $X$ is $G_{\delta}${\it-closed} in $Y$ if for every $y\in Y\backslash X$ there exists a $G_{\delta}$-set $G$ on $Y$ such that $y\in G\subset Y\backslash X$. Thus, we say that $\overline{X}^{G_{\delta}}$ is the $G_{\delta}${\it-closure} of $X$ in $Y$, if $\overline{X}^{G_{\delta}}$ is $G_{\delta}$-closed in $Y$ and every $G_{\delta}$-set on $Y$ which meets $\overline{X}^{G_{\delta}}$ also meets $X$. Moreover, $X$ is $G_{\delta}${\it-dense} in $Y$ if $Y=\overline{X}^{G_{\delta}}$.
\end{definition}

It is clear that it is possible to define the notion of $\mathcal{S}$-closed space for any family of subsets $\mathcal{S}$ of $Y$. In particular, for a (Tychonoff) space $Y$, $X$ is $G_{\delta}$-closed in $Y$ if and only if $X$ is {\it zero-closed}. Indeed, in a (Tychonoff) space every zero-set is $G_{\delta}$-set and, conversely, for every $G_{\delta}$-set $G$ and every $x\in G$ there exists a zero-set $Z$ such that $x\in Z\subset G$.
\smallskip

It is well-known that the Hewitt realcompactification $\upsilon X$ of a space $X$ is the $G_{\delta}$-closure of $X$ in its Stone-\v{C}ech-compactification $\beta X$ \cite[8.8]{gillman}. Moreover,  given $X$ a space and $\mathcal{L}\subset X$ a unital vector lattice separating points from closed sets of $X$, $H(\mathcal{L})$ is $G_{\delta}$-closed in $H(\mathcal{L}^{*})$. Indeed, if $\xi \in H(\mathcal{L}^{*})\backslash H(\mathcal{L})$, then there exists a function $f\in \mathcal{L}$ such that $f^{*}(\mathcal{L})=\infty$ (Theorem \ref{realcompact.subspace}). Then, $\xi$ belongs to the $G_{\delta}$-set $(f^{*})^{-1}(\{\infty\})$, which is disjoint from $H(\mathcal{L})$. As a consequence, the $G_{\delta}$-closure of $X$ in $H(\mathcal{L}^{*})$ is contained in $H(\mathcal{L})$. In fact, the $G_{\delta}$-closure of $X$ in $H(\mathcal{L}^{*})$  is a realcompactification of $X$ (\cite{hager1}).

\smallskip

Next, we consider the unital vector lattice $U_\mu(X)$ of a uniform space $(X,\mu)$. We denote by $ \mathcal{A}(U_\mu (X))$ the smallest uniformly closed subalgebra of $C(X)$ containing $U_\mu (X)$ which, in addition, is closed under inversion, that is, $1/f\in \mathcal{A}(U_\mu (X))$ for every $f\in \mathcal{A}(U_\mu (X))$ such that $f(x)\neq 0$ for every $x\in X$. By \cite[Theorem 2.5]{garrido.algebras}, $\mathcal{A}(U_\mu (X))$ is exactly the uniform closure of any of the following families of functions:

$$A_{0}(U_\mu (X))=\{f/g: f,g\in U_\mu (X), g(x)\neq 0 \text{ for every }x\in X\}$$

$$A_{0}(U_\mu ^*(X))=\{f/g: f,g\in U_\mu ^*(X), g(x)\neq 0 \text{ for every }x\in X\}.$$
\smallskip

Let us denote by $z_u$ the family of all zero-sets $Z$ of $X$ such that $Z=f^{-1}(0)$ for some $f\in U_\mu (X)$ (equivalently, $f\in U_\mu ^*(X)$).

\begin{definition} A $z_u$-{\it filter} of a space $X$ is a filter such that $z_u\cap \mathcal{F}$ is a filter-base of $\mathcal{F}$. A $z_u$-{\it ultrafilter} $\mathcal{F}$ is a  $z_u$-filter such that $z_u\cap \mathcal{F}$ is a maximal family in $z_u$, that is, $Z\in z_u$ belongs to $\mathcal{F}$ whenever $Z\cap F\neq \emptyset$ for every $F\in \mathcal{F}$.
\end{definition}

\begin{theorem} \label{Gdelta} {\rm (\cite{hager1},\cite{husek})} Let $(X,\mu)$ be a uniform space. The following statements are equivalent:
\begin{enumerate} 
\item $Y$ is the $G_{\delta}$-closure of $X$ in $s_\mu X$;

\item $Y$ is the intersection of all the cozero-sets in $s_\mu X$ containing $X$;

\item for every $\xi \in s_\mu X\backslash Y$ there exists some strictly positive function $f\in U_{\mu}^{*} (X)$ (or $f\in U_{\mu} (X)$) such that $f^{*}(\xi)=0$;

\item $Y$ is  the subset of all the clusters points in $H(U_{\mu}^{*} (X))$ of all the $z_u$-ultrafilters of $X$ satisfying the countable intersection property;

\item $Y$ is a realcompactification equivalent to $H(\mathcal{A}(U_{\mu} (X)))$.

\end{enumerate}



\end{theorem} 

\begin{remark}\label{rem.g-dense} The  realcompactification $H(\mathcal{A}(U_{\mu} (X)))$ is the smallest realcompactification of $X$ such that every function in $U_{\mu} (X)$ and every inverse function $1/f$ of every non-vanishing function $f$ in $U_{\mu} (X)$ can be continuously extended to it. In particular, $X$ is $G_{\delta}$-dense in $H(\mathcal{A}(U_{\mu}(X)))$.
\end{remark}

\begin{remark}By Theorem \ref{realcompact.subspace}, the realcompactification  $H(\mathcal{A}(U_{\mu} (X)))$ is a topological subspace of the compactification $\mathcal{A}^{*}(U_{\mu} (X))$ where $\mathcal{A}^{*}(U_{\mu} (X))=\mathcal{A}(U_{\mu} (X))\cap C^{*}(X)$. In general the compactification $H(\mathcal{A}^{*}(U_{\mu} (X)))$ is not equivalent to the Samuel compactification $s_\mu X$ (see \cite[Example I]{chekeev} and \cite{hager1}).
\end{remark}

\smallskip

Recall that a topological space $X$ is $z$-{\it embedded} in the space $Y$, whenever $X\subset Y$ and each zero-set of $X$ is the intersection with $X$ of a zero-set in $Y$. For instance, any metric space $(X,d)$ is clearly $z$-embedded in $s_d X$ and therefore in any set $Y$ such that $X\subset Y \subset s_d X$. Moreover a subspace $A$ of a space $X$ is $C$-{\it embedded} in $X$ if every real-valued continuous function on $A$ can be extended to a continuous function on $X$. Next, applying the result \cite[Corollary 3.6]{blair-hager} of Blair and Hager,  we have that under $G_{\delta}$-density assumption, $z$-embedding and $C$-embedding are equivalent properties. Then, the following result can be deduced.

\begin{theorem} \label{Gdelta.metric} For a metric space $(X,d)$, $ H(\mathcal{A}(U_d(X)))$ and $\upsilon X$ are equivalent realcompactifications. In particular, $\upsilon X$ is a topological subspace of $H(U_d(X))$, that is, $\upsilon X\subset H(U_d(X))$.
\begin{proof} As we have said above, $X$ is clearly $z$-embedded in $s_d X$ and also in $H(\mathcal{A}(U_{\mu}(X)))$. Now, since $X$ is $G_\delta$-dense in $H(\mathcal{A}(U_{\mu}(X)))$ (see Remark \ref{rem.g-dense}), by the above mentioned result, $X$ is $C$-embedded in $H(\mathcal{A}(U_{\mu}(X)))$. But the unique realcompactification in which $X$ is $C$-embedded is $\upsilon X$ \cite{gillman}. Therefore, $H(\mathcal{A}(U_{\mu}(X)))$ and $\upsilon X$ must be equivalent realcompactifications.

Finally, since $H(\mathcal{A}(U_{d}(X))\subset H(U_{d}(X))$, as $H(\mathcal{A}(U_{d}(X))$ is the $G_{\delta}$-closure of $X$ in $s_d X$, then $\upsilon X\subset  H(U_{d}(X))$.
\end{proof}
\end{theorem}

\begin{example} {\it There exists a complete metric space $(X,d)$ such that $H(U_d (X))$ and $H(\mathcal{A}(U_{d}(X))$ are not equivalent realcompactifications}. 
\begin{proof} Since for any metric space $\upsilon X=H(\mathcal{A}(U_d(X)))$, then any complete but not Bourbaki-complete realcompact metric space is such and example by Theorem \ref{measurable2}.  
\end{proof}
\end{example}


Motivated by Example \ref{example.betaD}, as not every Bourbaki-complete uniform space satisfies that the Hewitt realcompactification and the Samuel compactification are equivalent realcompactifications, we ask know when the realcompactifications $H(U_{\mu}(X))$ and $H(\mathcal{A}(U_{\mu}(X)))$ are equivalent. Clearly, by Theorem \ref{Gdelta.metric}, in the frame of metric spaces, this is the same than asking if for every Bourbaki-complete metric space, $H(U_{\mu}(X))$ and $\upsilon X$ are equivalent realcompactifications. Observe that,  by Theorem \ref{Gdelta}, this is also equivalent to ask if for every non-vanishing function $f\in U_{\mu}(X)$ the inverse $1/f$ can be continuously extended to $H(U_{\mu}(X))$. Thus, solving this questions,  should allow us to know a little bit more the subalgebra of functions $C(H(U_{\mu}(X)))$.
\smallskip



\begin{lemma}\label{lemma2f} Let $(X,\mu)$ be a uniform space such that $X=H(\mathcal{A}(U_{\mu}(X))$ (or in particular, Samuel realcompact) and let $(Y,\nu)$ be any uniform space. Then every uniformly continuous function $f\in U_{\mu \times \nu}(X\times Y)$ satisfying that $f(x,y)\neq 0$ for every $(x,y)\in X\times Y$ can be extended to a continuous function $f^{*}$ on $X\times H(\mathcal{A}(U_{\nu}(Y)))=H(\mathcal{A}(U_{\mu}(X)))\times H(\mathcal{A}(U_{\nu}(Y)))$ satisfying also that $f^{*}(x,\xi)\neq 0$ for every $(x,\xi)\in X\times H(\mathcal{A}(U_{\nu}(Y))$.
\begin{proof} Let $f\in U_{\mu\times\nu}(X\times Y)$ satisfying the additional property that $f(x,y)\neq 0$ for every $(x,y)\in X\times Y$. Then every function $f_x:(Y,\nu)\rightarrow (\Rset, d_u)$ defined by $f_x(y)=f(x,y)$,  is uniformly continuous and satisfies also that $f_x(y)\neq 0$ for every $y\in Y$. By Theorem \ref{Gdelta} every function $f_x$ can be extended to a continuous function $f^{*}_x :H(\mathcal{A}(U_{\mu}(X)))\rightarrow \Rset$ satisfying that $f^{*}_x(\xi)\neq 0$ for every $\xi\in H(U_{\nu}(Y))$.

Like in the proof of Lemma \ref{lemma1f}, the function $f^{*}:X\times H(\mathcal{A}(U_{\nu}(Y))) \rightarrow \Rset$ defined by $f^{*}(x,\xi)=f^{*}_x(\xi)$ is continuous. In addition, it is trivial that $f^{*}(x,\xi)\neq 0$ for every $(x,\xi)\in X\times H(\mathcal{A}(U_{\nu}(Y)))$.
\end{proof}
\end{lemma}

\begin{theorem} \label{product.Gdelta} Let $(X,\mu)$ Samuel realcompact space and $(Y,\nu)$ a uniform space satisfying that $H(U_{\nu}(Y))$ and $\upsilon Y$ are equivalent realcompactifications. Then all the following realcompactifications of $X\times Y$ are equivalent:

\begin{enumerate}

\item $H(\mathcal{A}(U_{\mu\times \nu}(X\times Y)))$

\item $H((U_{\mu\times \nu}(X\times Y))$

\item $\upsilon X\times \upsilon Y$

\item $H(\mathcal{A}(U_{\mu}(X)))\times H(\mathcal{A}(U_{\nu}(Y)))$

\item $H(U_{\mu}(X))\times H(U_{\nu}(Y))$

\item $X\times \upsilon Y$

\item $X\times H(\mathcal{A}(U_{\nu}(Y)))$

\item $X\times H(U_{\nu}(Y))$

\end{enumerate}

\begin{proof}

By Theorem \ref{Samuel.product}, we have that  $$H(U_{\mu\times \nu}(X\times Y))=H(U_{\mu}(X))\times H(U_{\nu}(Y))= X\times H(U_{\nu}(Y)).$$ Moreover, since $(X,\mu)$ is Samuel realcompact and $(Y,\nu)$ satisfies that $\upsilon Y=H(U_{\nu}(Y))$ then $$H(\mathcal{A}(U_{\mu}(X)))\times H(\mathcal{A}(U_{\nu}(Y))=X\times H(\mathcal{A}(U_{\nu}(Y)))=X\times H(U_{\nu}(Y))=$$ $$X\times \upsilon Y=\upsilon X\times \upsilon Y.$$ By Lemma \ref{lemma2f} $$H(\mathcal{A}(U_{\mu}(X)))\times H(\mathcal{A}(U_{\nu}(Y))) \geq  H(\mathcal{A}(U_{\mu\times \nu}(X\times Y))).$$ On the other hand, $H(\mathcal{A}(U_{\mu\times \nu}(X\times Y)))\geq H(U_{\mu\times \nu}(X\times Y))$ in general. Therefore, by the above equalities,

$$H(\mathcal{A}(U_{\mu\times \nu}(X\times Y)))= H(U_{\mu\times \nu}(X\times Y)).$$

\end{proof}

\end{theorem}

By the above result, it is clear that, for Example \ref{example.betaD}, $H(U_{\pi}(D \times \beta D)=\upsilon D \times \beta D= H(\mathcal{A}(U_{\pi}(D\times \beta D))$, where $(D\times \beta D, \pi)$ is of course a Bourbaki-complete space. However, we don't have a general answer for Bourbaki-complete uniform spaces parallel to Theorem \ref{positive}.

\begin{remark}
It is clear that in order to characterize those Bourbaki-complete metric spaces $(X,d)$ (or Bourbaki-complete uniform spaces $(X,\mu)$) satisfying that $\upsilon X$ and $H(U_{d}(X))$ are equivalent realcompactifications (or, for uniform spaces, $H(\mathcal{A}(U_{\mu}(X))$ and $H(U_\mu(X))$ are equivalent realcompactifications), we need a counterexample to the question. We strongly believe that this counterxample will be the uniformly 0-dimensional $(D^{\omega _0},\rho)$ where $D$ is a uniformly discrete space whose cardinal is not only Ulam-measurable cardinal but also $\omega _1$-strongly compact (see \cite{bagaria1} and \cite{bagaria2} for definitions). 
With the help of this couterexample, and Theorem \ref{positive} we would like first to characterize those complete uniformly 0-dimensional space satisfying our question, and then deduce from it the general case for all the Bourbaki-complete uniform spaces, with the help of the embeddings from the second part of this thesis.
\end{remark} 

\medskip

\medskip

\medskip

\bigskip

\bigskip

\bigskip

\begin{center} \Denarius \Denarius  \Denarius \end{center}
\bigskip

\bigskip

\newpage

\newpage
\mbox{}
\thispagestyle{empty}

\end{document}